\documentclass[preprint,12pt]{elsarticle}

\usepackage{mathrsfs}
\usepackage{amsmath}
\usepackage{amssymb}
\usepackage{amsthm}
\usepackage{epsfig}
\usepackage{graphics}
\usepackage{graphicx}
\usepackage{txfonts}
\usepackage{amsfonts}
\usepackage{multirow,tabularx}
\usepackage{caption}
\usepackage{float}
\usepackage{comment}
\usepackage{latexsym}
\usepackage{indentfirst}
\usepackage{url}
\usepackage{enumitem}
\usepackage{wrapfig}
\newcommand\RR{\mathbb{R}}
\newcommand\Rd{\mathbb{R}^d}
\newcommand\NN{\mathbb{N}}

\newcommand\PP{\mathbb{P}}
\newcommand\Mean{\text{E}}

\newcommand\vx{\boldsymbol{x}}   % x in R^d
\newcommand\vy{\boldsymbol{y}}   % y in R^d
\newcommand\vz{\boldsymbol{z}}
\newcommand\vc{\boldsymbol{c}}
\newcommand\vb{\boldsymbol{b}}

\newcommand\vv{\boldsymbol{v}}

\newcommand\vzero{\boldsymbol{0}}
\newcommand\vk{\boldsymbol{k}}
\newcommand\valpha{\boldsymbol{\alpha}}
\newcommand\vbeta{\boldsymbol{\beta}}

\newcommand\ve{\boldsymbol{e}}
\newcommand\vf{\boldsymbol{f}}
\newcommand\vt{\boldsymbol{t}}

\newcommand\vvartheta{\boldsymbol{\vartheta}}

\newcommand\vP{\boldsymbol{\mathcal{P}}}
\newcommand\vB{\boldsymbol{\mathcal{B}}}
\newcommand\vL{\boldsymbol{\mathcal{L}}}
\newcommand\vK{\mathsf{K}}
\newcommand\vA{\mathsf{A}}
\newcommand\vV{\mathsf{V}}
\newcommand\vD{\mathsf{D}}
\newcommand\vTheta{\mathsf{\Theta}}
\newcommand\vS{\boldsymbol{S}}

\newcommand\ud{\textup{d}}
\newcommand\Normal{\mathcal{N}}
\newcommand\Order{\mathcal{O}}
\newcommand\erfc{\textup{erfc}}

\newcommand\Span{\text{span}}
\newcommand\Diag{\text{diag}}
\newcommand\Filter{\mathcal{F}}
\newcommand\Borel{\mathcal{B}}
\newcommand\Domain{\mathcal{D}}
\newcommand\Aset{\mathcal{A}}
\newcommand\Eset{\mathcal{E}}
\newcommand\Var{\text{Var}}
\newcommand\Cov{\text{Cov}}

\newcommand\Power{\mathcal{P}}
\newcommand\measure{\mathit{m}}

\newcommand\Quadratic{\mathcal{Q}}

\newcommand\Phistar{\Phi^{*}}

\newcommand\Hilbert{\mathcal{H}}
\newcommand\Cont{\mathrm{C}}
\newcommand\Leb{\mathrm{L}}

\newcommand\Gaussian{\mathcal{G}}

\newcommand{\norm}[1]{\left\lVert#1\right\rVert}
\newcommand{\abs}[1]{\left\lvert#1\right\rvert}

\newtheorem{theorem}{Theorem}[section]
\newtheorem{lemma}[theorem]{Lemma}
\newtheorem{corollary}[theorem]{Corollary}
\newtheorem{proposition}[theorem]{Proposition}

\theoremstyle{definition}
\newtheorem{definition}[theorem]{Definition}
\newtheorem{example}[theorem]{Example}

\theoremstyle{remark}
\newtheorem{remark}[theorem]{Remark}

\numberwithin{equation}{section}
\numberwithin{figure}{section}
\numberwithin{table}{section}
\journal{ArXiVview}

\begin{document}

\begin{frontmatter}



\title{Kernel-based Methods for \\ Stochastic Partial Differential Equations}


\author{Qi Ye}\ead{qiye@syr.edu}

\address{Mathematics Department, Syracuse University, Syracuse, NY 13244}

\begin{abstract}
This article gives a new insight of kernel-based (approximation) methods to solve the high-dimensional stochastic partial differential equations.
We will combine the techniques of meshfree approximation and kriging interpolation to extend the kernel-based methods for the deterministic data to the stochastic data.
The main idea is to endow the Sobolev spaces with the probability measures induced by the positive definite kernels such that
the Gaussian random variables can be well-defined on the Sobolev spaces.
The constructions of these Gaussian random variables provide
the kernel-based approximate solutions of the stochastic models.
In the numerical examples of the stochastic Poisson and heat equations, we show that the approximate probability distributions are well-posed for various kinds of kernels such as the compactly supported kernels (Wendland functions) and the Sobolev-spline kernels (Mat\'ern functions).
\end{abstract}


\end{frontmatter}


\section{Introduction}\label{sec:Intr}
%---------------------------------------------------------------------------------------------------------------------

In this article, we will study with the approximate solutions of the stochastic partial differential equations (SPDEs) by the kernel-based (approximation) methods.
The SPDEs frequently arise from applications in areas such as physics, biology, engineering, economics, and finance.
Many analytical theorems of stochastic differential equations have been developed in~\cite{Oksendal2003,Chow2007} and their numerical algorithms are aslo a fast growing research area in~\cite{KloedenPlaten1992,BabuvskaTemponeZouraris2004,BabuvskaNobileTempone2010,JentzenKloeden2011,KuoSchwabSloan2012}.
Unfortunately, the current numerical tools often show the limited success in the high-dimensional equations or the complicated boundary conditions.

Recently, the kernel-based methods become a fundamental approach for scattered data approximation, statistical (machine) learning, engineering design, and numerical solutions of partial differential equations.
In particular, the research areas of the kernel-based methods cover the interdisciplinary fields of approximation theory and statistical learning such as meshfree approximation in~\cite{Buhmann2003,Wendland2005,Fasshauer2007} and kriging interpolation in~\cite{Stein1999,BerlinetThomas2004,HastieTibshiraniFriedman2009}.
Moreover, the kernel-based methods are known
by a variety of names in the monographs including radial basis functions, kernel-based collocation, smoothing splines, and Gaussian process regression.

In the studies of approximation theory and statistical learning, we know that
the kernel-based methods can be used to approximate the high-dimensional partial differential equations and estimate the simple stochastic models.
Naturally, therefore, we develop the kernel-based methods to approximate the high-dimensional SPDEs in
the paper~\cite{CialencoFasshauerYe2012} and the doctoral thesis~\cite{YeThesis2012}.
Now we propose to improve and complete the theorems and algorithms of the kernel-based methods for the high-dimensional stochastic data. The main idea is to combine the knowledge of approximation theory, statistical learning, probability theory, and stochastic analysis into one theoretical structure.
In this article, we will mainly focus on the mixture techniques of meshfree approximation and kriging interpolation for the constructions of the kernel-based approximate solutions of the stochastic models such that the kernel-based estimators have the both globally and locally geometrical meanings.
{
\center
\includegraphics[width=\textwidth, height=0.33\textwidth]{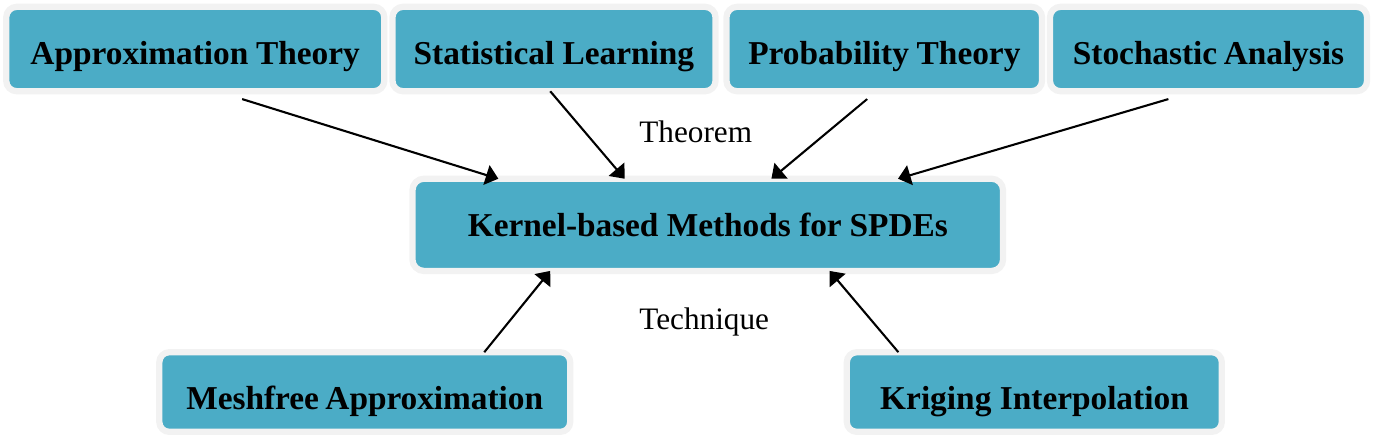}
}

Now we give the outlines of this article.
In Section~\ref{sec:Initial}, we firstly describe the initial ideas of the new insights of the kernel-based methods.
In the beginning of our researches, we study with the meshfree approximation for the high-dimensional interpolation by the positive definite kernels, for example,
the kernel-based interpolant induced by the Gaussian kernels in Figure~\ref{fig:1DExaKernel}. The reproducing properties also guarantee that the kernel-based interpolants are the globally optimal recovery in the reproducing kernel Hilbert spaces.
Here, we have a new idea to obtain the locally best estimators based on all globally interpolating paths by the statistics \& probability techniques
(see a simple example in Figure~\ref{fig:1DExaInitial} and Table~\ref{tab:InitialIdea}).
This indicates that we need some probability structures to measure the interpolating paths.
Moreover, we find that the kriging interpolation also provides the locally best linear unbiased prediction by the Gaussian fields such as the 1D example in Figure~\ref{fig:1DExaKriging}.
The recent paper~\cite{ScheuererSchabackSchlather2012} shows that the formulations of the both kernel-based interpolants and simple kriging predictions are the same. Thus, we guess that the meshfree approximation and the kriging interpolation could be strongly connected by one theoretical approach such that the global and local approximations could be obtained at the same time.
By the theorems of stochastic analysis, we know that the Brownian motion can be constructed on the continuous function space endowed with the Wiener measure (see \cite[Chapter~2]{KaratzasShreve1991} or Section~\ref{sec:Initial}). Then the Wiener measure and the Brownian motion provide a tool to measure the continuous interpolating paths.
It is also well-known that the Brownian motion is a Gaussian field and its covariance kernel is a min kernel which is a positive definite kernel.
Therefore, we will combine the knowledge of the kernel-based interpolants, the simple kriging predictions, and the Brownian motions together to renew the kernel-based methods.
More precisely, we will use the positive definite kernels to introduce the probability measures and the Gaussian fields on the Sobolev spaces such that the initial ideas can be generalized to measure all smooth interpolating paths.
Then the kernel-based probability structures of the Sobolev spaces will help us to obtain the kernel-based approximation for the deterministic and stochastic data.

In Section~\ref{sec:Gauss-PDK},
we will extend the initial ideas of a simple example of the 1D interpolating paths in Figure~\ref{fig:1DExaInitial} to all interpolating paths in the Sobolev spaces in Figure~\ref{fig:Sobolev-initial}.
Firstly, we will construct the Gaussian random variables by the chosen positive definite kernel $K$.
In this article, the \emph{Gaussian random variables} include Gaussian fields and normal random variables.
In Theorem~\ref{t:Gauss-PDK-L}, for any bounded linear functional $L$ such as $L:=\delta_{\vx}$ or $L:=\delta_{\vx}\circ\Delta$, the normal random variable $LS(\omega):=L\omega$ is well-defined on the Sobolev space $\Hilbert^m(\Domain)$ endowed with the probability measure $\PP_K$ induced by the kernel $K$.
Next, we will rethink the kernel-based approximation for the deterministic data by the kernel-based probability structures of the Sobolev spaces.
Then the constructions of the multivariate normal random variables $L_1S,\ldots,L_NS$ indicate a connection of kernel-based approximation to a kernel basis $L_{1,\vy}K(\cdot,\vy),\ldots,L_{N,\vy}K(\cdot,\vy)$ such as the kernel-based approximate functions are the linear combinations of the kernel basis (see Equation~(\ref{eq:deterministic-max-cond-prob}-\ref{eq:kernel-approx-model-2})).
Combining with the maximum likelihood estimation methods, we can obtain the locally optimal estimators (kernel-based estimators) which are also supported by the globally interpolating paths for the given data.

In Sections~\ref{sec:Stochastic-Data-Approx}-\ref{sec:par-SPDE},
we will extend the kernel-based methods for the deterministic problems in \cite{Wendland2005,Fasshauer2007,HonSchabackZhong2013} to the stochastic problems
such as the stochastic data interpolations, the elliptic SPDEs, and the parabolic SPDEs.
Same as in \cite{Wendland2005,Fasshauer2007,HonSchabackZhong2013},
we can also analyze their convergence in probability by the power functions and the fill distances.
Moreover, Section~\ref{sec:num-exp} shows the 3D, 2D, and 1D numerical examples of the stochastic data interpolations and the stochastic Poisson and heat equations by various kinds of positive definite kernels such as the Gaussian kernels, the compactly supported kernels (Wendland functions), and the Sobolev-spline kernels (Mat\'ern functions).
For reducing the complexity, we only look at the linear stochastic models here.
In fact, there are still many tools of statistical learning to solve the nonlinear stochastic models such as support vector machines with various loss functions in \cite{CuckerSmale2002,SteinwartChristmann2008}.
In Section~\ref{sec:Final}, we briefly describe the improvements and advance researches of the theorems and algorithms discussed in this article.

%---------------------------------------------------------------------------------------------------------------------
\section{Initial Ideas}\label{sec:Initial}
%---------------------------------------------------------------------------------------------------------------------

The approximation theory focuses on how a function $u$ can be approximated by a well-computable function $\hat{u}$.
Typically, the fundamental problem can be represented as follows. We have the data values $f_1,\ldots,f_N$ sampled from the function $u:[0,1]\to\RR$ at the distinct data points $X:=\left\{x_1,\ldots,x_N\right\}\subseteq[0,1]$, that is, $f_1:=u(x_1),\ldots,f_N:=u(x_N)$. An approximate function $\hat{u}:[0,1]\to\RR$ will be constructed to interpolate the given data values $f_j$ at $x_j$, that is, $\hat{u}(x_1)=f_1,\ldots,\hat{u}(x_N)=f_N$.
So, we can use this interpolant $\hat{u}$ to estimate $u$ at any unknown location $z\in[0,1]$, that is, $u(z)\approx \hat{u}(z)$.

\begin{figure}[h]
\center
\includegraphics[width=0.68\textwidth, height=0.38\textwidth]{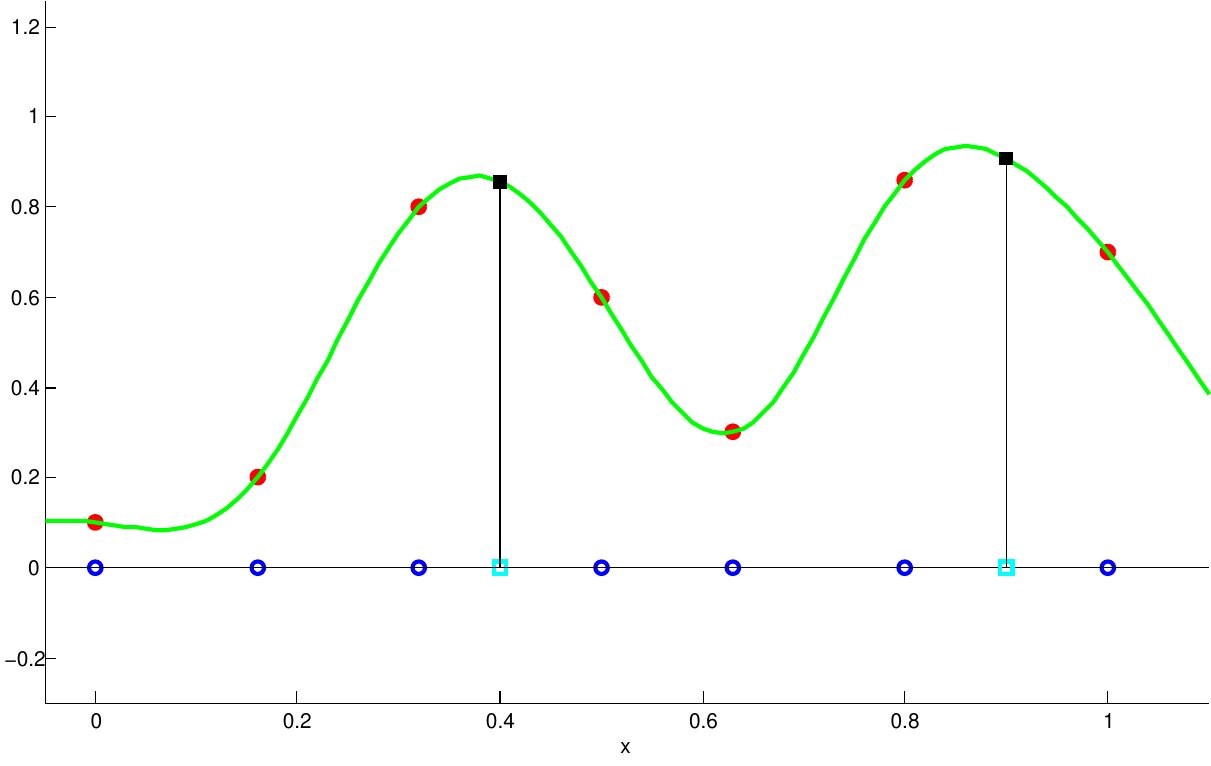}
\caption{\small The 1D example of meshfree approximation: The blue circles represent the data points $x_j$ and the red circles represent the data values $f_j$ at $x_j$ for $j=1,\ldots,7$. The green curve represents the kernel-based interpolant $\hat{u}$ induced by the Gaussian kernel with the shape parameter $\theta=6$. The black squares represent the estimate values $\hat{u}(z_1)$ and $\hat{u}(z_2)$ at the unknown locations $z_1:=0.4$ and $z_2:=0.9$ shown in cyan.}\label{fig:1DExaKernel}
\end{figure}

By the classical methods of polynomial and spline interpolation, an interpolant $\hat{u}$ will be constructed by the polynomials or the spline functions in \cite[Chapter~6]{KincaidWard2002}.
Recently, the kernel-based methods (radial basis functions) give a novel approximation tool to construct the \emph{kernel-based interpolant} $\hat{u}$ by a positive definite kernel $K:[0,1]\times[0,1]\to\RR$ (see Definition~\ref{d:PDK}), for example, the Gaussian kernel with the shape parameter $\theta>0$
\[
K(x,y):=e^{-\theta^2\abs{x-y}_2^2},\quad\text{for }x,y\in[0,1].
\]
To be more precise, the kernel-based interpolant $\hat{u}$ is a linear combination of the kernel basis $K(\cdot,x_1),\ldots,K(\cdot,x_N)$ such as
\[
\hat{u}(z):=\sum_{k=1}^Nc_kK(z,x_k),\quad\text{for }z\in[0,1],
\]
and the coefficients $\vc:=\left(c_1,\cdots,c_N\right)^T$ are computed by a well-posed linear system
\[
\vK_X\vc=\vf,
\]
where $\vK_X:=\left(K(x_j,x_k)\right)_{j,k=1}^{N,N}$ and $\vf:=\left(f_1,\cdots,f_N\right)^T$.
More details of kernel-based interpolation or called meshfree approximation are mentioned in the books \cite{Wendland2005,Fasshauer2007}.
Figure~\ref{fig:1DExaKernel} illustrates an example of the kernel-based interpolant $\hat{u}$ induced by the Gaussian kernels which is also the minimizer over the reproducing norms globally.

\begin{figure}[h]
\center
\includegraphics[width=0.68\textwidth, height=0.38\textwidth]{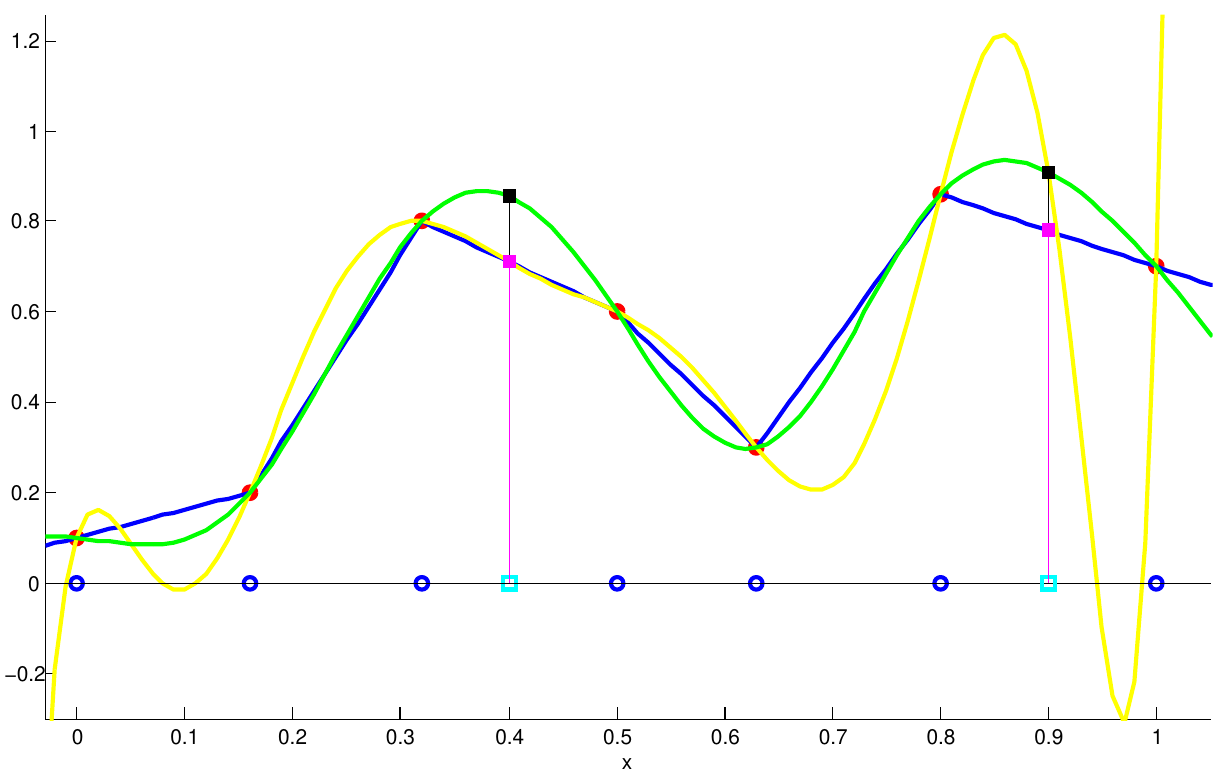}
\caption{\small The 1D example of the initial ideas:
The given data $f_j$ and $x_j$ (red and blue circles) are the same as in Figure~\ref{fig:1DExaKernel}.
The blue, green, and yellow curves represent the piecewise linear spline $\hat{u}_1$, the kernel-based interpolant $\hat{u}_2$, and the polynomial interpolant $\hat{u}_3$, respectively. The black and pink squares represent the different choices of the estimate values $\hat{u}_k(z_i)$ at the cyan squares $z_i$ for $k=1,2,3$ and $i=1,2$.}\label{fig:1DExaInitial}
\end{figure}

Usually, there are many choices of the interpolants $\hat{u}$ to approximate the unknown values $u(z)$.
So, we need to determine which estimator $\hat{u}(z)$ is the best.
Different from the classical approximation theory, we will choose the best estimator $\hat{v}$ based on all feasible interpolating paths by the statistics \& probability techniques.
Let us look at a simple example in Figure~\ref{fig:1DExaInitial} and Table~\ref{tab:InitialIdea} to study with the initial ideas of this article.
Figure~\ref{fig:1DExaInitial} has three interpolating paths, that is, the piecewise linear spline $\hat{u}_1$ (the blue curve), the kernel-based interpolant $\hat{u}_2$ (the green curve), and the polynomial interpolant $\hat{u}_3$ (the yellow curve).
We observe that there are two choices of the estimated values at $z_1,z_2$ given by $\hat{u}_1,\hat{u}_2,\hat{u}_3$ (see the black and pink squares in Figure~\ref{fig:1DExaInitial}).
Here, we view the interpolating paths $\hat{u}_1,\hat{u}_2,\hat{u}_3$ as the sample events. Then the happenings of the black and pink squares are supported by the sample events $\hat{u}_j$.
More precisely, the probabilities of the black and pink squares are counted by the numbers of the interpolating paths, for example, the probability of the black square at $z_2$ is endowed with $2/3$ because the both interpolating paths $\hat{u}_2$ and $\hat{u}_3$ pass it.
Naturally, we will choose the best estimators $\hat{v}_1$ and $\hat{v}_2$ to approximate $u(z_1)$ and $u(z_2)$, respectively, by the maximal probabilities in Table~\ref{tab:InitialIdea}.

\begin{table}[h]
\center
\begin{tabular}{c c c c}
\hline\hline
 Locations & Probabilities at Black & Probabilities at Pink & Best Estimators \\
\hline
 $z_1=0.4$ & $1/3$ (counted by $\hat{u}_2$) & $2/3$ (counted by $\hat{u}_1,\hat{u}_3$) & $\hat{v}_1:=\hat{u}_1(z_1)$ \\
\hline
 $z_2=0.9$ & $2/3$ (counted by $\hat{u}_2,\hat{u}_3$) & $1/3$ (counted by $\hat{u}_1$) & $\hat{v}_2:=\hat{u}_2(z_2)$ \\
\hline\hline
\end{tabular}
\caption{\small The initial ideas of the best estimators based on Figure~\ref{fig:1DExaInitial}.
}\label{tab:InitialIdea}
\end{table}

In Figure~\ref{fig:1DExaInitial}, we observe that the polynomial interpolating path $\hat{u}_3$ passes the both best estimators $\hat{v}_1$ and $\hat{v}_2$ at $z_1$ and $z_2$. By the classical methods, the interpolant $\hat{u}_3$ is not a good approximation which indicates that the estimators $\hat{v}_1$ and $\hat{v}_2$ can not be obtained at the same time.
But, the ill-posed problem of $\hat{u}_3$ just occurs globally and the best estimator can be obtained by $\hat{u}_3$ at some local points.
The extreme example in Figure~\ref{fig:1DExaInitial} let us rethink the classical approximation problems to connect the global interpolants and the local optimizers.
Generally speaking, we will look at all feasible interpolating paths and the best estimator is dependent of the largest probability counted by the interpolating paths massing at the unknown locations.
This indicates that we need a probability structure of the interpolating paths to measure various estimate values.

\begin{figure}[h]
\center
\includegraphics[width=0.68\textwidth, height=0.38\textwidth]{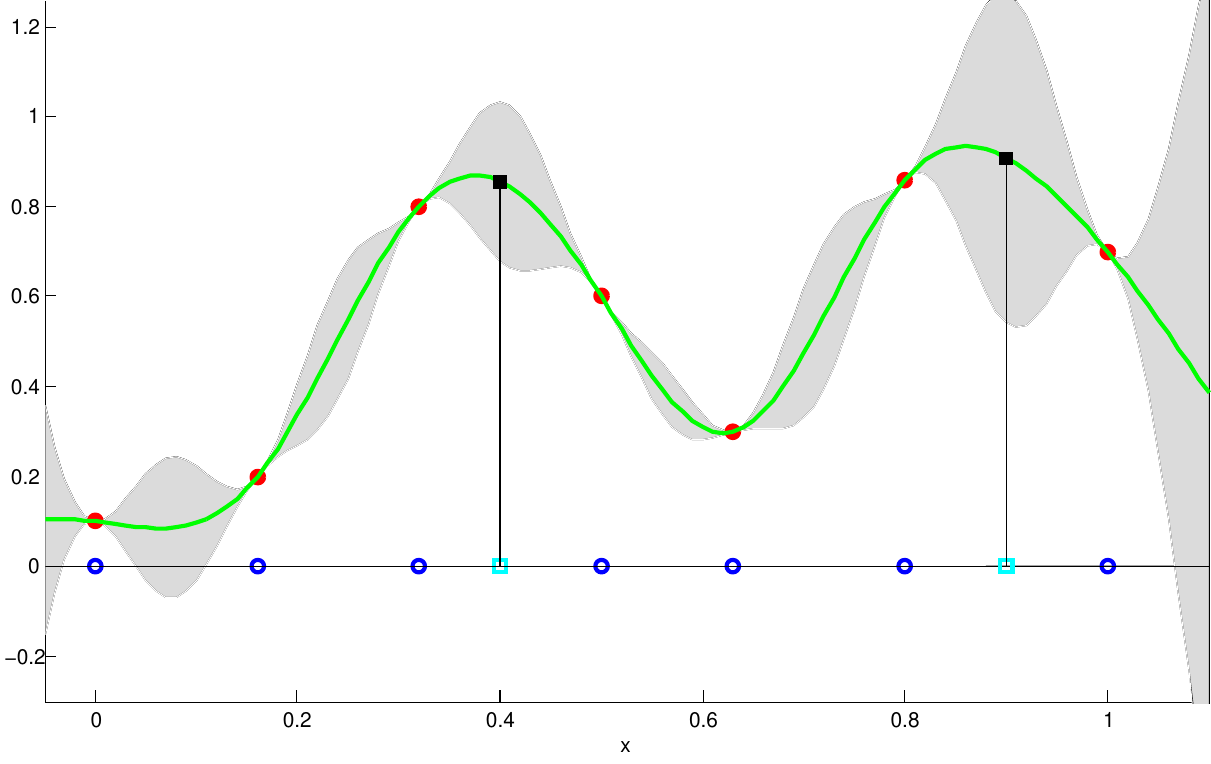}
\caption{\small The 1D example of kriging interpolation:
The observed data $f_j$ and $x_j$ (red and blue circles) are the same as in Figure~\ref{fig:1DExaKernel}.
The covariance kernel of the Gaussian field $S$ is the Gaussian kernel with the shape parameter $\theta=6$.
The best linear unbiased predictions $\hat{s}(z)$ of $S_{z}$ conditioned on the observed data values $f_j$ at $x_j$ are shown in the green curve
which runs along the means of the normally distributed confidence intervals of $99\%$ shown in gray.
}\label{fig:1DExaKriging}
\end{figure}

Moreover, we find that the kriging interpolation provides another way to obtain the locally best estimators by the Gaussian fields.
In statistical learning, originally in geostatistics, the kriging interpolation in~\cite{Stein1999} is modeled by a Gaussian field $S$ with a prior covariance kernel $K$. As Definition~\ref{d:Gaussian}, the Gaussian field $S$ composes of the deterministic domain $[0,1]$ and the probability space $\left(\Omega,\Filter,\PP\right)$, where $\PP$ is a probability measure defined on a measurable space $\left(\Omega,\Filter\right)$.
Roughly, the Gaussian field $S$ can be viewed as a map from $[0,1]\times\Omega$ into $\RR$.
Then $S_x$ is a normal random variable on the probability space $\left(\Omega,\Filter\right)$ for any $x\in[0,1]$.
For convenience, we suppose that the Gaussian field $S$ has the mean $0$ and the covariance kernel which is a Gaussian kernel.
In kriging interpolation, the data values $f_1,\ldots,f_N$ are viewed as the realized observations of the normal random variables $S_{x_1},\ldots,S_{x_N}$.
By the simple kriging methods, we can obtain the \emph{best linear unbiased prediction} $\hat{s}(z)$ of the Gaussian field $S$ at any unobserved location $z$ conditioned on the observed data values $f_j$ at $x_j$, that is,
\[
\hat{s}(z):=\Mean\left(S_{z}|S_{x_1}=f_1,\ldots,S_{x_N}=f_N\right)=\vk_X(z)\vK_X^{-1}\vf,\quad\text{for }z\in[0,1],
\]
where $\vk_X(z):=\left(K(z,x_1),\cdots,K(z,x_N)\right)^T$.
For example in Figure~\ref{fig:1DExaKriging}, if the shape parameter of the Gaussian kernel is endowed with $\theta=6$, then we can obtain the simple kriging prediction $\hat{s}(z)$ locally which is also consistent with the kernel-based interpolant $\hat{u}$ in Figure~\ref{fig:1DExaKernel}, that is, $\hat{s}(z)=\hat{u}(z)$.
Recently, the paper~\cite{ScheuererSchabackSchlather2012} compares the spatial-data interpolations
between the deterministic problems (meshfree approximation) and the stochastic problems (kriging interpolation) and the paper~\cite{ScheuererSchabackSchlather2012} also shows that the representations of these both estimators are the same.
Therefore, we conjecture that the interpolating paths in Figure~\ref{fig:1DExaInitial} could be equivalently transferred into some Gaussian fields such that the best estimators could be measured by the related Gaussian random variables.

Fortunately, the constructions of Brownian motions inspire the connections of the interpolating paths and the Gaussian fields.
It is well-known that the standard Brownian motion $W$ is a Gaussian field with the mean $0$ and the covariance kernel $K(t,s):=\min\left\{t,s\right\}$ which is also a positive definite kernel.
\cite[Chapter~2]{KaratzasShreve1991} provides various constructions of Brownian motions and one kind of the constructions is defined on the continuous function space $\Cont[0,\infty)$.
More precisely, the Wiener measure $\PP_{\ast}$ is well-posed on the sample space $\left(\Omega_{\ast},\Filter_{\ast}\right)$ composed of the function space $\Cont[0,\infty)$ and the Borel $\sigma$-algebra $\Borel\left(\Cont[0,\infty)\right)$.
By \cite[Theorem~4.20]{KaratzasShreve1991}, the coordinate mapping process $W_t(\omega):=\omega(t)$ for $t\in[0,\infty)$ and $\omega\in\Omega_{\ast}$ is a standard Brownian motion on the probability space $\left(\Omega_{\ast},\Filter_{\ast},\PP_{\ast}\right)$.
This shows that we can connect the continuous paths to the Brownian motions.
Moreover, we find that the initial condition $Y_0=y_0$ of the simple stochastic ordinary differential equation $\ud Y_t=\ud B_t$
is equivalent to the interpolation at the origin (see \cite[Section~5.2]{Oksendal2003}).
By the construction of the Brownian motions, we obtain an extension of the interpolating paths $\hat{u}_1,\hat{u}_2,\hat{u}_3$ in Figure~\ref{fig:1DExaInitial}
to all interpolating paths in $\Cont[0,\infty)$ such as the interpolation
$\Aset_X(\vf):=\left\{\omega\in\Cont[0,\infty):\omega(x_1)=f_1,\ldots,\omega(x_N)=f_N\right\}$ can be measured by the multivariate normal random variables $W_{x_1},\ldots,W_{x_N}$.

Therefore, we believe that the meshfree approximation and the kriging interpolation can be strongly connected by the Gaussian fields with the
analogous structures of the Brownian motions.
In \cite{CialencoFasshauerYe2012,YeThesis2012},
we extend the initial ideas in Figure~\ref{fig:1DExaInitial} to all interpolating paths in a reproducing kernel Hilbert space $\Hilbert_\Phi([0,1])$ (see Definition~\ref{d:RKHS}).
The theorems in \cite{CialencoFasshauerYe2012,YeThesis2012} guarantee that
the Gaussian field $S_x(\omega):=\omega(x)$ is well-defined on $\Hilbert_\Phi([0,1])$ similar as the Brownian motion defined on $\Cont[0,\infty)$. Moreover,
the covariance kernel of this Gaussian field $S$ is the integral-type kernel $\Phistar$ of the reproducing kernel $\Phi$, that is,
$\Phistar(x,y):=\int_0^1\int_0^1\Phi(x,t)\Phi(y,t)\ud t$.

In the following section, we will improve the theorems in \cite{CialencoFasshauerYe2012,YeThesis2012} to endow the Sobolev spaces with the probability measures induced by the positive definite kernels such that
the interpolating paths can be measured by the Gaussian random variables (see Theorem~\ref{t:Gauss-PDK-L} and Lemma~\ref{l:Gauss-PDK}).
In Section~\ref{sec:ker-app},
the kernel-based probability structures of the Sobolev spaces will provide the best estimators induced by the positive definite kernels (see Figure~\ref{fig:Sobolev-initial} which is the generalization of the initial ideas in Figure~\ref{fig:1DExaInitial} and Table~\ref{tab:InitialIdea}).

%---------------------------------------------------------------------------------------------------------------------
\section{Gaussian Random Variables and Positive Definite Kernels}\label{sec:Gauss-PDK}
%---------------------------------------------------------------------------------------------------------------------

In this section, we firstly study with the constructions of various multivariate normal random variables defined on the $\Leb_2$-based Sobolev space $\Hilbert^m(\Domain)$ of the degree $m$ by the deterministic bounded linear functionals and the given positive definite kernels.

For convenience of the proofs, we let $\Domain$ be a \emph{regular and compact} domain of the $d$-dimensional real space $\Rd$ in this article (the details of the regularity such as a Lipschitz domain are mentioned in~\cite[Section 4.1]{AdamsFournier2003}).

Next, we will discuss how to use these multivariate normal random variables to approximate the unknown value $Lu$ by the given data information
\[
f_1:=L_1u,\ldots,f_N:=L_Nu,
\]
where the target function $u\in\Hilbert^m(\Domain)$ and $L,L_1,\ldots,L_N$ are the bounded (continuous) linear functionals on $\Hilbert^m(\Domain)$, for example, the point evaluation function $\delta_{\vx}$, the partial derivative $\delta_{\vx}\circ D^{\valpha}$, or the integral $\int_{\Domain}$.

\subsection{Constructing Gaussian Random Variables by Positive Definite Kernels}\label{sec:ConstrGaussPDK}

Now we generalize \cite[Theorem~3.1]{CialencoFasshauerYe2012} to
endow the Sobolev spaces with the probability measures induced by the positive definite kernels such that
we can obtain the normal random variables indexed by the given bounded linear functionals.

%////////////////////////////////////////////////////////////////////////////////////////////////////////////////////////
\begin{theorem}\label{t:Gauss-PDK-L}
Suppose that the positive definite kernel $K\in\Cont^{2m,1}\left(\Domain\times\Domain\right)$ for $m>d/2$.
Let $L$ be a bounded linear functional on the Sobolev space $\Hilbert^m(\Domain)$. Then there exists a probability measure $\PP_K$ on the measurable space
\[
\left(\Omega_m,\Filter_m\right):=\left(\Hilbert^m(\Domain),\Borel\left(\Hilbert^m(\Domain)\right)\right),
\]
such that the normal random variable
\[
LS(\omega):=L\omega,\quad\text{for }\omega\in\Omega_m,
\]
is well-defined on the probability space $\left(\Omega_m,\Filter_m,\PP_K\right)$
and this random variable $LS$ has the mean $0$ and the variance $L_{\vx}L_{\vy}K(\vx,\vy)$.
Moreover, the probability measure $\PP_K$ is independent of the bounded linear functional $L$.
\end{theorem}
%////////////////////////////////////////////////////////////////////////////////////////////////////////////////////////

%////////////////////////////////////////////////////////////////////////////////////////////////////////////////////////
\begin{remark}\label{r:Gauss-PDK-L}
In Theorem~\ref{t:Gauss-PDK-L},
the collection $\Borel\left(\Hilbert^m(\Domain)\right)$ represents the \emph{Borel $\sigma$-algebra} in the Sobolev space $\Hilbert^m(\Domain)$ and
$\omega\in\Omega_m$ represents the \emph{sample path} (trajectory).
The space $\Cont^{2m,1}\left(\Domain\times\Domain\right)\subseteq\Cont^{2m}\left(\Domain\times\Domain\right)$
consists of all functions which have the continuous derivatives up to order $2m$ and of which the $2m$th partial derivatives satisfy the Lipschitz condition.
Moreover, the notations $L_{\vx}$ and $L_{\vy}$ denote the linear operator $L$ associated to the first and second arguments of $\vx$ and $\vy$, respectively, that is, $L_{\vx}K(\vx,\vy)=L\left(K(\cdot,\vy)\right)$ and $L_{\vy}K(\vx,\vy)=L\left(K(\vx,\cdot)\right)$.
\end{remark}
%////////////////////////////////////////////////////////////////////////////////////////////////////////////////////////

Before the proofs of Theorem~\ref{t:Gauss-PDK-L}, we review some basic concepts of positive definite kernels, reproducing kernels, and Gaussian fields.

%////////////////////////////////////////////////////////////////////////////////////////////////////////////////////////
\begin{definition}[{\cite[Definition~6.24]{Wendland2005}}]\label{d:PDK}
A symmetric kernel $K:\Domain\times\Domain\to\RR$ is called \emph{positive definite} if, for any $N\in\NN$ and any distinct points $X:=\left\{\vx_1,\ldots,\vx_N\right\}\subseteq\Domain$, the quadratic form
\[
\sum_{j,k=1}^Nc_jc_kK(\vx_j,\vx_k)>0,\quad \text{for all }\vc:=\left(c_1,\cdots,c_N\right)^T\in\RR^N\setminus\{\vzero\}.
\]
\end{definition}
%////////////////////////////////////////////////////////////////////////////////////////////////////////////////////////
Definition~\ref{d:PDK} assures that all positive definite kernels are \emph{symmetric} in this article.
Obviously, all associated matrixes $\vK_{X}$ of the positive definite kernel $K$ are \emph{strictly} positive definite because $\vc^T\vK_{X}\vc>0$ for all $\vc\in\RR^N\setminus\{\vzero\}$.

%////////////////////////////////////////////////////////////////////////////////////////////////////////////////////////
\begin{definition}[{\cite[Definition~10.1]{Wendland2005}}]\label{d:RKHS}
A kernel $K:\Domain\times\Domain\to\RR$ is called a \emph{reproducing kernel} of a \emph{reproducing kernel Hilbert space} $\Hilbert_K(\Domain)$ composing of functions $f:\Domain\to\RR$ if
\[
\text{(i) }K(\cdot,\vy)\in\Hilbert_K(\Domain)\text{ and (ii) }f(\vy)=\left(f,K(\cdot,\vy)\right)_{\Hilbert_K(\Domain)},
\]
for all $\vy\in\Domain$ and all $f\in\Hilbert_K(\Domain)$,
where $\left(\cdot,\cdot\right)_{\Hilbert_K(\Domain)}$ is an inner product of the Hilbert space $\Hilbert_K(\Domain)$.
\end{definition}
%////////////////////////////////////////////////////////////////////////////////////////////////////////////////////////
\cite[Theorem~10.10]{Wendland2005} guarantees that any positive definite kernel is a reproducing kernel and its reproducing kernel Hilbert space exists uniquely.

%////////////////////////////////////////////////////////////////////////////////////////////////////////////////////////
\begin{example}[{\cite[Example~5.7]{FasshauerYe2011Dist}}]\label{ex:Sobolev-spline}
The typical example of positive definite kernels and reproducing kernels is the \emph{Sobolev-spline kernel (Mat\'ern function)} $G_m$ of the degree $m>d/2$,
that is,
\begin{equation}\label{eq:Sobolev-spline-kernel}
G_m(\vx,\vy):=\frac{2^{1-m-d/2}}{\pi^{d/2}\Gamma(m)}\norm{\vx-\vy}_2^{m-d/2}\mathcal{K}_{d/2-m}
\left(\norm{\vx-\vy}_2\right),\quad\text{for }\vx,\vy\in\Rd,
\end{equation}
where $\Gamma$ is the Gamma function and $\mathcal{K}_{\nu}$ is the modified Bessel function of
the second kind of order $\nu$.
According to the discussions in \cite{FasshauerYe2011Dist,FasshauerYe2012DiffBound},
the Sobolev-spline kernel $G_m$ is a positive definite kernel and its reproducing kernel Hilbert space $\Hilbert_{G_m}(\Rd)$
is equivalent to the Sobolev space $\Hilbert^m(\Rd)$.
Since the domain $\Domain$ is regular, \cite[Corollary~10.48]{Wendland2005} (the restrictions of reproducing kernel Hilbert spaces and Sobolev spaces) also guarantees that $\Hilbert_{G_m}(\Domain)$ and $\Hilbert^m(\Domain)$ are isomorphic.
This indicates that the spaces $\Hilbert_{G_m}(\Domain)=\Hilbert^m(\Domain)$ and the Borel $\sigma$-algebras $\Borel\left(\Hilbert_{G_m}(\Domain)\right)=\Borel\left(\Hilbert^m(\Domain)\right)$.
The condition of $m>d/2$ is sufficient to assure that $G_m\in\Cont\left(\Domain\times\Domain\right)$ and the point evaluation function $\delta_{\vx}$ is continuous on $\Hilbert^m(\Domain)$ by the Sobolev imbedding theorem~\cite[Theorem~4.12]{AdamsFournier2003}.
\end{example}
%////////////////////////////////////////////////////////////////////////////////////////////////////////////////////////

%////////////////////////////////////////////////////////////////////////////////////////////////////////////////////////
\begin{definition}[{\cite[Definition~3.28]{BerlinetThomas2004}}]\label{d:Gaussian}
A stochastic field $S:\Domain\times\Omega\to\RR$ defined on a probability space $\left(\Omega,\Filter,\PP\right)$ is called a \emph{Gaussian field} with a mean $0$ and a covariance kernel $K:\Domain\times\Domain\to\RR$ if, for any $N\in\NN$ and any distinct points $X:=\left\{\vx_1,\ldots,\vx_N\right\}\subseteq\Domain$, the random vector $\vS_{X}:=\left(S_{\vx_1},\cdots,S_{\vx_N}\right)^T$ is a multivariate normal random vector with the mean $\vzero$ and the covariance matrix $\vK_{X}$, that is, $\vS_X\sim\Normal\left(\vzero,\vK_{X}\right)$.
\end{definition}
%////////////////////////////////////////////////////////////////////////////////////////////////////////////////////////

%////////////////////////////////////////////////////////////////////////////////////////////////////////////////////////
\begin{remark}\label{r:Gaussian}
In stochastic analysis \cite{KaratzasShreve1991} and probability theory \cite{Shiryaev1996}, the \emph{measurable space} $\left(\Omega,\Filter\right)$ is called a \emph{sample space} and the \emph{$\sigma$-algebra} $\Filter$ in $\Omega$ is called a \emph{filtration}.
Next, we illustrate some specific notations of the Gaussian field $S$.
For any fixed point $\vx\in\Domain$, the symbol $S_{\vx}$ represents a random variable defined on the probability space $\left(\Omega,\Filter,\PP\right)$.
In another hands, for any fixed sample $\omega\in\Omega$, the symbol $\vx\mapsto S_{\vx}(\omega)$ or $S(\omega)$ represents a deterministic function defined on the domain $\Domain$. Since the mean of $S$ is equal to $0$, we have
$\Mean\big(S_{\vx}\big)=\Mean\big(S_{\vy}\big)=0$; hence
\[
K(\vx,\vy)=\Cov\big(S_{\vx},S_{\vy}\big)
=\Mean\big(S_{\vx}S_{\vy}\big)-\Mean\big(S_{\vx}\big)\Mean\big(S_{\vy}\big)
=\Mean\big(S_{\vx}S_{\vy}\big),
\]
for any $\vx,\vy\in\Domain$. This indicates that the covariance matrix $\vK_{X}$ of the random vector $\vS_{X}$ can be computed by
\[
\vK_X
=
\left(K(\vx_j,\vx_k)\right)_{j,k=1}^{N,N}=
\left(\Cov\big(S_{\vx_j},S_{\vx_k}\big)\right)_{j,k=1}^{N,N}=
\left(\Mean\big(S_{\vx_j}S_{\vx_k}\big)\right)_{j,k=1}^{N,N}.
\]
In this article, all equalities of random variables and stochastic fields are equal almost surely without any specific illustration.
\end{remark}
%////////////////////////////////////////////////////////////////////////////////////////////////////////////////////////

To prove the theorems and lemmas, we need to study with the properties of the positive definite kernel $K\in\Cont^{2m,1}\left(\Domain\times\Domain\right)$ given in Theorem~\ref{t:Gauss-PDK-L}.
Since $\Domain$ is compact and $K$ is symmetric and continuous, the Mercer's theorem~\cite[Theorem~4.49]{SteinwartChristmann2008} guarantees that
there exist a countable set of \emph{eigenvalues} $\lambda_1\geq\lambda_2\geq\cdots>0$ and orthonormal \emph{eigenfunctions} $\left\{e_n\right\}_{n=1}^{\infty}$ in $\Leb_2(\Domain)$ such that
\[
\lambda_ne_n(\vx)=\int_{\Domain}K(\vx,\vy)e_n(\vy)\ud\vy,\quad\text{for all }n\in\NN,
\]
and the positive definite kernel $K$ possesses the absolutely and uniformly convergent representation
\begin{equation}\label{eq:PDK-eigen-rep}
K(\vx,\vy)=\sum_{n=1}^{\infty}\lambda_ne_n(\vx)e_n(\vy),\quad \text{for }\vx,\vy\in\Domain.
\end{equation}
Since $K\in\Cont^{2m,1}\left(\Domain\times\Domain\right)$, we have
\[
D^{\valpha}e_n(\vx)=\lambda_n^{-1}\int_{\Domain}D^{\valpha}_{\vx}K(\vx,\vy)e_n(\vy)\ud\vy,
\quad\text{for }\valpha\in\NN_0^d\text{ with }\abs{\valpha}\leq 2m,
\]
where
$D^{\valpha}:=\partial^{\valpha}/\partial\vx^{\valpha}$ is the partial derivative of order $\valpha$; hence
$\left\{e_n\right\}_{n=1}^{\infty}\subseteq\Cont^{2m}(\Domain)$.
This indicates that
the representation
\[
D_{\vx}^{\valpha}D_{\vy}^{\vbeta}K(\vx,\vy)=\sum_{n=1}^{\infty}\lambda_nD^{\valpha}e_n(\vx)D^{\vbeta}e_n(\vy),
\quad\text{for }\vx,\vy\in\Domain,
\]
converges absolutely and uniformly for any $\valpha,\vbeta\in\NN_0^d$ with $\abs{\valpha}+\abs{\vbeta}\leq 2m$.
The compactness of the domain $\Domain$ assures that $\Cont^{2m,1}\left(\Domain\times\Domain\right)\subseteq\Hilbert^{2m}\left(\Domain\times\Domain\right)$ and $L_{\vx}L_{\vy}K(\vx,\vy)$ is well-posed for the bounded linear functional $L$ on $\Hilbert^{m}(\Domain)$. Moreover
$L_{\vx}L_{\vy}K(\vx,\vy)$ also possesses the convergent representation
\[
L_{\vx}L_{\vy}K(\vx,\vy)=\sum_{n=1}^{\infty}\lambda_nL_{\vx}\left(e_n(\vx)\right)L_{\vy}\left(e_n(\vy)\right)
=\sum_{n=1}^{\infty}\lambda_n\left(Le_n\right)^2.
\]

For the proofs of Theorem~\ref{t:Gauss-PDK-L}, we need the generalization of \cite[Lemma~2.2]{CialencoFasshauerYe2012} which guarantees that
there exists a probability measure $\PP_{\Phistar}$ induced by the integral-type kernel $\Phistar$ of the reproducing kernel $\Phi$ such that the Gaussian field $S_{\vx}(\omega):=\omega(\vx)$ is well-defined on the
reproducing kernel Hilbert space $\Hilbert_\Phi(\Domain)$.
Roughly speaking, we will extend the original relationships
\[
\Hilbert_\Phi(\Domain)\longleftrightarrow \Phi\longleftrightarrow \Phistar\longleftrightarrow \PP_{\Phistar}\longleftrightarrow S,
\]
in \cite[Lemma~2.2]{CialencoFasshauerYe2012} to another general forms
\[
\Hilbert^m(\Domain)\cong\Hilbert_{G_m}(\Domain)\longleftrightarrow G_m\longleftrightarrow K\in\Cont^{2m,1}\left(\Domain\times\Domain\right)\longleftrightarrow \PP_{K}\longleftrightarrow S,
\]
in Lemma~\ref{l:Gauss-PDK} which shows that there exists a probability measure $\PP_K$ induced by the given positive definite kernel $K$ such that the Gaussian field $S_{\vx}(\omega):=\omega(\vx)$ is well-defined on the Sobolev space $\Hilbert^m(\Domain)$.

Usually, it is difficult to obtain the probability measure $\PP_K$ directly on the Sobolev space $\Hilbert^m(\Domain)$.
In the proofs of \cite[Lemma~2.2]{CialencoFasshauerYe2012}, a Gaussian field $\xi$, which is easily constructed by the integral-type kernel $\Phistar$, is a primary element to introduce the probability measure $\PP_{\Phistar}$,
and the main technique is based on the theorems in \cite{LukicBeder2001}.
Same as this idea, we will also use \cite[Lemma~2.1 and Theorem~3.2]{LukicBeder2001} to verify Lemma~\ref{l:Gauss-PDK}, that is, the extensions of the original proofing process
$\Phistar\rightarrow\xi\rightarrow\PP_{\xi}=\PP_{\Phistar}\rightarrow S$
to
$K\rightarrow\xi\rightarrow\PP_{\xi}=\PP_{K}\rightarrow S$.
So, the proofs of Lemma~\ref{l:Gauss-PDK} will be separated into two steps: the first step is to construct a Gaussian field $\xi$ by the given positive definite kernel $K$, and
we will introduce the probability measure $\PP_K$ by this Gaussian field $\xi$ in the next step.

For convenience, we repeat \cite[Lemma~2.1 and Theorem~3.2]{LukicBeder2001} in Lemma~\ref{l:Gauss-RKHS-sample} which is consistent with the formats of this article.

%////////////////////////////////////////////////////////////////////////////////////////////////////////////////////////
\begin{lemma}[{\cite[Lemma~2.1 and Theorem~3.2]{LukicBeder2001}}]\label{l:Gauss-RKHS-sample}
Suppose that a Gaussian field $\xi$ defined on a probability space $\left(\Omega,\Filter,\PP\right)$ belongs to a reproducing kernel Hilbert space $\Hilbert_{G}(\Domain)$ almost surely, that is, $\PP\left(\xi\in\Hilbert_{G}(\Domain)\right)=1$. Then the probability measure
\[
\PP_{\xi}(A):=\PP\left(\xi^{-1}(A)\right),\quad\text{for }A\in\Borel\left(\Hilbert_{G}(\Domain)\right),
\]
is well-posed on the measurable space $\left(\Hilbert_{G}(\Domain),\Borel\left(\Hilbert_{G}(\Domain)\right)\right)$.
Moreover, the Gaussian field
\[
S_{\vx}(\omega):=\omega(\vx),\quad\text{for }\vx\in\Domain\text{ and }\omega\in\Hilbert_{G}(\Domain),
\]
is well-defined on the probability space $\big(\Hilbert_{G}(\Domain),\Borel\left(\Hilbert_{G}(\Domain)\right),\PP_{\xi}\big)$
and the means and covariance kernels of the Gaussian fields $S$ and $\xi$ are the same.
\end{lemma}
%////////////////////////////////////////////////////////////////////////////////////////////////////////////////////////

%////////////////////////////////////////////////////////////////////////////////////////////////////////////////////////
\begin{remark}\label{r:Almost-Surely}
Here, we call a stochastic field $\xi$ belongs to a function space $\Hilbert$ almost surely if the function $\vx\mapsto\xi_{\vx}(\omega)$ belongs to $\Hilbert$ for $\omega\in\Omega$ almost surely, or the probability of the set $A:=\left\{\omega\in\Omega:\xi(\omega)\in\Hilbert\right\}$ is equal to $1$. In \cite{LukicBeder2001}, the stochastic field $\xi$ can be viewed as a measurable map from $\Omega$ to $\Hilbert_{G}(\Domain)$ such that $\PP\left(\xi_{\vx}\leq z\right)=\PP_{\xi}\left(S_{\vx}\leq z\right)$ for any $\vx\in\Domain$ and any $z\in\RR$.
This shows that $S$ and $\xi$ have the same probability distributions.
In fact, Lemma~\ref{l:Gauss-RKHS-sample} for the Gaussian fields is just a typical case of the theorems in \cite{LukicBeder2001} which can cover more general stochastic fields.
\end{remark}
%////////////////////////////////////////////////////////////////////////////////////////////////////////////////////////

Now we construct a Gaussian field $\xi$ with the mean $0$ and the covariance kernel $K$.
The Kolmogorov's extension theorem \cite[Theorem~2.3]{Shiryaev1996} guarantees that there exist a
countable independent standard normal random variables $\left\{\zeta_n\right\}_{n=1}^{\infty}$ defined on a probability space $\left(\Omega,\Filter,\PP\right)$.
A example of $\PP$ is the infinite-dimensional Gaussian measure placed on $\RR^{\NN}$ (see \cite[Section~2.3]{Shiryaev1996}).
Combining the random variables $\left\{\zeta_n\right\}_{n=1}^{\infty}$ with the eigenvalues $\left\{\lambda_n\right\}_{n=1}^{\infty}$ and eigenfunctions $\left\{e_n\right\}_{n=1}^{\infty}$ of the given positive definite kernel $K$ in Theorem~\ref{t:Gauss-PDK-L}, we construct a stochastic field on the probability space $\left(\Omega,\Filter,\PP\right)$ such as
\begin{equation}\label{eq:Gauss-xi}
\xi_{\vx}:=\sum_{n=1}^{\infty}\zeta_n\sqrt{\lambda_n}e_n(\vx),\quad \text{for }\vx\in\Domain.
\end{equation}
Since $\zeta_n\sim\text{i.i.d.}\Normal(0,1)$, we have
\[
\Mean\big(\zeta_n\big)=0,\quad
\Mean\big(\zeta_n^2\big)=1,\quad \Mean\big(\zeta_k\zeta_n\big)=\Mean\big(\zeta_k\big)\Mean\big(\zeta_n\big)=0\text{ when }k\neq n,
\]
for all $k,n\in\NN$.
Notes that
\[
\Mean\big(\xi_{\vx}^2\big)=\sum_{k,n=1}^{\infty}\Mean\big(\zeta_k\zeta_n\big)\sqrt{\lambda_k\lambda_n}e_k(\vx)e_n(\vx)
=\sum_{n=1}^{\infty}\lambda_ne_n(\vx)^2=K(\vx,\vx)<\infty;
\]
hence the stochastic field $\xi$ is well-defined.

%////////////////////////////////////////////////////////////////////////////////////////////////////////////////////////
\begin{lemma}\label{l:Gauss-xi}
The stochastic field $\xi$ given in Equation~\eqref{eq:Gauss-xi} is a Gaussian field with the mean $0$ and the covariance kernel $K$.
\end{lemma}
%////////////////////////////////////////////////////////////////////////////////////////////////////////////////////////
\begin{proof}
Since the linear combination of normal random variables is still normal (see~\cite[Theorem~A.17 and A.19]{Oksendal2003}), the random variable $\xi_{\vx}$ is normal for any $\vx\in\Domain$.
Next, we compute the mean and the covariance kernel of $\xi$.
Take any $\vx,\vy\in\Domain$. Equation~\eqref{eq:Gauss-xi} assures that
\[
\Mean\left(\xi_{\vx}\right)=\sum_{n=1}^{\infty}\Mean\left(\zeta_n\right)\sqrt{\lambda_n}e_n(\vx)=0,
\]
and
\begin{align*}
\Cov\big(\xi_{\vx},\xi_{\vy}\big)
=\Mean\big(\xi_{\vx}\xi_{\vy}\big)
=\sum_{k,n=1}^{\infty}\Mean\big(\zeta_k\zeta_n\big)\sqrt{\lambda_k\lambda_n}e_k(\vx)e_n(\vy)
=\sum_{n=1}^{\infty}\lambda_ne_n(\vx)e_n(\vy)=K(\vx,\vy).
\end{align*}
\end{proof}

Lemma~\ref{l:Gauss-xi} shows that the Gaussian field $\xi$ is a centered Gaussian field, that is, a Gaussian field with the mean $0$, and the covariance kernel $K$ of $\xi$ belongs to $\Cont^{2m,1}\left(\Domain\times\Domain\right)$; hence the Kolmogorov-\v{C}entsov continuity theorem in \cite[Section~2.2.B]{KaratzasShreve1991} guarantees that:

%////////////////////////////////////////////////////////////////////////////////////////////////////////////////////////
\begin{lemma}\label{l:Gauss-Sobolev}
The Gaussian field $\xi$ given in Equation~\eqref{eq:Gauss-xi} belongs to $\Cont^m(\Domain)$ almost surely, that is, $\PP\left(\xi\in\Cont^m(\Domain)\right)=1$.
\end{lemma}
%////////////////////////////////////////////////////////////////////////////////////////////////////////////////////////

%////////////////////////////////////////////////////////////////////////////////////////////////////////////////////////
\begin{remark}\label{r:Gauss-Sobolev}
In fact, by the Karhunen representation theorem~\cite[Theorem~3.41]{BerlinetThomas2004},
Equation~\eqref{eq:Gauss-xi} can be also seen as the Karhunen-Lo\`{e}ve expansion of the Gaussian field $\xi$.
Now we take any $\valpha\in\NN_0^d$ with $\abs{\valpha}\leq m$ to
construct a stochastic field
\[
\xi_{\valpha,\vx}:=\sum_{n=1}^{\infty}\zeta_n\sqrt{\lambda_n}D^{\valpha}e_n(\vx),\quad \text{for }\vx\in\Domain.
\]
Since the mean square
\[
\Mean\big(\xi_{\valpha,\vx}^2\big)
=\sum_{n=1}^{\infty}\lambda_nD^{\valpha}e_n(\vx)^2
=D^{\valpha}_{\vz_1}D^{\valpha}_{\vz_2}K\left(\vz_1,\vz_2\right)|_{\vz_1=\vz_2=\vx}<\infty,
\]
the stochastic field $\xi_{\valpha}$ is well-defined. Same as the properties of the Karhunen-Lo\`{e}ve expansion, the expansion of $\xi_{\valpha}$ is uniformly convergent on the compact domain $\Domain$ because $D^{\valpha}_{\vx}D^{\valpha}_{\vy}K\in\Cont^{0,1}\left(\Domain\times\Domain\right)$. Combining with $D^{\valpha}\xi=\xi_{\valpha}$, the expansion of $\xi$ is also convergent in $\Cont^m(\Domain)$.

Even though the eigenfunctions $e_n\in\Cont^{2m}(\Domain)$ for all $n\in\NN$, we still can not determine whether $D^{\valpha}\xi\in\Cont(\Domain)$ when $\abs{\valpha}>m$ because $D^{\valpha}_{\vz_1}D^{\valpha}_{\vz_2}K\left(\vz_1,\vz_2\right)|_{\vz_1=\vz_2=\vx}$ may not exist for all $\vx\in\Domain$.
This indicates that the smoothing sample paths of the Gaussian fields can not be determined if their covariance kernels are non-smooth.
\end{remark}
%////////////////////////////////////////////////////////////////////////////////////////////////////////////////////////

By the smoothness of $\xi$, we can further check that:
%////////////////////////////////////////////////////////////////////////////////////////////////////////////////////////
\begin{lemma}\label{l:Gauss-Sobolev-RKHS}
The Gaussian field $\xi$ given in Equation~\eqref{eq:Gauss-xi} belongs to the reproducing kernel Hilbert space $\Hilbert_{G_m}(\Domain)$ almost surely, that is, $\PP\left(\xi\in\Hilbert_{G_m}(\Domain)\right)=1$, where
the Sobolev-spline kernel $G_m$ of the degree $m$ is given in Equation~\eqref{eq:Sobolev-spline-kernel}.
\end{lemma}
%////////////////////////////////////////////////////////////////////////////////////////////////////////////////////////
\begin{proof}
Lemma~\ref{l:Gauss-Sobolev} assures that $\xi\in\Cont^m(\Domain)$ almost surely.
Since $\Domain$ is compact, we have $\Cont^m(\Domain)\subseteq\Hilbert^m(\Domain)$; hence $\xi\in\Hilbert^m(\Domain)$ almost surely.
Moreover, the discussions in Example~\ref{ex:Sobolev-spline} show that $\Hilbert^m(\Domain)$ is equivalent to $\Hilbert_{G_m}(\Domain)$.
This assures that $\xi\in\Hilbert_{G_m}(\Domain)$ almost surely.
\end{proof}

Combining with Lemmas~\ref{l:Gauss-RKHS-sample},~\ref{l:Gauss-xi}, and~\ref{l:Gauss-Sobolev-RKHS}, we can complete the proofs of Lemma~\ref{l:Gauss-PDK} for the constructions of the probability measure $\PP_K$ in Theorem~\ref{t:Gauss-PDK-L}.

%////////////////////////////////////////////////////////////////////////////////////////////////////////////////////////
\begin{lemma}\label{l:Gauss-PDK}
Suppose that the positive definite kernel $K\in\Cont^{2m,1}\left(\Domain\times\Domain\right)$ for $m>d/2$. Then there exists a probability measure $\PP_K$ on the measurable space
\[
\left(\Omega_m,\Filter_m\right):=\left(\Hilbert^m(\Domain),\Borel\left(\Hilbert^m(\Domain)\right)\right),
\]
such that the Gaussian field
\[
S_{\vx}(\omega):=\omega(\vx),\quad\text{for }\vx\in\Domain\text{ and }\omega\in\Omega_m,
\]
is well-defined on the probability space $\left(\Omega_m,\Filter_m,\PP_K\right)$ and this Gaussian field $S$ has the mean $0$ and the covariance kernel $K$.
\end{lemma}
%////////////////////////////////////////////////////////////////////////////////////////////////////////////////////////
\begin{proof}
Firstly, by Example~\ref{ex:Sobolev-spline}, we have $\Omega_m=\Hilbert^m(\Domain)=\Hilbert_{G_m}(\Domain)$ and $\Filter_m=\Borel\left(\Hilbert^m(\Domain)\right)=\Borel\left(\Hilbert_{G_m}(\Domain)\right)$.

By Lemma~\ref{l:Gauss-xi}, the Gaussian field $\xi$ in Equation~\eqref{eq:Gauss-xi} has the mean $0$ and the covariance kernel $K$, and Lemma~\ref{l:Gauss-Sobolev-RKHS} provides that $\xi\in\Hilbert_{G_m}(\Domain)$ almost surely.
Therefore, Lemmas~\ref{l:Gauss-RKHS-sample} guarantees that the Gaussian field $\xi$ can be used to introduce the probability measure $\PP_K$ on the measurable space $\left(\Omega_m,\Filter_m\right)$, that is,
\[
\PP_{K}(A):=\PP\left(\xi^{-1}(A)\right),
\quad\text{for }A\in\Filter_m,
\]
such that $S_{\vx}(\omega):=\omega(\vx)$ is a Gaussian field with the mean $0$ and the covariance kernel $K$ placed on the probability space $\left(\Omega_m,\Filter_m,\PP_K\right)$.
\end{proof}

%////////////////////////////////////////////////////////////////////////////////////////////////////////////////////////
\begin{remark}\label{r:Gauss-PDK}
The probability measure $\PP_K$ in Lemma~\ref{l:Gauss-PDK} can be seen as the generalization of the
the Wiener measure $\PP_{\ast}$ on the continuous function space $\Cont[0,\infty)$ for the Brownian motion $W$ (see the discussions in Section~\ref{sec:Initial}), that is,
$\Cont[0,\infty)\to\Hilbert^m(\Domain)$, $W\to S$, and $\PP_{\ast}\to\PP_K$.
There may be another methods to introduce the probability measure $\PP_K$ directly
by the cylinder sets or the cylindrical $\sigma$-algebra in the Sobolev spaces such as the Wiener measures on the continuous function spaces in \cite{KaratzasShreve1991} and
the Gaussian measures on the reproducing kernel Hilbert spaces in \cite{Janson1997}.
This means that there may be another generalizations of Lemma~\ref{l:Gauss-PDK}, for example,
$\Cont^{2m,1}\left(\Domain\times\Domain\right)\to\Hilbert^{2m}\left(\Domain\times\Domain\right)$ for non-smooth kernels or $\Hilbert^m(\Domain)\to\mathcal{W}^m_p(\Domain)$ for Sobolev Banach spaces.
In this article, we do not discuss another constructions and proofs of the probability measure $\PP_K$ deeply.
\end{remark}
%////////////////////////////////////////////////////////////////////////////////////////////////////////////////////////

The Gaussian field $S$ can be viewed as the invariant element of the original Gaussian field $\xi$;
hence we can also obtain the \emph{Karhunen-Lo\`{e}ve expansion} of $S$ same as the discussions of $\xi$ in Remark~\ref{r:Gauss-Sobolev}, that is,
\begin{equation}\label{eq:KL-expansion-S}
S=\sum_{n=1}^{\infty}\eta_n\sqrt{\lambda_n}e_n,
\end{equation}
where $\left\{\eta_n\right\}_{n=1}^{\infty}$ are the i.i.d. standard normal random variables defined on the probability space $\left(\Omega_m,\Filter_m,\PP_K\right)$. Here, the random variables $\left\{\eta_n\right\}_{n=1}^{\infty}$ can be also thought as the invariant elements of $\left\{\zeta_n\right\}_{n=1}^{\infty}$ in Equation~\eqref{eq:Gauss-xi}.
The both random coefficients $\left\{\eta_n\right\}_{n=1}^{\infty}$ and $\left\{\zeta_n\right\}_{n=1}^{\infty}$ are just defined on the different probability spaces.

Same as the smoothness of $\xi$, the Gaussian field $S$ belongs to $\Cont^m(\Domain)$ almost surely, that is,
$\PP_K\left(S\in\Cont^m(\Domain)\right)=1$, and the expansion of $S$ in Equation~\eqref{eq:KL-expansion-S} is also convergent in $\Cont^m(\Domain)$.
Roughly, the probability measure $\PP_K$ vanishes the non-smooth paths in $\Hilbert^m(\Domain)$.

Finally, we verify the main theorem by Lemma \ref{l:Gauss-PDK} as follows.

\begin{proof}[{\bf \emph{Proof of Theorem~\ref{t:Gauss-PDK-L}}}]
Firstly, Lemma~\ref{l:Gauss-PDK} guarantees that the probability measure $\PP_K$ induced by the positive definite kernel $K$ is well-posed on the measurable space $\left(\Omega_m,\Filter_m\right)$ and the random variable $LS(\omega):=L\omega$ is also well-defined for the bounded linear functional $L$ on the Sobolev space $\Hilbert^m(\Domain)$.

Next, we will use the Karhunen-Lo\`{e}ve expansion of the Gaussian field $S$ in Equation~\eqref{eq:KL-expansion-S} to complete the proofs.
As the above discussions of Equation~\eqref{eq:KL-expansion-S}, the expansion $S=\sum_{n=1}^{\infty}\eta_n\sqrt{\lambda_n}e_n$ is convergent in $\Cont^m(\Domain)$; hence
the compactness of $\Domain$ assures that the expansion $S=\sum_{n=1}^{\infty}\eta_n\sqrt{\lambda_n}e_n$ is also convergent in $\Hilbert^m(\Domain)$.
Since the linear functional $L$ is bounded on $\Hilbert^m(\Domain)$, we have
\begin{equation}\label{eq:KL-expansion-S-L}
LS=\sum_{n=1}^{\infty}\eta_n\sqrt{\lambda_n}Le_n.
\end{equation}
Moreover, since $\eta_n\sim\text{i.i.d.}\Normal(0,1)$ for all $n\in\NN$, the random variable $LS$ is a linear combination of the normal random variables. Therefore, the random variable $LS$ is a normal random variable with the mean
\[
\Mean\left(LS\right)=\sum_{n=1}^{\infty}\Mean\left(\eta_n\right)\sqrt{\lambda_n}Le_n=0,
\]
and the variance
\[
\Var\left(LS\right)=\Mean\left(LS\right)^2
=\sum_{k,n=1}^{\infty}\Mean\left(\eta_k\eta_n\right)\sqrt{\lambda_k\lambda_n}Le_kLe_n
=\sum_{n=1}^{\infty}\lambda_n\left(Le_n\right)^2=L_{\vx}L_{\vy}K(\vx,\vy).
\]
This completes the proofs of the theorem.
\end{proof}

{\bf Applications of Theorem~\ref{t:Gauss-PDK-L}:}
In the approximation problems, we usually study with the finite many bounded linear functionals $L_1,\ldots,L_N$ on the Sobolev space $\Hilbert^m(\Domain)$. Thus, we will look at the the multivariate normal random variables
$L_1S,\ldots,L_NS$ defined as in Theorem~\ref{t:Gauss-PDK-L}.
Since the probability measure $\PP_K$ is independent of $L_1,\ldots,L_N$, all random variables $L_1S,\ldots,L_NS$ are placed on the same probability space $\left(\Omega_m,\Filter_m,\PP_K\right)$.
For convenience, we define the new notations
\[
\vL S:=\left(L_1S,\cdots,L_NS\right)^T,\quad \vL:=\left(L_1,\cdots,L_N\right)^T.
\]
Now we compute the means and covariances of these random variables by the convergent representation $L_jS=\sum_{n=1}^{\infty}\eta_n\sqrt{\lambda_n}L_je_n$ given in Equation~\eqref{eq:KL-expansion-S-L},
that is,
\[
\Mean\big(L_jS\big)=\sum_{n=1}^{\infty}\Mean\big(\eta_n\big)\sqrt{\lambda_n}L_je_n=0,
\]
and
\[
\Cov\big(L_jS,L_kS\big)=\Mean\big(L_jSL_kS\big)
=\sum_{n=1}^{\infty}\lambda_nL_je_nL_ke_n
=L_{j,\vx}L_{k,\vy}K(\vx,\vy),
\]
for all $j,k=1,\ldots,N$.
This indicates that the multivariate normal random vector $\vL S$ has the mean $\vzero$ and the covariance matrix
\[
\vK_{\vL}:=
\left(\Cov\big(L_jS,L_kS\big)\right)_{j,k=1}^{N,N}
=
\left(L_{j,\vx}L_{k,\vy}K(\vx,\vy)\right)_{j,k=1}^{N,N}.
\]

%////////////////////////////////////////////////////////////////////////////////////////////////////////////////////////
\begin{corollary}\label{c:Gauss-PDK-multi-L}
Suppose that the positive definite kernel $K\in\Cont^{2m,1}\left(\Domain\times\Domain\right)$ for $m>d/2$.
Let $\vL$ be composed of finite many bounded linear functionals $L_1,\ldots,L_N$ on the Sobolev space $\Hilbert^m(\Domain)$. Then the multivariate normal random vector
\[
\vL S(\omega):=\left(L_1S(\omega),\cdots,L_NS(\omega)\right)^T
=\left(L_1\omega,\cdots,L_N\omega\right)^T,\quad\text{for }\omega\in\Omega_m,
\]
is well-defined on the probability space $\left(\Omega_m,\Filter_m,\PP_K\right)$ given in Theorem~\ref{t:Gauss-PDK-L} and
and this random vector $\vL S$ has the mean $\vzero$ and the covariance matrix $\vK_{\vL}$.
\end{corollary}
%////////////////////////////////////////////////////////////////////////////////////////////////////////////////////////

%////////////////////////////////////////////////////////////////////////////////////////////////////////////////////////
\begin{remark}\label{r:Gauss-PDK-multi-L}
Even though the kernel $K$ is positive definite, we can not determine whether the covariance matrix $\vK_{\vL}$ is strictly positive definite.
But, we can assure that $\vK_{\vL}$ is always positive definite.
Thus, the \emph{pseudo inverse} $\vK_{\vL}^{\dag}$ of $\vK_{\vL}$ is well-defined by the eigen-decomposition of $\vK_{\vL}=\vV\vD\vV^T$, that is, $\vK_{\vL}^{\dag}=\vV\vD^{\dag}\vV^T$ (see \cite[Section~5.4]{KincaidWard2002}). Here $\vD$ and $\vV$ are composed of the nonnegative eigenvalues and orthonormal eigenvectors of $\vK_{\vL}$ and $\vD^{\dag}$ is taken by the reciprocal of each nonzero element on the diagonal of $\vD$, for example, $\vD^{\dag}:=\Diag\big(\lambda_1^{-1},\cdots,\lambda_k^{-1},0,\cdots,0\big)$ when $\vD:=\Diag\big(\lambda_1,\cdots,\lambda_k,0,\cdots,0\big)$.
\end{remark}
%////////////////////////////////////////////////////////////////////////////////////////////////////////////////////////

In the following section, we will discuss how to approximate $Lu$ by the given information $\vL u$.
Naturally, the relationships of $LS$ and $\vL S$ will be required in the following approximations.
So, we need to compute the conditional probability density function $p_{L|\vL}$ of $LS$ given $\vL S$.
By the basic probability theory, we have
\begin{equation}\label{eq:cond-join-prob-1}
p_{L|\vL}(v|\vv)=\frac{p_{L,\vL}(v,\vv)}{p_{\vL}(\vv)},\quad\text{for }v\in\RR\text{ and }\vv\in\RR^N,
\end{equation}
where $p_{L,\vL}$ and $p_{\vL}$ are the joint probability density functions of $\left(LS,\vL S\right)$ and $\vL S$, respectively.
Moreover,
Corollary~\ref{c:Gauss-PDK-multi-L} provides that
\[
\left(LS,\vL S\right)
\sim
\Normal\left(\vzero,\vK_{L,\vL}\right),
\quad
\vL S\sim
\Normal\left(\vzero,\vK_{\vL}\right),
\]
where the covariance matrix
\[
\vK_{L,\vL}
:=
\begin{pmatrix}
L_{\vx}L_{\vy}K(\vx,\vy)&L\vk_{\vL}^T\\
L\vk_{\vL}&\vK_{\vL}
\end{pmatrix},
\]
and the vector
\[
L\vk_{\vL}:=\left(L_{\vx}L_{1,\vy}K(\vx,\vy),\cdots,L_{\vx}L_{N,\vy}K(\vx,\vy)\right)^T,
\]
is computed by the kernel basis
\[
\vk_{\vL}(\vx):=\left(L_{1,\vy}K(\vx,\vy),\cdots,L_{N,\vy}K(\vx,\vy)\right)^T,\quad\text{for }\vx\in\Domain;
\]
hence
\begin{equation}\label{eq:cond-join-prob-2}
p_{L,\vL}(v,\vv):=
\frac{1}{\sqrt{\det_{\dag}\left(2\pi\vK_{L,\vL}\right)}}\exp\left(-\frac{1}{2}\vt^T\vK_{L,\vL}^{\dag}\vt\right),
\quad\text{for }\vt:=\begin{pmatrix}
v\\\vv
\end{pmatrix}\in\RR^{N+1},
\end{equation}
and
\begin{equation}\label{eq:cond-join-prob-3}
p_{\vL}(\vv):=\frac{1}{\sqrt{\det_{\dag}\left(2\pi\vK_{\vL}\right)}}\exp\left(-\frac{1}{2}\vv^T\vK_{\vL}^{\dag}\vv\right),
\quad\text{for }\vv\in\RR^N,
\end{equation}
where $\det_{\dag}$ is the \emph{pseudo determinant}, that is, the product of all nonzero eigenvalues of a positive definite matrix.
Combing with Equation~(\ref{eq:cond-join-prob-1}-\ref{eq:cond-join-prob-3}), we have:

%////////////////////////////////////////////////////////////////////////////////////////////////////////////////////////
\begin{corollary}\label{c:Gauss-PDK-L-pdf}
Let the normal random variable $LS$ and the multivariate normal random vector $\vL S$ be given in Theorem~\ref{t:Gauss-PDK-L} and Corollary~\ref{c:Gauss-PDK-multi-L}. Then the conditional probability density function $p_{L|\vL}$ of $LS$ given $\vL S$ can be written as
\[
p_{L|\vL}(v|\vv)=\frac{1}{\sigma_{L|\vL}\sqrt{2\pi}}\exp\left(-\frac{\left(v-\mu_{L|\vL}(\vv)\right)^2}{2\sigma_{L|\vL}^2}\right),
\quad\text{for }v\in\RR\text{ and }\vv\in\RR^N,
\]
where the mean
\[
\mu_{L|\vL}(\vv):=L\vk_{\vL}^T\vK_{\vL}^{\dag}\vv,
\]
and the standard deviation
\begin{equation}\label{eq:standard-dev}
\sigma_{L|\vL}:=\sqrt{L_{\vx}L_{\vy}K(\vx,\vy)-L\vk_{\vL}^T\vK_{\vL}^{\dag}L\vk_{\vL}}.
\end{equation}
\end{corollary}
%////////////////////////////////////////////////////////////////////////////////////////////////////////////////////////
Thus $\sigma_{L|\vL}^2$ is the variance of the conditional probability density function $p_{L|\vL}$. In Sections~\ref{sec:Stochastic-Data-Approx} and~\ref{sec:ell-SPDE}, we will show that the standard deviation $\sigma_{L|\vL}$ is equivalent to the (generalized) power functions in meshfree approximation.

%////////////////////////////////////////////////////////////////////////////////////////////////////////////////////////
\begin{remark}\label{r:notation}
Roughly, we can view the vector operator $\vL$ and the matrix operator $\vL_{\vx}\vL_{\vy}$ as the gradient and the Hessian matrix, respectively.
Thus, another good notations of the kernel basis $\vk_{\vL}(\vx)$ and the interpolating matrix $\vK_{\vL}$
can be rewritten as $\vL_{\vy}K(\vx,\vy)$ and $\vL_{\vx}\vL_{\vy}K(\vx,\vy)$, respectively.
\end{remark}
%////////////////////////////////////////////////////////////////////////////////////////////////////////////////////////

\subsection{Kernel-based Approximation for Deterministic Data}\label{sec:ker-app}

Now we study with the renewal kernel-based approximation by the multivariate normal random variables $LS,L_1S,\ldots,L_NS$
given in Theorem~\ref{t:Gauss-PDK-L} and Corollary~\ref{c:Gauss-PDK-multi-L}.
Let the vector
\[
\vf:=\left(f_1,\cdots,f_N\right)^T,
\]
be composed of the given data information $f_1=L_1u,\ldots,f_N=L_Nu$ evaluated by some deterministic function $u\in\Hilbert^m(\Domain)$ for $m>d/2$ and a vector bounded linear functional $\vL=\left(L_1,\cdots,L_N\right)^T$ on $\Hilbert^m(\Domain)$.
Then $\vL u=\vf$ and $Lu$ is well-defined for any bounded linear functional $L$ on $\Hilbert^m(\Domain)$.
Given the positive definite kernel $K\in\Cont^{2m,1}\left(\Domain\times\Domain\right)$,
we will construct the best estimator (kernel-based estimator) $\hat{v}=L\hat{u}_{\vL}$ or the kernel-based approximate function $\hat{u}_{\vL}$ to approximate the unknown value $Lu$ or the target function $u$.

\begin{wrapfigure}{r}{0.5\textwidth}
\vspace{-10pt}
\center
\includegraphics[width=0.5\textwidth, height=0.48\textwidth]{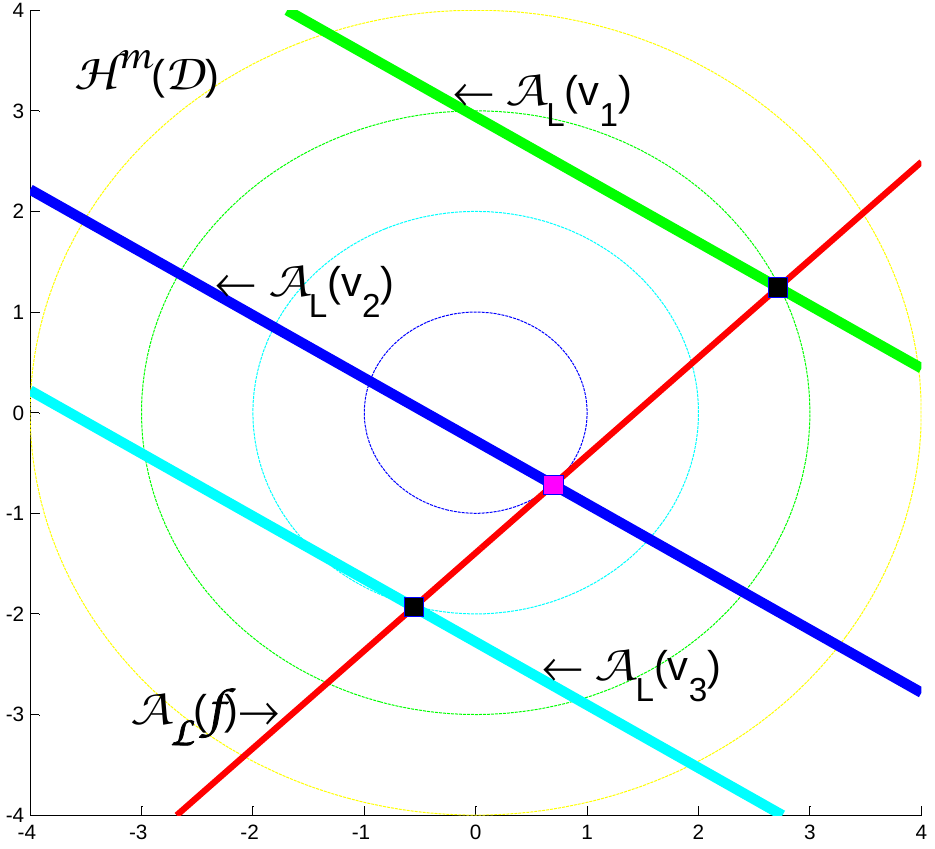}
\vspace{-20pt}
\caption{\small Probability on Sobolev spaces}\label{fig:Sobolev-initial}
\vspace{0pt}
%\vspace{10pt}
\end{wrapfigure}

In this article, we will rethink the classical approximation problems by the kernel-based probability structures of the Sobolev spaces such as Figure~\ref{fig:Sobolev-initial}.
Theorem~\ref{t:Gauss-PDK-L} provides that the Sobolev space $\Hilbert^m(\Domain)=\Omega_m$ can be endowed with the probability measure $\PP_K$ induced by the positive definite kernel $K$.
Since the probability measure $\PP_K$ is placed on the Borel $\sigma$-algebra $\Borel\left(\Hilbert^m(\Domain)\right)=\Filter_m$,
the probability of the Sobolev spaces is dependent of the Sobolev norms such as
the probability is largest at the origin and the probability decreases to $0$ when the Sobolev norm tends to $\infty$.
These kernel-based probability structures are consistent with the common senses of the initial guess at $0$ with no information.
The probability will help us to measure the estimate values based on all feasible interpolating paths in the Sobolev spaces similar as the initial ideas in Section~\ref{sec:Initial}.

Let us look at the interpretive example in Figure~\ref{fig:Sobolev-initial} which can be thought as the generalization of the initial ideas in Figure~\ref{fig:1DExaInitial}.
The green, blue, and cyan lines represent the collections of the sample paths for various estimate values $v\in\RR$, that is,
\[
\Aset_L(v):=\left\{\omega\in\Hilbert^m(\Domain): L\omega=v\right\}.
\]
The red line represents all feasible interpolating paths,
that is,
\[
\Aset_{\vL}(\vf):=\left\{\omega\in\Hilbert^m(\Domain): \vL\omega=\vf\right\}.
\]
The black and pink squares, which can be roughly thought as the generalizations of the black and pink squares in Figure~\ref{fig:1DExaInitial}, represent the intersections $\Aset_L(v)\cap\Aset_{\vL}(\vf)$.
The yellow, green, cyan, and blue dashed circles
represent various ranges of the probability on $\Hilbert^m(\Domain)$.
Since the blue dashed circle is closed to the origin,
the probability shown in the blue dashed circle is larger than the others.
Then the best estimator $\hat{v}$ is given by the value $v_2$ because the probability of $v_2$ shown in the pink square is largest, for example, the probability of $v_2$ is larger than $v_1,v_3$ shown in the black squares.

Since $L$ and $\vL$ are bounded on $\Hilbert^m(\Domain)$, the sets $\Aset_L(v)$ and $\Aset_{\vL}(\vf)$ are closed in $\Hilbert^m(\Domain)$;
hence $\Aset_L(v),\Aset_{\vL}(\vf)\in\Borel\left(\Hilbert^m(\Domain)\right)$.
Here, we think that the interpolation $\Aset_{\vL}(\vf)$ has happened because the interpolating data $\vf$ have been given.
Then the estimate value $v$ can be measured by the probability of $\Aset_L(v)$ conditioned on the interpolation $\Aset_{\vL}(\vf)$.
According to the constructions of $LS$ and $\vL S$,
we have
\[
\Aset_L(v)=\left\{\omega\in\Omega_m: LS(\omega)=v\right\},\quad \Aset_{\vL}(\vf)=\left\{\omega\in\Omega_m:\vL S(\omega)=\vf\right\}.
\]
This indicates that the sets $\Aset_L(v)$ and $\Aset_{\vL}(\vf)$ can be equivalently transferred into
$LS=v$ and $\vL S=\vf$, respectively.
This shows that $LS$ and $\vL S$ can be used to compute the conditional probability $\PP_K\left(\Aset_L(v)|\Aset_{\vL}(\vf)\right)$.

Next, we show how to obtain the best estimator $\hat{v}$ by the techniques of statistical learning.
According to the \emph{maximum likelihood estimation} \cite[Section 6.4]{Stein1999},
the best estimator $\hat{v}$ is to maximize the conditional probability, that is,
the maximizer $\hat{v}$ of the optimization problem
\begin{equation}\label{eq:deterministic-max-cond-prob}
\max_{v\in\RR}\PP_K\left(\Aset_L(v)|\Aset_{\vL}(\vf)\right)
=\max_{v\in\RR}\PP_K\left(LS=v|\vL S=\vf\right).
\end{equation}
Since Corollary~\ref{c:Gauss-PDK-L-pdf} provides the conditional probability density function $p_{L|\vL}$ of $L S$ given $\vL S$,
the optimal solution $\hat{v}$ of the maximum problem~\eqref{eq:deterministic-max-cond-prob} can be obtained by another equivalent optimization problem of $p_{L|\vL}$, that is,
\begin{equation}\label{eq:deterministic-approx}
\hat{v}:=
\underset{v\in\RR}{\text{argmax}}~p_{L|\vL}\left(v|\vf\right)
=\mu_{L|\vL}\left(\vf\right)
=L\vk_{\vL}^T\vK_{\vL}^{\dag}\vf.
\end{equation}
Here, the best estimator $\hat{v}$ is called the \emph{kernel-based estimator} of the value $Lu$.

%////////////////////////////////////////////////////////////////////////////////////////////////////////////////////////
\begin{remark}
Generally speaking, the conditional probability density function does not mean the exact conditional probability.
In statistics, the mode of the probability density function indicates the maximum probability.
Thus, the optimization problems~\eqref{eq:deterministic-max-cond-prob} and~\eqref{eq:deterministic-approx} are equivalent
(see \cite{Stein1999,HastieTibshiraniFriedman2009}).
\end{remark}
%////////////////////////////////////////////////////////////////////////////////////////////////////////////////////////

In another ways,
by the \emph{Bayesian estimation} \cite[Section 4]{SchabackWendland2006},
the best estimator $\hat{v}$ can be also computed by the conditional expectation of $LS$ given the interpolation $\vL S=\vf$, that is,
\[
\hat{v}:=\Mean\left(LS|\vL S=\vf\right)=\int_{\RR}vp_{L|\vL}\left(v|\vf\right)\ud v
=\mu_{L|\vL}\left(\vf\right).
\]
Roughly, the Bayesian estimator can be thought as the averages of the estimate values $v$ based on the probability measures on the Sobolev spaces.
Since $LS$ and $\vL S$ have the normal distributions, the best estimator $\hat{v}$ is the same for the both maximum-likelihood and Bayesian methods.

Moreover, the best estimator (kernel-based estimator) $\hat{v}$
can be rewritten as the similar forms of the Hermite-Birkhoff interpolation.
So, we will construct a function $\hat{u}_{\vL}\in\Hilbert^m(\Domain)$ to compute the best estimator $\hat{v}$ by the bounded linear functional $L$ such as $\hat{v}=L\hat{u}_{\vL}$.
Equation~\eqref{eq:deterministic-approx} assures that $\hat{u}_{\vL}$ can be written as a linear combination of the kernel basis
$\vk_{\vL}$, that is,
\begin{equation}\label{eq:kernel-approx-model-1}
\hat{u}_{\vL}(\vx):=\vk_{\vL}(\vx)^T\vc,\quad\text{for }\vx\in\Domain,
\end{equation}
where the coefficients $\vc$ are solved by a linear system
\begin{equation}\label{eq:kernel-approx-model-2}
\vK_{\vL}\vc=\vf.
\end{equation}
We further find that the approximate function $\hat{u}_{\vL}$ is independent of the bounded linear functional $L$.
Since the point evaluation function $\delta_{\vx}$ is a bounded linear functional on $\Hilbert^m(\Domain)$ for $m>d/2$,
the function values $u(\vx)=\delta_{\vx}\left(u\right)$ can be approximated by the estimators $\delta_{\vx}\left(\hat{u}_{\vL}\right)=\hat{u}_{\vL}(\vx)$.
Thus, we say $\hat{u}_{\vL}$ a \emph{kernel-based approximate function} of the target function $u$.

Finally, we show the convergence of the kernel-based estimators.
Suppose that we have the countable data information $f_1:=L_1u,\ldots,f_n:=L_nu,\ldots$ such that
\[
\Aset_{\vL_{\infty}}\left(\vf_{\infty}\right)=
\left\{\omega\in\Hilbert^m(\Domain): \vL_{\infty}\omega=\vf_{\infty}\right\}
=\left\{u\right\},
\]
where $\vL_{\infty}:=\left(L_1,\cdots,L_n,\cdots\right)^T$
and $\vf_{\infty}:=\left(f_1,\cdots,f_n,\cdots\right)^T$.
For example, the operator $\vL_{\infty}$ is composed of the point evaluation functions $\delta_{\vx_n}$, where the data points $\left\{\vx_n\right\}_{n=1}^{\infty}$
is dense in the domain $\Domain$.
Here, since $\Cont(\Domain)$ is imbedded into $\Hilbert^m(\Domain)$,
there is a unique function $u$ to interpolate all the given data.
Let $\vL_{n}:=\left(L_1,\cdots,L_n\right)^T$
and $\vf_{n}:=\left(f_1,\cdots,f_n\right)^T$ for all $n\in\NN$.
Then we can obtain the kernel-based estimator $L\hat{u}_{\vL_n}$ for the given data $\vf_n$ same as in Equation~\eqref{eq:deterministic-approx}.
Since
\[
\Aset_{\vL_1}\left(\vf_1\right)\supseteq\cdots\supseteq\Aset_{\vL_n}\left(\vf_n\right)
\supseteq\cdots\supseteq\overset{\infty}{\underset{n=1}{\cap}}\Aset_{\vL_n}\left(\vf_n\right)=\Aset_{\vL_{\infty}}\left(\vf_{\infty}\right),
\]
we have
\[
\lim_{n\to\infty}L\hat{u}_{\vL_n}
=
\lim_{n\to\infty}\Mean\left(LS|\Aset_{\vL_n}\left(\vf_n\right)\right)
=
\Mean\left(LS|\Aset_{\vL_{\infty}}\left(\vf_{\infty}\right)\right)
=Lu.
\]
In particular $\hat{u}_{\vL_n}(\vx)=\delta_{\vx}\left(\hat{u}_{\vL_n}\right)\to\delta_{\vx}\left(u\right)=u(\vx)$ when $n\to\infty$ for any $\vx\in\Domain$.
This assures that the kernel-based approximate function $\hat{u}_{\vL_n}$ is also convergent to the target function $u$ when $n\to\infty$.

\emph{Comments:}
In approximation theory, we mainly focus on the constructions of the globally best interpolants.
In statistical learning, we usually learn the locally random variables by another correlated random elements.
For example, the meshfree approximation gives the globally optimal solutions while the kriging interpolation provides the locally optimal estimators.
In this article, we try to combine the knowledge of meshfree approximation and kriging interpolation in one theoretical structure such that we can obtain the best estimators both solved by the locally random variables $LS,\vL S$ and supported by the globally interpolating paths $\Aset_L(v),\Aset_{\vL}(\vf)$.
The meaning of the best is dependent of the kernel-based probability structures of the Sobolev spaces here.
This new idea also improves the meshfree approximation and the kriging interpolation as follows:
\begin{itemize}

\item In meshfree approximation, we usually suppose that the matrix $\vK_{\vL}$ is well-condition such that the Lagrangian basis $\ve_{\vL}:=\vK_{\vL}^{-1}\vk_{\vL}$ is well-defined. Unfortunately, the matrix $\vK_{\vL}$ could be nonsingular or ill-condition in the practical applications. In numerical analysis, we can still solve the ill-condition problems by the least-square techniques such as the coefficients are given by $\vc:=\vK_{\vL}^{\dag}\vf$ (see the pseudoinverse minimal solutions in \cite[Theorem~5.4.2]{KincaidWard2002}). But, the exact geometrical meanings of the least-square solutions are still unclear for the interpolation problems. Here, the probability measure $\PP_K$ provides another way to explain the least squares and the generalized Lagrangian basis $\ve_{\vL}:=\vK_{\vL}^{\dag}\vk_{\vL}$.
    Even though $\hat{u}_{\vL}=\ve_{\vL}^T\vf$ may not belong to $\Aset_{\vL}(\vf)$ or satisfy the interpolation conditions, the least-square interpolant $\hat{u}_{\vL}$ can be still thought as the best adjacent element of $\Aset_{\vL}(\vf)$.

\item In kriging interpolation, we can only consider the interpolating or spatial data related to the point evaluation functions $\delta_{\vx}$. Here, the random variables can be supported by the interpolating paths in the Sobolev spaces such that the kriging interpolation is still well-defined by another operators, for example, the differential operators $\delta_{\vx}\circ\nabla$ and $\delta_{\vx}\circ\Delta$. This shows that the meshfree approximation implies the generalized kriging interpolation.

\end{itemize}

In the following sections, we will continue to extend the kernel-based methods for the deterministic problems to the stochastic problems by the same manners.

%---------------------------------------------------------------------------------------------------------------------
\section{Stochastic Data Interpolations}\label{sec:Stochastic-Data-Approx}
%---------------------------------------------------------------------------------------------------------------------

In this section, we will extend the meshfree approximation~\cite{Wendland2005,Fasshauer2007} for the deterministic data to the stochastic data.
Hence, let us look at the random data values $f_1,\ldots,f_N$ interpolated at the distinct data points
\[
X:=\left\{\vx_1,\ldots,\vx_N\right\}\subseteq\Domain\subseteq\Rd.
\]
In Section~\ref{sec:ker-app}, the deterministic data are obtained by a deterministic function.
Here, we suppose that the stochastic data are simulated by a stochastic model
\begin{equation}\label{eq:stoch-model}
u(\vx)=\Gamma\left(\vx,\vartheta_{\vx}\right),\quad\text{for }\vx\in\Domain,
\end{equation}
where $\Gamma$ is a deterministic function and $\xi$ is a Gaussian field with the mean $\vzero$ and the known covariance kernel $\Phi:\Domain\times\Domain\to\RR$ on a probability space $\left(\Omega,\Filter,\PP\right)$.
Then
\[
f_1:=u\left(\vx_1\right)=\Gamma\left(\vx_1,\vartheta_{\vx_1}\right),\ldots,
f_N:=u\left(\vx_N\right)=\Gamma\left(\vx_N,\vartheta_{\vx_N}\right),
\]
are the random variables defined on the probability space $\left(\Omega,\Filter,\PP\right)$.
By the Monte Carlo methods in \cite{Glasserman2004}, we can easily simulate
the multivariate normal random vector
\[
\vvartheta_X:=\left(\vartheta_{\vx_1},\cdots,\vartheta_{\vx_N}\right)^T\sim\Normal\left(\vzero,\vTheta_X\right),
\]
where the covariance matrix $\vTheta_X:=\left(\Phi(\vx_j,\vx_k)\right)_{j,k=1}^{N,N}$.
Then we can obtain the probability distributions of the random vector
\[
\vf:=\left(f_1,\cdots,f_N\right)^T.
\]

Same as in Section~\ref{sec:Gauss-PDK}, we suppose that the domain $\Domain$ is regular and compact and the stochastic model $u\in\Hilbert^m(\Domain)$ for $m>d/2$.

%////////////////////////////////////////////////////////////////////////////////////////////////////////////////////////
\begin{remark}
Some papers may require $u\in\Cont^m(\Domain)\text{ or }\Hilbert^m(\Domain)$ almost surely.
Usually, the smoothness of $u$ can be guaranteed by the smoothness of $\Gamma$ and $\Phi$, for example,
if $\Gamma\in\Cont^{2m}\left(\Domain\times\RR\right)$ and $\Phi\in\Cont^{2m,1}\left(\Domain\times\Domain\right)$ then $u\in\Cont^m(\Domain)\subseteq\Hilbert^m(\Domain)$ almost surely. For convenience, we can ignore the non-smooth or non-Sobolev-normed sample paths of the stochastic model $u$ in this section.
\end{remark}
%////////////////////////////////////////////////////////////////////////////////////////////////////////////////////////

Given a derivative operator
\begin{equation}\label{eq:L-derivative}
L:=\delta_{\vx}\circ D^{\valpha},\quad\text{for }\vx\in\Domain\text{ and }\abs{\valpha}<m-d/2,
\end{equation}
we try to compute the probability distributions of the random variable $Lu$.
According to the Sobolev imbedding theorem, the derivative operator $L$ is a bounded linear functional on the Sobolev space $\Hilbert^m(\Domain)$; hence $Lu$ is well-defined on the probability space $\left(\Omega,\Filter,\PP\right)$.
However, it may be difficult to simulate the probability distribution directly when $\Gamma$ is a nonlinear function.

Hence, we need to use the easily simulated stochastic data to approximate the probability distributions of $Lu$.
In stochastic analysis, the initial conditions of the stochastic ordinary differential equations can be deterministic or stochastic such as
the existence and uniqueness theorem for the stochastic ordinary differential equations \cite[Theorem~5.2.1]{Oksendal2003}.
This inspires us to extend the kernel-based approximation in Section~\ref{sec:ker-app} to construct the best estimator (kernel-based estimator)
of $Lu$ in the following steps.

Firstly, we choose a positive definite kernel $K\in\Cont^{2m}\left(\Domain\times\Domain\right)$.
Let the vector operator
\[
\vL:=\left(L_1,\cdots,L_N\right)^T,
\]
be composed of the point evaluation functions
\[
L_1:=\delta_{\vx_1},\ldots,L_N:=\delta_{\vx_N}.
\]
Obviously, all point evaluation functions are bounded on $\Hilbert^m(\Domain)$ because $m>d/2$.
Thus, by Theorem~\ref{t:Gauss-PDK-L} and Corollary~\ref{c:Gauss-PDK-multi-L}, we can construct the multivariate normal random variables $LS$ and $\vL S$ on $\Omega_m=\Hilbert^m(\Domain)$ under the probability measure $\PP_K$ induced by the kernel $K$.
This indicates that $LS$ and $\vL S$ are correlated on the probability space $\left(\Omega_m,\Filter_m,\PP_K\right)$ given in Theorem~\ref{t:Gauss-PDK-L}. But $Lu$ and $\vf$ are correlated on another probability space $\left(\Omega,\Filter,\PP\right)$.

Therefore, we need to combine the different probability spaces $\left(\Omega_m,\Filter_m,\PP_K\right)$ and $\left(\Omega,\Filter,\PP\right)$
into one probability space such that we can discuss $LS,\vL S$ and $Lu,\vf$ together.
Then we define a tensor product probability space
\begin{equation}\label{eq:product-prob}
\hat{\Omega}:=\Omega_m\times\Omega,\quad \hat{\Filter}:=\Filter_m\otimes\Filter,\quad \hat{\PP}:=\PP_K\times\PP,
\end{equation}
such that all original random variables on $\left(\Omega_m,\Filter_m,\PP_K\right)$ and $\left(\Omega,\Filter,\PP\right)$ can be extended naturally onto $\big(\hat{\Omega},\hat{\Filter},\hat{\PP}\big)$. To be more precisely, the extensions of the original random variables $V_1:\Omega_{m}\to\RR$ and
$V_2:\Omega\to\RR$ defined by
\[
V_1\left(\omega_1\times\omega_2\right):=V_1\left(\omega_1\right),\quad
V_2\left(\omega_1\times\omega_2\right):=V_2\left(\omega_2\right),\quad
\text{for }\omega_1\times\omega_2\in\hat{\Omega},
\]
preserve the original probability distributions and the extensions of $V_1$ and $V_2$ are independent on $\big(\hat{\Omega},\hat{\Filter},\hat{\PP}\big)$ because the two probability spaces
$\left(\Omega_m,\Filter_m,\PP_K\right)$ and $\left(\Omega,\Filter,\PP\right)$ are independent.
Thus, the extensions of $LS$ and $\vL S$ keep the same probability distributions on $\big(\hat{\Omega},\hat{\Filter},\hat{\PP}\big)$.
This indicates that the conditional probability density function of the extension of $LS$ given $\vL S$ is still
equal to $p_{L|\vL}$ given in Corollary~\ref{c:Gauss-PDK-L-pdf}.
Moreover, the extensions of $\left(LS,\vL S\right)$ and $\vf$ are independent.

{\bf Kernel-based Estimators and Kernel-based Approximate Functions:}
Obviously, we find that
\[
\vk_{\vL}(\vx)=\vk_X(\vx)=\left(K\left(\vx,\vx_1\right),\cdots,K\left(\vx,\vx_N\right)\right)^T,\quad
\vK_{\vL}=\vK_{X}=\left(K(\vx_j,\vx_k)\right)_{j,k=1}^{N,N}.
\]
Therefore, same as in Equations~(\ref{eq:deterministic-max-cond-prob}-\ref{eq:deterministic-approx}),
we can obtain the best estimator $\hat{v}$ of $Lu$ by the maximum likelihood estimation methods, that is,
\begin{equation}\label{eq:kernel-approx-stoch-model-0}
\hat{v}:=
\underset{v\in\RR}{\text{argmax}}~\hat{\PP}\left(LS=v|\vL S=\vf\right)
=\underset{v\in\RR}{\text{argmax}}~p_{L|\vL}\left(v|\vf\right)
=L\vk_{X}^T\vK_{X}^{-1}\vf.
\end{equation}
Same as in Section~\ref{sec:ker-app}, the best estimator $\hat{v}$ is called the kernel-based estimator of $Lu$.
Moreover, we can create the following algorithm to produce the thousands $p$ samples of $\hat{v}$ to approximate the probability distributions of $Lu$, that is,
\[
\tag{A1}
\begin{split}
&\text{Initialize:}\\
&\ve^T:=L\vk_{X}^T\vK_{X}^{-1}=D^{\valpha}\vk_{X}(\vx)^T\vK_{X}^{-1},\\
&\text{Repeat }i=1,\ldots,p\\
&\quad{}\text{Simulate }\vvartheta_X^{(i)}:=\left(\vartheta_{\vx_1}^{(i)},\cdots,\vartheta_{\vx_N}^{(i)}\right)^T\sim\Normal\left(\vzero,\vTheta_X\right),\\
&\quad{}\vf^{(i)}:=\left(\Gamma\big(\vx_1,\vartheta_{\vx_1}^{(i)}\big),\cdots,\Gamma\big(\vx_N,\vartheta_{\vx_N}^{(i)}\big)\right)^T,\\
&\quad{}\hat{v}^{(i)}:=\ve^T\vf^{(i)}.
\end{split}
\]
Here $\vvartheta^{(1)},\ldots,\vvartheta^{(p)}$ can be seen as the simulated duplications of the multivariate normal random vector $\vvartheta$ by the Monte Carlo methods.
For example, the mean and variance of $Lu$ can be approximated by
\[
\Mean\left(Lu\right)\approx\hat{\mu}:=\frac{1}{p}\sum_{i=1}^p\hat{v}^{(i)},
\quad
\Var\left(Lu\right)\approx\hat{\sigma}^2:=
\frac{1}{p}\sum_{i=1}^p\left(\hat{v}^{(i)}-\hat{\mu}\right)^2.
\]

Same as in Equations~(\ref{eq:kernel-approx-model-1}-\ref{eq:kernel-approx-model-2}), we can represent the best estimator $\hat{v}$ in Equation~\eqref{eq:kernel-approx-stoch-model-0} by the kernel-based approximate function $\hat{u}_X$ with the derivative operator $L$, that is,
$\hat{v}=L\hat{u}_X$.
Since $\hat{v}=L\vk_{X}^T\vK_{X}^{-1}\vf$, the kernel-based approximate function $\hat{u}_X$ can be written as a linear combination of the kernel basis $K\left(\cdot,\vx_1\right),\ldots,K\left(\cdot,\vx_N\right)$ such as
\begin{equation}\label{eq:kernel-approx-stoch-model-1}
\hat{u}_X(\vx):=\sum_{k=1}^Nc_kK\left(\vx,\vx_k\right),\quad
\text{for }\vx\in\Domain,
\end{equation}
and the random coefficients $\vc:=\left(c_1,\cdots,c_N\right)^T$ are solved by a well-posed random linear system
\begin{equation}\label{eq:kernel-approx-stoch-model-2}
\vK_X\vc=\vf.
\end{equation}
It is clear that $\hat{u}_X\in\Hilbert^m(\Domain)$ and $\vL\hat{u}_X=\vf$.
But, since the random vector $\vf$ may not be normal, the random function $\hat{u}_X$ may not be Gaussian.

In the following, we continue to study with the random parts of the kernel-based estimator $L\hat{u}_X$.
It is obvious that the random parts of $L\hat{u}_X$ is only dependent of the random vector $\vf$.
In probability theory, we can transfer $\vf$ equivalently onto a finite-dimensional
probability space (see \cite[Section~1.4]{Shiryaev1996}). To be more precise, we can view
$\vf$ as a random vector placed on the finite-dimensional probability space $\big(\RR^N,\Borel(\RR^N),\measure_{\vf}\big)$, where the probability measure $\measure_{\vf}$ is introduced by the probability density function $p_{\vf}$ of the random vector $\vf$, that is, $\measure_{\vf}(\ud\vv):=p_{\vf}(\vv)\ud\vv$.
Thus, the kernel-based estimator $L\hat{u}_X$ has the same probability distributions on all probability spaces $\big(\Omega,\Filter,\PP\big)$, $\big(\hat{\Omega},\hat{\Filter},\hat{\PP}\big)$, and
$\big(\RR^N,\Borel(\RR^N),\measure_{\vf}\big)$.

{\bf Error Analysis:}
Finally, we propose to verify the convergence of the kernel-based estimator $L\hat{u}_X$ in probability.

It is well known that the convergence of kriging interpolation is dependent of standard deviations and
the convergence of meshfree approximation is dependent of power functions.
In the beginning,
we show that the standard deviation $\sigma_{L|\vL}$
of the conditional probability density function $p_{L|\vL}$ is equal to the \emph{power functions} $\Power_X^{\valpha}(\vx)$.
Let a quadratic form $\Quadratic:\RR^N\to\RR$ be
\[
\Quadratic(\vv):=
D^{\valpha}_{\vz_1}D^{\valpha}_{\vz_2}K\left(\vz_1,\vz_2\right)|_{\vz_1=\vz_2=\vx}-2D^{\valpha}\vk_X(\vx)^T\vv+\vv^T\vK_{X}\vv,
\quad\text{for }\vv\in\RR^N.
\]
By \cite[Definition~11.2]{Wendland2005}, the power function $\Power_X^{\valpha}(\vx)$ is the minimum of $\Quadratic$, that is,
\[
\Power_X^{\valpha}(\vx):=\min_{\vv\in\RR^N}\sqrt{\Quadratic\left(\vv\right)}.
\]
Comparing with Equation~\eqref{eq:standard-dev}, we have
\begin{equation}\label{eq:stoch-model-error-1}
\Power_X^{\valpha}(\vx)=\sqrt{\Quadratic\left(\vK_{X}^{-1}D^{\valpha}\vk_X(\vx)\right)}=\sigma_{L|\vL}.
\end{equation}
Moreover, since $K\in\Cont^{2m}(\Domain\times\Domain)$,
\cite[Theroem~11.13]{Wendland2005} (errors estimates for power functions) provides that
\begin{equation}\label{eq:stoch-model-error-2}
\Power_X^{\valpha}(\vx)=\Order\left(h_X^{m-\abs{\valpha}}\right),
\quad
\text{when }h_X\text{ is small enough}.
\end{equation}
Here $h_X$ is the fill distance of the data points $X$ for the domain $\Domain$, that is,
\[
h_X:=\sup_{\vx\in\Domain}\min_{k=1,\ldots,N}\norm{\vx-\vx_k}_2,
\]
or the fill distance $h_X$ denotes the radius of the largest ball in the domain $\Domain$ and without any data points $X$.
Combining Equations~\eqref{eq:stoch-model-error-1} and~\eqref{eq:stoch-model-error-2}, we have
\begin{equation}\label{eq:stoch-model-error-3}
\sigma_{L|\vL}=\Order\left(h_X^{m-\abs{\valpha}}\right),
\quad
\text{when }h_X\text{ is small enough}.
\end{equation}

%////////////////////////////////////////////////////////////////////////////////////////////////////////////////////////
\begin{remark}
According to the smoothness of $K$,
Equation~\eqref{eq:stoch-model-error-2} can be checked by the Taylor expansion of $K$.
More details of the upper bounds of the power functions can be found in \cite[Chapter~11]{Wendland2005} and \cite[Chapter~14]{Fasshauer2007}.
\end{remark}
%////////////////////////////////////////////////////////////////////////////////////////////////////////////////////////

To investigate the error $\epsilon>0$, we estimate the probability of
$\abs{Lu-L\hat{u}_X}<\epsilon$ or $\abs{Lu-L\hat{u}_X}\geq\epsilon$ firstly. In probability theory, we call that $L\hat{u}_X$ converges to $Lu$ in probability if $\hat{\PP}\left(\abs{Lu-L\hat{u}_X}<\epsilon\right)\to1$ or
$\hat{\PP}\left(\abs{Lu-L\hat{u}_X}\geq\epsilon\right)\to0$.
Here, since $u\left(\omega_2\right)\in\Hilbert^m(\Domain)$ for any $\omega_2\in\Omega$, we have $u\left(\cdot,\omega_2\right)\in\hat{\Omega}$; hence
the value $Lu(\omega_2)$ is dependent of the sample $u\left(\cdot,\omega_2\right)\in\hat{\Omega}$.
This indicates that $\abs{Lu-L\hat{u}_X}<\epsilon$ or $\abs{Lu-L\hat{u}_X}\geq\epsilon$ can be viewed as an event on the probability space $\big(\hat{\Omega},\hat{\Filter},\hat{\PP}\big)$.
For the proofs of the convergence, we will compute the probability $\hat{\PP}\left(\abs{Lu-L\hat{u}_X}\geq\epsilon\right)$ as follows:

%////////////////////////////////////////////////////////////////////////////////////////////////////////////////////////
\begin{lemma}\label{l:stoch-model-error}
Suppose that $L\hat{u}_X$ is the kernel-based estimator of $Lu$ in Equations~\eqref{eq:kernel-approx-stoch-model-0} or (\ref{eq:kernel-approx-stoch-model-1}-\ref{eq:kernel-approx-stoch-model-2}). Then we have
\begin{equation}\label{eq:stoch-model-error-erfc}
\hat{\PP}\left(\abs{Lu-L\hat{u}_X}\geq\epsilon\right)=\erfc\left(\frac{\epsilon}{\sqrt{2}\sigma_{L|\vL}}\right),\quad
\text{for any }\epsilon>0,
\end{equation}
where the variance $\sigma_{L|\vL}^2=L_{\vx}L_{\vy}K(\vx,\vy)-L\vk_{X}^T\vK_{X}^{-1}L\vk_{X}$ is same as in Corollary~\ref{c:Gauss-PDK-L-pdf} and $\erfc$ is the complementary error function, that is, $\erfc(z):=2\pi^{-1/2}\int_z^{\infty}e^{-t^2}\ud t$.
\end{lemma}
%////////////////////////////////////////////////////////////////////////////////////////////////////////////////////////
\begin{proof}
Let the set
\[
\Eset:=\left\{\omega_1\times\omega_2\in\hat{\Omega}: \abs{L\omega_1-L\hat{u}_X\left(\omega_2\right)}\geq\epsilon\text{ subject to }\vL\omega_1=\vf\left(\omega_2\right)\right\}.
\]
The main idea of the proof is to use the probability $\hat{\PP}\left(\Eset\right)$ to estimate the probability $\hat{\PP}\left(\abs{Lu-L\hat{u}_X}\geq\epsilon\right)$. Generally speaking, we will evaluate the probability of the kernel-based estimator $L\hat{u}_X$ against the error $\epsilon$ when the interpolations are true.

The constructions of the random variables $LS$ and $\vL S$ provides that
\[
\Eset=\left\{\omega_1\times\omega_2\in\hat{\Omega}: \abs{LS\left(\omega_1\right)-L\hat{u}_X\left(\omega_2\right)}<\epsilon\text{ subject to  }\vL S\left(\omega_1\right)=\vf\left(\omega_2\right)\right\}.
\]
Thus, combing with the independence of $(LS,\vL S)$ and $(L\hat{u}_X,\vf)$, we have
\begin{align*}
&~\hat{\PP}\left(\Eset\right)=
\int_{\RR^N}\hat{\PP}\left(\abs{LS-L\hat{u}_X(\vv)}\geq\epsilon|\vL S=\vv\right)\measure_{\vf}\left(\ud\vv\right)\\
=&\int_{\RR^N}\int_{\abs{v-L\hat{u}_X(\vv)}\geq\epsilon}p_{L|\vL}(v|\vv)p_{\vf}(\vv)\ud v\ud\vv
=\erfc\left(\frac{\epsilon}{\sqrt{2}\sigma_{L|\vL}}\right).
\end{align*}

Moreover,
since $u\left(\cdot,\omega_2\right)\in\hat{\Omega}$  and $\vL u\left(\omega_2\right)=\vf\left(\omega_2\right)$ for all $\omega_2\in\Omega$,
we can assure that
$\abs{Lu\left(\omega_2\right)-L\hat{u}_X\left(\omega_2\right)}\geq\epsilon$ if and only if $u\left(\cdot,\omega_2\right)\in\Eset$.
Therefore,
\[
\hat{\PP}\left(\abs{Lu-L\hat{u}_X}\geq\epsilon\right)
=\hat{\PP}\left(\Eset\right)
=\erfc\left(\frac{\epsilon}{\sqrt{2}\sigma_{L|\vL}}\right).
\]
\end{proof}

Different from kriging interpolation, we will obtain the convergence of the kernel-based estimators by the techniques of meshfree approximation.
Combining Equations~(\ref{eq:stoch-model-error-3}-\ref{eq:stoch-model-error-erfc}), we have
\begin{equation}\label{eq:stoch-model-error-fill-distance}
\hat{\PP}\left(\abs{Lu-L\hat{u}_X}\geq\epsilon\right)=\Order\left(\frac{h_X^{m-\abs{\valpha}}}{\epsilon}\right),\quad
\text{when }h_X\text{ is small enough};
\end{equation}
hence
\[
\lim_{h_X\to0}\hat{\PP}\left(\abs{Lu-L\hat{u}_X}\geq\epsilon\right)=0.
\]
Therefore, we can conclude that:

%////////////////////////////////////////////////////////////////////////////////////////////////////////////////////////
\begin{proposition}\label{p:stoch-model-error}
Suppose that $L\hat{u}_X$ is the kernel-based estimator of $Lu$ in Equations~\eqref{eq:kernel-approx-stoch-model-0} or (\ref{eq:kernel-approx-stoch-model-1}-\ref{eq:kernel-approx-stoch-model-2}). Then $L\hat{u}_X$ converges to $Lu$ in probability when the fill distance $h_X\to0$.
\end{proposition}
%////////////////////////////////////////////////////////////////////////////////////////////////////////////////////////

%////////////////////////////////////////////////////////////////////////////////////////////////////////////////////////
\begin{remark}
Obviously $Lu$ and $L\hat{u}_X$ are well-posed on the both probability spaces $\left(\Omega,\Filter,\PP\right)$ and $\big(\hat{\Omega},\hat{\Filter},\hat{\PP}\big)$.
Proposition~\ref{p:stoch-model-error} provides the convergence of $L\hat{u}_X$ under $\hat{\PP}$. Since $\hat{\PP}$ is the product probability measure composed of $\PP_K$ and $\PP$, the convergence of $L\hat{u}_X$ is also well-posed under $\PP$.
But, this does not imply that the convergence of $L\hat{u}_X$ is exactly true because $\PP_K$ vanishes all non-smooth paths.
Roughly, we say that $L\hat{u}_X$ converges weakly to $Lu$.
More details of the various kinds of the convergence of the sequences of random variables are mentioned in probability theory (see \cite[Section~2.10]{Shiryaev1996}).
\end{remark}
%////////////////////////////////////////////////////////////////////////////////////////////////////////////////////////

Since the convergence in probability implies the convergence in distribution by
\cite[Theorem~2.2]{Shiryaev1996}, we have:

%////////////////////////////////////////////////////////////////////////////////////////////////////////////////////////
\begin{corollary}\label{c:stoch-model-error}
Suppose that $L\hat{u}_X$ is the kernel-based estimator of $Lu$ in Equations~\eqref{eq:kernel-approx-stoch-model-0} or (\ref{eq:kernel-approx-stoch-model-1}-\ref{eq:kernel-approx-stoch-model-2}).
Let $g\in\Cont(\RR)$ be a bounded function.
Then
\[
\lim_{h_X\to0}\Mean\left(g\left(L\hat{u}_X\right)\right)
=\Mean\left(g\left(Lu\right)\right).
\]
In particular, if $\abs{Lu}\leq C_L$ for a deterministic constant $C_L>0$, then
\[
\lim_{h_X\to0}\Mean\left(L\hat{u}_X\right)
=\Mean\left(Lu\right),\quad
\lim_{h_X\to0}\Var\left(L\hat{u}_X\right)
=\Var\left(Lu\right).
\]
\end{corollary}
%////////////////////////////////////////////////////////////////////////////////////////////////////////////////////////

Combing with \cite[Theorem~3.2]{Shiryaev1996} (weak convergence of probability distributions),
Corollary~\ref{c:stoch-model-error} also assures that
the cumulative distribution function of $L\hat{u}_X$ converges to the cumulative distribution function of $Lu$ when $h_X\to0$.
This shows that the probability distributions of $Lu$ can be approximated by the probability distributions of $L\hat{u}_X$.

\emph{Comments:}
In this section, we show that the meshfree approximation for the deterministic interpolations can be extended to the stochastic interpolations.
Typically, we find that the kernel-based approximate function $\hat{u}_X$ given in Equations~(\ref{eq:kernel-approx-stoch-model-1}-\ref{eq:kernel-approx-stoch-model-2}) is consistent with the classical formats of meshfree approximation, that is, a linear combination of the kernel basis $K\left(\cdot,\vx_1\right),\ldots,K\left(\cdot,\vx_N\right)$ (see~\cite{Wendland2005,Fasshauer2007}). In approximation theory, the kernel-based approximate function is solved to minimize the reproducing norms globally, that is,
\[
\min_{f\in\Hilbert_K(\Domain)}\norm{f}_{\Hilbert_K(\Domain)}\text{ subject to }
f\left(\vx_1\right)=f_1,\ldots,f\left(\vx_N\right)=f_N,
\]
(see \cite[Theorem~13.2]{Wendland2005}).
In statistical learning, the kernel-based approximate function is obtained by the maximizing probabilities locally such as Equation~\eqref{eq:kernel-approx-stoch-model-0}.
Roughly speaking, the kernel-based methods gather the global and local solutions in one theoretical approach.

Since we show that the standard deviations and the power functions are the same (see Equation~\eqref{eq:stoch-model-error-1}),
we find that the error estimates of kriging interpolation can be analyzed by the techniques of meshfree approximation.
The paper~\cite{ScheuererSchabackSchlather2012} firstly illustrates the equivalent concepts for $\abs{\valpha}=0$.
In this section, we verify that it is true for all feasible $\valpha$.
This let us obtain the convergent rates of the random variables by the fill distances shown in meshfree approximation.
The fill distance is a common sense in numerical analysis. But, the fill distance is a novel concept in statistics. Many current researches of statistical learning focus on the number $N$ of the data information.
This gives a new way to design the optimal estimators of the stochastic models.

%---------------------------------------------------------------------------------------------------------------------
\section{Elliptic Stochastic Partial Differential Equations}\label{sec:ell-SPDE}
%---------------------------------------------------------------------------------------------------------------------

In this section, we will solve the elliptic SPDEs by the kernel-based methods.
Same as in Section~\ref{sec:Gauss-PDK}, we let $\Domain\subseteq\Rd$ be a regular and compact domain.
Then the boundary $\partial\Domain$ of $\Domain$ is also regular and compact.
Now we look at a SPDE
\begin{equation}\label{eq:ellptic-SPDE}
\begin{cases}
Pu=\Gamma\left(\cdot,\vartheta\right),&\text{in }\Domain,\\
Bu=g,&\text{on }\partial\Domain,
\end{cases}
\end{equation}
where $\Gamma$ and $g$ are the deterministic functions and $\vartheta$ is a Gaussian field with the mean $\vzero$ and the known covariance kernel $\Phi:\Domain\times\Domain\to\RR$ on a probability space $\left(\Omega,\Filter,\PP\right)$.
For comparing the kernel-based approximate solutions of the deterministic PDEs in \cite[Section~16.3]{Wendland2005} and \cite[Chapter~38]{Fasshauer2007} easily, we only discuss the
uniformly elliptic differential operator of the 2nd order with the constant coefficients and the Dirichlet's boundary conditions, that is,
\[
P:=-\nabla^T\vA\nabla+\vb^T\nabla+c,\quad B:=I|_{\partial\Domain},
\]
where $\nabla$ is a gradient operator, $I$ is an identity operator, $\vA\in\RR^{d\times d}$ is a strictly positive definite matrix, $\vb\in\RR^d$, and $c\in\RR$. We further suppose that the solution $u\in\Hilbert^m(\Domain)$ for $m>2+d/2$.

Before the constructions of the kernel-based approximate solutions of the SPDE~\eqref{eq:ellptic-SPDE},
we firstly illustrate the symbols in this section. Let $K\in\Cont^{2m,1}\left(\Domain\times\Domain\right)$ be a positive definite kernel and
\[
L:=\delta_{\vx},\quad\text{for }\vx\in\Domain.
\]
Since there are two regions $\Domain$ and $\partial\Domain$, we will choose the distinct data points in the domain $\Domain$ and the boundary $\partial\Domain$, respectively, that is,
\[
X:=\left\{\vx_1,\ldots,\vx_N\right\}\subseteq\Domain,\quad
Z:=\left\{\vz_1,\ldots,\vz_M\right\}\subseteq\partial\Domain.
\]
Then, by the Sobolev imbedding theorem and the boundary trace imbedding theorem \cite[Theorem~5.36]{AdamsFournier2003},
the linear vector operator
\[
\vL:=\left(\delta_{\vx_1}\circ P,\cdots,\delta_{\vx_N}\circ P,\delta_{\vz_1}\circ B,\cdots,\delta_{\vz_M}\circ B\right)^T,
\]
is bounded on $\Hilbert^m(\Domain)$. This indicates that the kernel basis $\vk_{\vL}$ and the covariance (interpolating) matrix $\vK_{\vL}$ can be written as
\[
\vk_{\vL}(\vx)=\left(P_{\vy}K\left(\vx,\vx_1\right),\cdots,P_{\vy}K\left(\vx,\vx_N\right),
B_{\vy}K\left(\vx,\vz_1\right),\cdots,B_{\vy}K\left(\vx,\vz_M\right)
\right)^T,
\]
and
\[
\vK_{\vL}=
\begin{pmatrix}
\left(P_{\vx}P_{\vy}K(\vx_j,\vx_k)\right)_{j,k=1}^{N,N}&\left(P_{\vx}B_{\vy}K(\vx_j,\vz_k)\right)_{j,k=1}^{N,M}\\
\left(B_{\vx}P_{\vy}K(\vz_j,\vx_k)\right)_{j,k=1}^{M,N}&\left(B_{\vx}B_{\vy}K(\vz_j,\vz_k)\right)_{j,k=1}^{M,M}
\end{pmatrix}.
\]
Since all coefficients of the differential operator $P$ are constant,
\cite[Corollary~16.12]{Wendland2005} assures that $\vK_{\vL}$ is a strictly positive definite matrix.

Moreover, we can obtain the stochastic data simulated by the right-hand sides of the SPDE~\eqref{eq:ellptic-SPDE}.
To be more precise,
we can simulate the Gaussian field $\vartheta$ at $X$ by the Monte Carlo methods, that is,
\[
\vvartheta_X:=\left(\vartheta_{\vx_1},\cdots,\vartheta_{\vx_N}\right)^T\sim\Normal\left(\vzero,\vTheta_X\right).
\]
Notes that the stochastic data
\[
\gamma_1:=\Gamma\left(\vx_1,\vartheta_{\vx_1}\right),\ldots,\gamma_N:=\Gamma\left(\vx_N,\vartheta_{\vx_N}\right),~
g_1:=g(\vz_1),\ldots,g_M:=g(\vz_M).
\]
For convenience, we let
\[
\vf:=\left(\gamma_1,\cdots,\gamma_N,g_1,\cdots,g_M\right)^T.
\]

{\bf Kernel-based Approximate Solutions:}
Next, by the same manners of Equations~(\ref{eq:deterministic-max-cond-prob}-\ref{eq:deterministic-approx})
or~\eqref{eq:kernel-approx-stoch-model-0}, we can obtain the kernel-based estimator $\hat{v}$ of $Lu=u(\vx)$, that is,
\[
u(\vx)
=Lu\approx\hat{v}:=
L\vk_{\vL}^T\vK_{\vL}^{-1}\vf=\vk_{\vL}(\vx)^T\vK_{\vL}^{-1}\vf.
\]
Therefore, the kernel-based approximate solution $\hat{u}_{XZ}$ can be written as
\begin{equation}\label{eq:kernel-based-sol-elliptic-SPDE-1}
\hat{u}_{XZ}(\vx)
=\sum_{k=1}^Nc_kP_{\vy}K\left(\vx,\vx_k\right)+\sum_{k=1}^Mc_{N+k}B_{\vy}K\left(\vx,\vz_k\right),\quad\text{for }\vx\in\Domain,
\end{equation}
where the random coefficients $\vc:=\left(c_1,\cdots,c_{N+M}\right)^T$ are solved by the random linear system
\begin{equation}\label{eq:kernel-based-sol-elliptic-SPDE-2}
\vK_{\vL}\vc=\vf.
\end{equation}
Obviously $\hat{u}_{XZ}\in\Cont^m(\Domain)\subseteq\Hilbert^m(\Domain)$. We can further find that
the kernel basis $\vk_{\vL}$ of $\hat{u}_{XZ}$ are deterministic and the random coefficients $\vc$ dominate the stochastic structures of $\hat{u}_{XZ}$. Thus, we can design
the following algorithm to obtain the thousands $p$ sample paths of $\hat{u}_{XZ}$ to simulate the probability distributions of $u$, that is,
\[
\tag{A2}
\begin{split}
&\text{Initialize:}\\
&\ve^T:=\vk_{\vL}^T\vK_{\vL}^{-1},\\
&\text{Repeat }i=1,\ldots,p\\
&\quad{}\text{Simulate }\vvartheta_X^{(i)}:=\left(\vartheta_{\vx_1}^{(i)},\cdots,\vartheta_{\vx_N}^{(i)}\right)^T\sim\Normal\left(\vzero,\vTheta_X\right),\\
&\quad{}\vf^{(i)}:=
\left(\Gamma\big(\vx_1,\vartheta_{\vx_1}^{(i)}\big),\cdots,\Gamma\big(\vx_N,\vartheta_{\vx_N}^{(i)}\big),g\big(\vz_1\big),\cdots,g\big(\vz_M\big)\right)^T,\\
&\quad{}\hat{u}_{XZ}^{(i)}:=\ve^T\vf^{(i)}.
\end{split}
\]

{\bf Error Analysis:}
Finally, we study with the error analysis of the kernel-based approximate solution $\hat{u}_{XZ}$.
Same as Equation~\eqref{eq:stoch-model-error-erfc} in Lemma~\ref{l:stoch-model-error}, for any $\epsilon>0$, we have
\begin{equation}\label{eq:elliptic-SPDE-error-1}
\hat{\PP}\left(\abs{u(\vx)-\hat{u}_{XZ}(\vx)}\geq\epsilon\right)=
\hat{\PP}\left(\abs{Lu-L\hat{u}_{XZ}}\geq\epsilon\right)
=\erfc\left(\frac{\epsilon}{\sqrt{2}\sigma_{L|\vL}}\right),
\end{equation}
where $\hat{\PP}=\PP_K\times\PP$ is the product probability measure given in Equation~\eqref{eq:product-prob}
and $\sigma_{L|\vL}$ is the standard deviation defined as in Equation~\eqref{eq:standard-dev}.
This indicates that the convergence of $\hat{u}_{XZ}$ is dependent of $\sigma_{L|\vL}$.
Now we verify that the standard deviation $\sigma_{L|\vL}$ is equal to the \emph{generalized power function} $\Power_{\vL}(L)$. \cite[Section~16.1]{Wendland2005} shows that the generalized power function $\Power_{\vL}(L)$ is defined by
\[
\Power_{\vL}(L):=\min_{\Lambda\in\Span\left\{\vL\right\}}\norm{L-\Lambda}_{\Hilbert_K(\Domain)'}
=\sqrt{L_{\vx}L_{\vy}K(\vx,\vy)-L\vk_{\vL}^T\vK_{\vL}^{-1}L\vk_{\vL}},
\]
where $\Hilbert_K(\Domain)'$ is the dual space of the reproducing kernel Hilbert space $\Hilbert_K(\Domain)$.
Thus, we have
\begin{equation}\label{eq:elliptic-SPDE-error-2}
\sigma_{L|\vL}=\Power_{\vL}(L).
\end{equation}
According to the theorems in \cite[Section~16.3]{Wendland2005}, we can obtain the upper bounds of $\Power_{\vL}(L)$, that is,
\begin{equation}\label{eq:elliptic-SPDE-error-3}
\Power_{\vL}(L)=\Power_{\vL}(\delta_{\vx})=\Order\big(h_X^{m-2}\big)+\Order\big(h_Z^{m}\big),\quad
\text{when }h_X,h_Z\text{ are small enough},
\end{equation}
where the fill distances
\[
h_X:=\sup_{\vx\in\Domain}\min_{k=1,\ldots,N}\norm{\vx-\vx_k}_2,\quad
h_Z:=\sup_{\vz\in\partial\Domain}\min_{k=1,\ldots,M}d_2(\vz,\vz_k),
\]
and $d_2:\partial\Domain\times\partial\Domain\to[0,\infty)$ is the standard distance function on the manifolds.
For example, if the boundary $\partial\Domain$ is the unit sphere, then $d_2(\vx,\vy):=\cos^{-1}\big(\vx^T\vy\big)$.
For convenience, we transfer Equation~\eqref{eq:elliptic-SPDE-error-3} to
\begin{equation}\label{eq:elliptic-SPDE-error-4}
\Power_{\vL}(L)=\Order\left(h_{XZ}^{m-2}\right),
\quad\text{when }h_{XZ}\text{ is small enough},
\end{equation}
where
\[
h_{XZ}:=\max\left\{h_X,h_Z\right\}.
\]

%////////////////////////////////////////////////////////////////////////////////////////////////////////////////////////
\begin{remark}
The rough proofs of the convergent rates of the generalized power functions are checked by the error bounds
\[
\Power_{\vL}\left(\delta_{\vx}\circ P\right)\leq
\Power_{\vP}\left(\delta_{\vx}\circ P\right)
=\Order\big(h_X^{m-2}\big),\quad
\Power_{\vL}\left(\delta_{\vz}\circ B\right)\leq
\Power_{\vB}\left(\delta_{\vz}\circ B\right)
=\Order\big(h_Z^{m}\big),
\]
where $\vP:=\left(\delta_{\vx_1}\circ P,\cdots,\delta_{\vx_N}\circ P\right)^T$ and $\vB:=\left(\delta_{\vz_1}\circ B,\cdots,\delta_{\vz_M}\circ B\right)^T$ (see \cite[Theroem~16.10 and 16.11]{Wendland2005}).
Then the good designs of the data points $X$ and $Z$ are obviously $h_{X}^{m-2}\approx h_Z^{m}$.
In this article, we ignore the proofs of the error bounds of the (generalized) power functions.
The deep discussions of the convergent rates of the power functions can be found in many well-known publications of meshfree approximation such as the books \cite{Wendland2005,Fasshauer2007}.
\end{remark}
%////////////////////////////////////////////////////////////////////////////////////////////////////////////////////////

Combining Equations~\eqref{eq:elliptic-SPDE-error-1},~\eqref{eq:elliptic-SPDE-error-2}, and~\eqref{eq:elliptic-SPDE-error-4},
we have
\[
\hat{\PP}\left(\abs{u(\vx)-\hat{u}_{XZ}(\vx)}\geq\epsilon\right)=\Order\left(\frac{h_{XZ}^{m-2}}{\epsilon}\right),
\quad\text{when }h_{XZ}\text{ is small enough};
\]
hence
\[
\lim_{h_{XZ}\to0}\hat{\PP}\left(\abs{u(\vx)-\hat{u}_{XZ}(\vx)}\geq\epsilon\right)=0.
\]
Moreover, by the compactness of the domain $\Domain$,
we can even conclude that
\[
\lim_{h_{XZ}\to0}\hat{\PP}\left(\norm{u-\hat{u}_{XZ}}_{\Leb_{\infty}(\Domain)}\geq\epsilon\right)=0.
\]
Therefore, we have:
%////////////////////////////////////////////////////////////////////////////////////////////////////////////////////////
\begin{proposition}\label{p:elliptic-SPDE-error}
Suppose $\hat{u}_{XZ}$ is the kernel-based approximate solution of the elliptic SPDE~\eqref{eq:ellptic-SPDE} in Equations~(\ref{eq:kernel-based-sol-elliptic-SPDE-1}-\ref{eq:kernel-based-sol-elliptic-SPDE-2}).
Then $\hat{u}_{XZ}$ converges to $u$ uniformly in probability when $h_{XZ}\to0$.
\end{proposition}
%////////////////////////////////////////////////////////////////////////////////////////////////////////////////////////

\emph{Comments:}
In this section, we generalize
the kernel-based methods for the deterministic PDEs to the stochastic PDEs.
We show that the formulas and the error bounds of the kernel-based approximate solutions for the elliptic SPDEs
are consistent with the classical results of meshfree approximation.

In the following, we will compare the kernel-based methods with another current popular numerical methods for the elliptic SPDEs
such as the Galerkin finite element methods~\cite{BabuvskaTemponeZouraris2004,KuoSchwabSloan2012}
and the stochastic collocation methods~\cite{BabuvskaNobileTempone2010}.
\begin{itemize}
\item Both the Galerkin finite element methods and the stochastic collocation methods use the polynomial basis to obtain the numerical solutions of the SPDEs. But, the kernel-based approximate solutions can be constructed by the non-polynomial basis.

\item By the Galerkin finite element methods or the stochastic collocation methods, we need to choose the typical grid points to construct the meshes. But, the kernel-based methods are the meshfree methods and the data points can be placed at rather arbitrarily scattered locations. This indicates that the random designs of data points are still feasible for the kernel-based methods such as Sobol points.

\item The kernel-based methods are robust for any high-dimensional SPDE with the complex boundaries.

%\item In the computational processes of kernel-based approximation, the data points are only needed in the deterministic domains. But, both the Galerkin finite element methods and the stochastic collocation methods require the scattered points separately in the deterministic domains and the sample spaces of the probability spaces.

\item Usually, both the Galerkin finite element methods and the stochastic collocation methods need to know the Karhunen-Lo\`{e}ve expansion of the given random term $\vartheta$ such that we can truncate the original probability spaces to the finite dimensional probability spaces for the computations. But, we can simulate $\vartheta$ directly to construct the kernel-based approximate solutions.

\item The covariance (interpolating) matrixes for the kernel-based methods are not affected by the random term $\vartheta$. This indicates that we can construct the efficient kernel-based algorithms to obtain the thousands of sample paths to simulate the probability distributions.
\end{itemize}

%---------------------------------------------------------------------------------------------------------------------
\section{Parabolic Stochastic Partial Differential Equations}\label{sec:par-SPDE}
%---------------------------------------------------------------------------------------------------------------------

We know that the numerical analysis of the kernel-based methods for the parabolic PDEs is a delicate and non-trivial question.
In this section, we will extend the
the kernel-based methods for the deterministic parabolic PDE in~\cite{HonSchabackZhong2013} to the stochastic parabolic SPDE driven by the time and space white noises.
The recent paper~\cite{HonSchabackZhong2013} mainly focuses on the 1D parabolic equations.
For convenience of the comparison with \cite{HonSchabackZhong2013}, we only investigate
the 1D white noises.

Let $W$ be a time and space white noise with the mean $0$ and the spatial covariance kernel $\Phi:[0,1]\times[0,1]\to\RR$ defined on a probability space $\left(\Omega,\Filter,\left\{\Filter_t\right\}_{t\geq0},\PP\right)$, that is, $\Mean\left(W_{t}(x)\right)=0$ and $\Cov\left(W_{t}(x),W_{s}(y)\right)=\min\left\{t,s\right\}\Phi(x,y)$ for $x,y\in[0,1]$ and $t,s\geq0$.
The white noise $W$ does not exist the derivatives at the time $t$; but $W$ can be smooth at the space $x$.
The spatial covariance kernel $\Phi$ is only related to the space $x$.
For example in \cite[Section~3.2]{Chow2007}, the time and space white noise $W$ is constructed by a sequence of the i.i.d. standard scalar Brownian motions $\left\{W_n\right\}_{n=1}^{\infty}$, that is,
\[
W_t(x):=\sum_{n=1}^{\infty}\frac{W_{n,t}}{2n^2\pi^2}\sin\left(n\pi x\right);
\]
hence the spatial covariance kernel $\Phi$ has the form
\[
\Phi(x,y)=\frac{1}{4n^4\pi^4}\sin\left(n\pi x\right)\sin\left(n\pi y\right).
\]

Now we look at a parabolic SPDE driven by the white noise $W$,
\begin{equation}\label{eq:parabolic-SPDE}
\begin{cases}
~~~\ud U_t=\Delta U_t\ud t+\ud W_t,&\text{in }[0,1],~0\leq t\leq T,\\
U_{t}(0)=U_{t}(1)=0,&\text{on }\{0,1\},~0\leq t\leq T,\\
~~~~U_0=u^0,
\end{cases}
\end{equation}
where $\Delta:=\ud^2/\ud x^2$ is a Laplace differential operator and $u^0\in\Hilbert^m([0,1])$ for $m>2+1/2$.
Suppose that the solution $U_t\in\Hilbert^{m}\left([0,1]\right)$ for all $0\leq t\leq T$.

{\bf Discrete Kernel-based Approximate Solutions:}
Combining with the explicit Euler schemes, we will use
a positive definite kernel $K\in\Cont^{2m,1}\left([0,1]\times[0,1]\right)$ to construct the discrete kernel-based approximate solutions of the SPDE~\eqref{eq:parabolic-SPDE} in the following steps.
\begin{itemize}
\item[(S1)] Let $t_i:=iT/n$, $\delta t:=t_{i}-t_{i-1}=T/n$, and $\delta W_i:=W_{t_{i}}-W_{t_{i-1}}$ for $i=1,\ldots,n$. Then $\delta W_i$ is a Gaussian field with the mean $0$ and the covariance kernel $\delta t\Phi$ defined on the probability space $\left(\Omega,\Filter_{t_i},\PP\right)$.
By the explicit Euler schemes, we discretize the SPDE~\eqref{eq:parabolic-SPDE} at the discrete time $t_i$, that is,
\begin{equation}\label{eq:explict-Euler}
U_{t_i}-U_{t_{i-1}}=\Delta U_{t_{i-1}}\delta t+\delta W_i,\quad \text{for }i=1,\ldots,n.
\end{equation}
We continue to approximate the values of $U_{t_i}$ at the space data points $X:=\left\{x_k\right\}_{k=1}^{N}$ and $Z:=\left\{x_0,x_{N+1}\right\}$ such as
    $0=x_0<x_1<\cdots x_{N}<x_{N+1}=1$.
\item[(S2)] Let $u^{i-1}:=U_{t_{i-1}}$, $u^{i}:=U_{t_{i}}$, and $\vartheta^i:=\delta W_i$. If we already have the information
     \[
     r^{i-1}:=\left(\delta t\Delta+I\right)u^{i-1},
     \]
     at the previous time step $t_{i-1}$, then
     the Euler scheme~\eqref{eq:explict-Euler} provides that we can obtain the solution $u^{i}$ by the simulations of the Gaussian field $\vartheta^i$, that is,
    \[
    u^{i}(x)=r^{i-1}(x)+\vartheta^i_x,\quad\text{for }x\in(0,1),\quad u^{i}(0)=u^{i}(1)=0,
    \]
    because the white noise increment $\delta W_i$ is independent of $U_{t_{i-1}}$ at the current time step $t_i$.
    Next, we need to approximate
     \[
     r^{i}:=\left(\delta t\Delta+I\right)u^{i},
     \]
     for the computations at the next time step $t_{i+1}$. Let
     \[
     L:=\delta_{\vx_j}\circ\left(\delta t\Delta+I\right),\quad\text{for }1\leq j\leq N.
     \]
     Then $L$ is a bounded linear functional on $\Hilbert^{m}([0,1])$ by the Sobolev Imbedding Theorem.
     Now we construct the kernel-based estimator $L\hat{u}_{XZ}^i$ of $Lu^i$ by the chosen positive definite kernel $K$, that is,
     \[
     r^i(x_j)=Lu^i\approx L\hat{u}_{XZ}^i.
     \]
     Firstly, we simulate the multivariate normal random vector
    \[
    \vvartheta^i_{X}:=\left(\vartheta^i_{x_1},\cdots,\vartheta^i_{x_{N}}\right)^T\sim\Normal\left(\vzero,\delta t\vTheta_{X}\right).
    \]
    Since
    \[
    u^{i}(x_0)=0,
    u^{i}(x_1)=r^{i-1}(x_1)+\vartheta^i_{x_1},\ldots,u^{i}(x_{N})=r^{i-1}(x_{N})+\vartheta^i_{x_{N}},
    u^{i}(x_{N+1})=0,
    \]
    we can
    obtain the stochastic data evaluated by $u^i$ at $X$ such as
    \[
    f_0:=u^{i}(x_0),f_1:=u^{i}(x_1),\ldots,f_N:=u^{i}(x_N),f_{N+1}:=u^{i}(x_{N+1})
    \]
    and
    \[
    \vf:=\left(f_0,f_1,\cdots,f_N,f_{N+1}\right)^T.
    \]
    Thus, by Equations~(\ref{eq:kernel-approx-stoch-model-1}-\ref{eq:kernel-approx-stoch-model-2}), we have the kernel-based approximate function   \[
    \hat{u}^{i}_{XZ}(x):=\vk_{XZ}^T(x)\vK_{XZ}^{-1}\vf,\quad\text{for }x\in[0,1],
    \]
    where $\vk_{XZ}(x):=\left(K(x,x_0),\cdots,K(x,x_{N+1})\right)^T$ and $\vK_{XZ}:=\left(K(x_k,x_l)\right)_{k,l=0}^{N+1,N+1}$.
    This indicates that
    \[
    L\hat{u}_{XZ}^i:=\hat{u}_{XZ}^{i}(x_j)+\delta t\Delta \hat{u}_{XZ}^{i}(x_j)
    =f_j+\delta t\Delta\vk_{XZ}^T(x_j)\vK_{XZ}^{-1}\vf.
    \]
\item[(S3)] Repeat the step (S2) for all $i=1,\ldots,n$. Here $u^{0}$ is given and the estimation of $r^{n}=\left(\delta t\Delta+I\right)u^{n}$ is not necessary.
\end{itemize}

Moreover, the algorithms of the discrete kernel-based approximate solutions given in the above step (S1-S3) can be written as follows:
\[
\tag{A3}
\begin{split}
&\text{Initialize:}\\
&\vf:=\left(0,u^{0}(x_1),\ldots,u^{0}(x_{N}),0\right)^T,\\
&\text{Repeat i=1,\ldots,n}\\
&\quad{}\text{Simulate }\vvartheta^i_{X}:=\left(\vartheta^i_{x_1},\cdots,\vartheta^i_{x_{N}}\right)^T\sim\Normal\left(\vzero,\delta t\vTheta_{X}\right),\\
&\quad{}\hat{u}^{i}_j:=f_j+\delta t\Delta\vk_{XZ}(x_j)^T\vK_{XZ}^{-1}\vf+\vartheta^i_{x_j},\text{for }j=1,\ldots,N,\\
&\quad{}\text{update }\vf:=\left(0,\hat{u}^{i}_1,\ldots,\hat{u}^{i}_N,0\right)^T\text{ when }i<n,\\
\end{split}
\]

%////////////////////////////////////////////////////////////////////////////////////////////////////////////////////////
\begin{remark}
Algorithm (A3) is different from Algorithms (A1) and (A2). Here, Algorithm (A3) only produces one sample path and we need to repeat Algorithm (A3) to obtain the thousands $p$ sample paths to approximate the probability distributions of the solutions $U_t$.
\end{remark}
%////////////////////////////////////////////////////////////////////////////////////////////////////////////////////////

{\bf Error Analysis:}
Finally, we study with the convergence of the kernel-based approximate solutions of the SPDE~\eqref{eq:parabolic-SPDE}.
Same as in Equation~\eqref{eq:product-prob}, we define
the tensor product probability space
\[
\hat{\Omega}:=\Omega_m\times\Omega,\quad \hat{\Filter}:=\Filter_m\otimes\Filter,\quad
\hat{\Filter}_t:=\Filter_m\otimes\Filter_t\text{ for }t\geq0,
\quad
\hat{\PP}:=\PP_K\times\PP,
\]
such that the convergence of the kernel-based estimators is well-posed on this probability space by Proposition~\ref{p:stoch-model-error}.

Now we look at the local errors of the kernel-based approximate solutions.
The It\^o-Taylor expansion of $U_t$ guarantees that
\[
U_{t_1}=U_{t_0}+\Delta U_{t_0}\int_{t_0}^{t_1}\ud s+\int_{t_0}^{t_1}\ud W_s+R,
\]
and the remainder
\[
\Mean\big(R^2\big)=\Order\big(\delta t^{3}\big).
\]
Thus, we can obtain
the local truncation errors at time in probability, that is,
\[
U_{t_{i}}-U_{t_{i-1}}-\Delta U_{t_{i-1}}\delta t-\delta W_i\overset{\hat{\PP}}{=}\Order\big(\delta t^{3/2}\big),\quad
\text{when }\delta t\text{ is small enough},
\]
for $i=1,\ldots,n$. This indicates that
\begin{equation}\label{eq:heat-SPDE-local-error-1}
u^{i}(x_j)-u^{i-1}(x_j)-\Delta u^{i-1}(x_j)\delta t-\vartheta^i_{x_j}\overset{\hat{\PP}}{=}\Order\big(\delta t^{3/2}\big),\quad
\text{when }\delta t\text{ is small enough},
\end{equation}
for $i=1,\ldots,n$ and $j=1,\ldots,N$.
Here, the notation $u-\hat{u}\overset{\hat{\PP}}{=}\Order\left(\delta\right)$ means that $\hat{u}$ converges to $u$ in probability when $\delta\to0$.

%////////////////////////////////////////////////////////////////////////////////////////////////////////////////////////
\begin{remark}
Roughly, the It\^o-Taylor expansion is based on the iterated application of the It\^o formula.
Since the white noises do not have the continuous derivatives at time,
the convergent orders of the Euler schemes of the SPDEs are lower than the PDEs.
More details of the Euler schemes of the stochastic differential equations can be found in \cite[Section 10.2]{KloedenPlaten1992} and~\cite[Section~6.3]{JentzenKloeden2011}.
\end{remark}
%////////////////////////////////////////////////////////////////////////////////////////////////////////////////////////

Moreover, Equation~\eqref{eq:stoch-model-error-fill-distance} provides another local errors at space in probability
\begin{equation}\label{eq:heat-SPDE-local-error-2}
\Delta u^{i-1}(x_j)-\Delta\hat{u}^{i-1}_{XZ}(x_j)\overset{\hat{\PP}}{=}\Order\left(h_{XZ}^{m-2}\right),
\quad\text{when }h_{XZ}\text{ is small enough},
\end{equation}
for $i=1,\ldots,n$ and $j=1,\ldots,N$. Here $h_{XZ}=\max_{k=1,\ldots,N+1}\abs{x_{k}-x_{k-1}}/2$.

Next, we estimate the global errors
\[
\ve^{i}:=\begin{pmatrix}U_{t_{i}}(x_0)\\U_{t_{i}}(x_1)\\\vdots\\U_{t_{i}}(x_{N})\\U_{t_{i}}(x_{N+1})\end{pmatrix}-
\begin{pmatrix}0\\\hat{u}^{i}_1\\\vdots\\\hat{u}^{i}_N\\0\end{pmatrix},
\quad\text{for }i=1,\ldots,n.
\]
Combining the local errors in Equations~\eqref{eq:heat-SPDE-local-error-1} and \eqref{eq:heat-SPDE-local-error-2}, we have
\begin{equation}\label{eq:heat-SPDE-global-error-1}
\ve^{i}\overset{\hat{\PP}}{=}\ve^{i-1}+\delta t\vK_{XZ}''\vK_{XZ}^{-1}\ve^{i-1}+\Order\left(\delta th_{XZ}^{m-2}\right)+\Order\left(\delta t^{3/2}\right),
\end{equation}
when $\delta t,h_{XZ}$ are small enough. Here $\vK_{XZ}'':=\left(\Delta_{x}K(x_k,x_l)\right)_{k,l=0}^{N+1,N+1}$.
According to \cite[Theorem~7.2]{HonSchabackZhong2013}, which is verified by the sampling inequality in~\cite{RiegerZwicknaglSchaback2010},
the spectral radius of $\vK_{XZ}''\vK_{XZ}^{-1}$ satisfies
\[
\rho\left(\vK_{XZ}''\vK_{XZ}^{-1}\right)=\Order\left(h_{XZ}^{-2}\right),
\quad\text{when }h_{XZ}\text{ is small enough}.
\]
Since the explicit Euler schemes are used here, we naturally need the \emph{Courant-Friedrichs-Lewy condition}, that is,
\[
\frac{\delta t}{h_{XZ}^2}=\Order\left(1\right),\quad\text{when }\delta t,h_{XZ}\text{ are small enough}.
\]
Then, by the induction of Equation~\eqref{eq:heat-SPDE-global-error-1}, we notes that
\[
\frac{1}{\sqrt{N}}\norm{\ve^{n}}_2\overset{\hat{\PP}}{=}
\Order\left(\frac{\delta th_{XZ}^{m-2}}{\delta t}\right)+\Order\left(\frac{\delta t^{3/2}}{\delta t}\right)
=\Order\left(h_{XZ}^{m-2}\right)+\Order\left(\delta t^{1/2}\right),
\]
when $\delta t,h_{XZ}$ are small enough; hence we can conclude that
\[
\lim_{\delta t,h_{XZ}\to0}\hat{u}^{i}_j\overset{\hat{\PP}}{=}U_{t_i}(x_j),
\]
for $i=1,\ldots,n$ and $j=1,\ldots,N$.

%////////////////////////////////////////////////////////////////////////////////////////////////////////////////////////
\begin{proposition}\label{p:parabolic-SPDE-error}
Suppose that $\hat{u}^{i}_j$ is the discrete kernel-based approximate solution of the parabolic SPDE~\eqref{eq:parabolic-SPDE} in Algorithm (A3). If the Courant-Friedrichs-Lewy condition is well-posed, then $\hat{u}^{i}_j$ converges to $U_{t_i}(x_j)$ in probability when $\delta t,h_{XZ}\to0$, for all $i=1,\ldots,n$ and $j=1,\ldots,N$.
\end{proposition}
%////////////////////////////////////////////////////////////////////////////////////////////////////////////////////////

\emph{Comments:}
In this section, we only discuss the 1D parabolic SPDEs.
Moreover, we can update Algorithm (A3) to the high-dimensional domains.
Same as the numerical experiments in \cite{HonSchabackZhong2013},
the boundary $\left\{0,1\right\}$ can be extended to the discrete data points $Z\subseteq\partial\Domain$.
However, the paper \cite{HonSchabackZhong2013} has not given the proofs of the convergence of the high-dimensional parabolic PDEs and the technique points of the proofs could be the spectral radius of the associated matrix $\vK_{XZ}''\vK_{XZ}^{-1}$.
So, we do not investigate the high-dimensional parabolic SPDEs currently.

%---------------------------------------------------------------------------------------------------------------------
\section{Numerical Examples}\label{sec:num-exp}
%---------------------------------------------------------------------------------------------------------------------

In this section, we will give the 3D, 2D, and 1D numerical examples of the kernel-based estimators and the kernel-based approximate solutions in Sections~\ref{sec:Stochastic-Data-Approx}-\ref{sec:par-SPDE}.
The kernel-based algorithms will be constructed by the Gaussian kernels, the compactly supported kernels (Wendland functions), and the Sobolev-spline kernels (Mat\'ern functions).

%---------------------------------------------------------------------------------------------------------------------
\subsection{Stochastic Data Interpolations}\label{sec:num-3D-model}
%---------------------------------------------------------------------------------------------------------------------

Let the data points $X$ be the Halton points in the unit cube $\Domain:=[0,1]^3$.
Suppose that the stochastic data $\vf$ at the data points $X$ are obtained by the simple 3D stochastic model
\begin{equation}\label{eq:num-stoch-model}
u(\vx)=\vartheta_{\vx}^2,\quad\text{for }\vx:=\left(x_1,x_2,x_3\right)\in\Domain,
\end{equation}
where $\vartheta_{\vx}:=\zeta\phi(\vx)$ is composed of $\phi(\vx):=\sin\left(\pi x_1\right)\sin\left(2\pi x_2\right)\sin\left(3\pi x_3\right)$ and $\zeta\sim\Normal(0,1)$.
Then $\vartheta$ is a Gaussian field with the mean $0$ and the covariance kernel $\Phi(\vx,\vy):=\phi(\vx)\phi(\vy)$.
Notes that $u\in\Cont^{\infty}(\Domain)\subseteq\Hilbert^4(\Domain)$ and $\delta_{\vx}\circ\Delta$ is a bounded linear functional of $\Hilbert^4(\Domain)$; hence the target random variable $\Delta u(\vx)$ is well-defined for any $\vx\in\Domain$.

We will use a Gaussian kernel with a shape parameter $\theta>0$
\[
K_{\theta}(\vx,\vy):=e^{-\theta^2\norm{\vx-\vy}_2^2},
\quad\text{for }\vx,\vy\in\Domain,
\]
to construct the kernel-based estimator $\Delta\hat{u}_X(\vx)$ of $\Delta u(\vx)$ given in Equations~\eqref{eq:kernel-approx-stoch-model-0} or (\ref{eq:kernel-approx-stoch-model-1}-\ref{eq:kernel-approx-stoch-model-2}). Here, we can view $L:=\delta_{\vx}\circ\Delta$.

%////////////////////////////////////////////////////////////////////////////////////////////////////////////////////////
\begin{remark}
Usually, the shape parameters of the kernels are used to control the shapes of the kernel basis.
The shape parameters $\theta$ of the Gaussian kernels are chosen empirically and are based on the personal experiences. In this article, we do not investigate the optimal shape parameters.
\end{remark}
%////////////////////////////////////////////////////////////////////////////////////////////////////////////////////////

\begin{figure*}
\begin{minipage}{0.5\textwidth}
\center
\includegraphics[width=\textwidth, height=\textwidth]{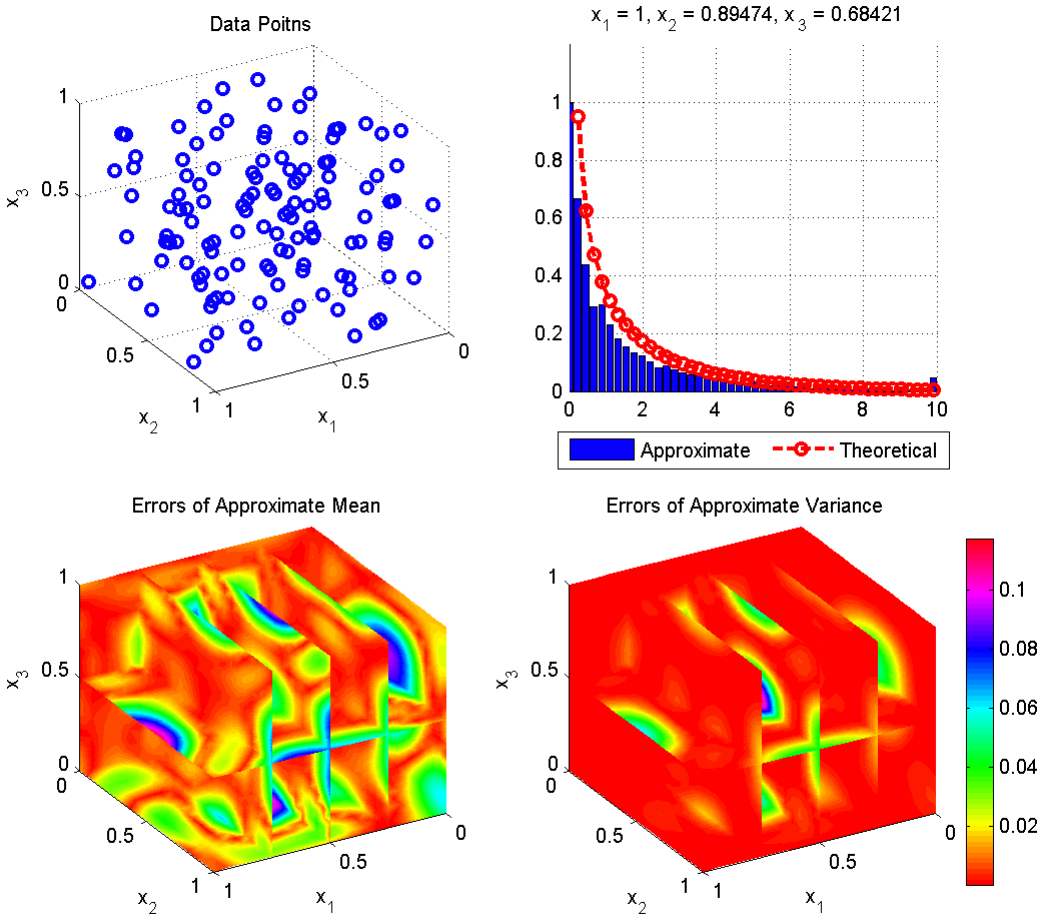}
{\small Approximate probability distributions.}
\end{minipage}
\begin{minipage}{0.5\textwidth}
\center
\includegraphics[width=\textwidth, height=\textwidth]{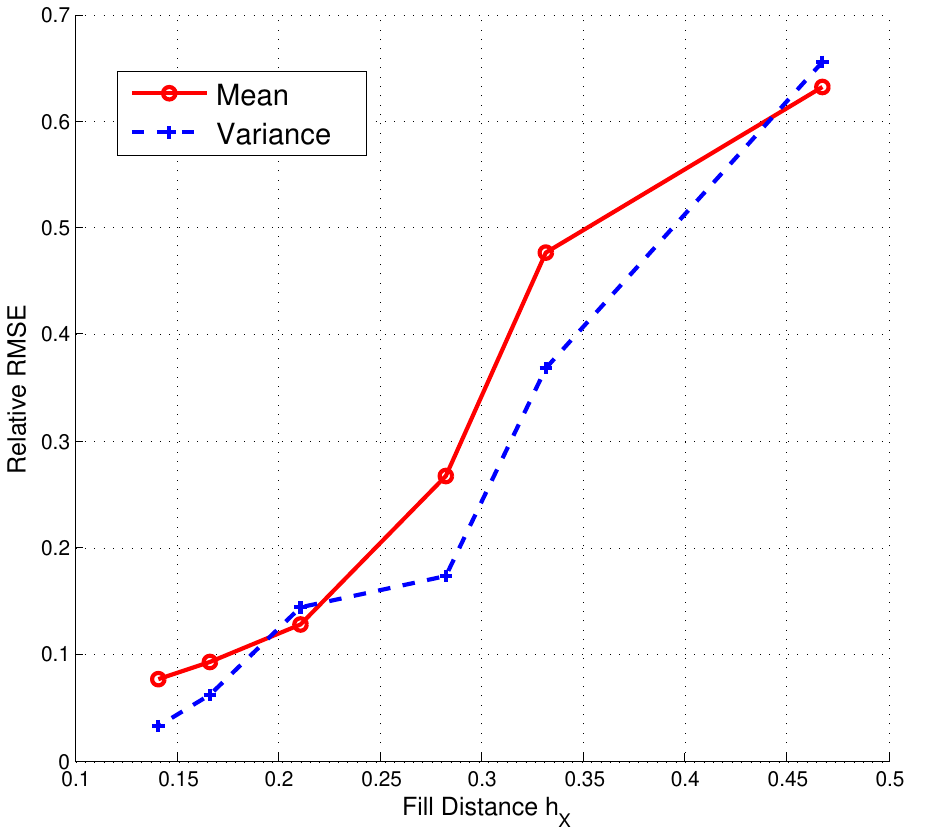}
{\small Convergence of approximate means and variances.}
\end{minipage}
\caption{The kernel-based estimator $\Delta\hat{u}_X(\vx)$ of the target random variable $\Delta u(\vx)$ in Equation~\eqref{eq:num-stoch-model} for $\vx\in[0,1]^3$:
The kernel-based estimators $\Delta\hat{u}_X(\vx)$ are constructed by the Gaussian kernel $K_{\theta}$ with the shape parameter $\theta=2.6$ (left) and $\theta=6$ (right) in Equations~\eqref{eq:kernel-approx-stoch-model-0} or (\ref{eq:kernel-approx-stoch-model-1}-\ref{eq:kernel-approx-stoch-model-2}).
The approximate probability distributions of $\Delta\hat{u}_X(\vx)$ are simulated by the $p=10000$ samples.
The left parts are the numerical experiments including the $N=125$ Halton Points $X$ shown in blue points (the top left panel), the approximate and theoretical probability density functions of $\Delta\hat{u}_{X}(\vx)$ and $\Delta u(\vx)$ for $\vx=(1,0.894,0.684)$ (the top right panel), and the approximate means and variances of $\Delta\hat{u}_{X}(\vx)$ (the bottom panels), where the right color bar represents the values of relative absolute errors. The right parts are the numerical experiments of the relative RMSE of the approximate means and variances of $\Delta\hat{u}_{X}(\vx)$ for different Halton points.
Here, the root mean square error is denoted by $\text{RMSE}:=\sqrt{\int_{\Domain}\abs{u(\vx)-\hat{u}(\vx)}^2\ud\vx}$
where $\hat{u}(\vx)$ is the estimator of $u(\vx)$.
}\label{fig:Kernel-3D-Stoch}
\end{figure*}

Notes that the unit cube $\Domain$ is not a Lipschitz domain.
However, by Figure~\ref{fig:Kernel-3D-Stoch}, the approximate probability distributions of $\Delta\hat{u}_X(\vx)$ are still convergent to the Chi-squared probability distributions of $\Delta u(\vx)$.
This indicates that the regularity of the domains may not be the necessary conditions of the kernel-based estimators.

%---------------------------------------------------------------------------------------------------------------------
\subsection{Stochastic Poisson Equations}\label{sec:num-poisson-SPDE}
%---------------------------------------------------------------------------------------------------------------------

Let the domain $\Domain\subseteq\RR^2$ be a circle centered at origin with the radius $1/2$, that is, $\Domain:=\left\{\vx\in\RR^2:\norm{\vx}_2\leq1/2\right\}$.
Denote that
\begin{align*}
&\psi_1(\vx):=\sin\left(\pi(x_1-1/2)\right)\sin\left(\pi(x_2-1/2)\right),\\
&\psi_2(\vx):=\sin\left(2\pi(x_1-1/2)\right)\sin\left(2\pi(x_2-1/2)\right),\\
&\varphi_1(\vx):=\cos\big(2\pi\norm{\vx}_2^2\big),\quad \varphi_2(\vx):=\sin\big(4\pi\norm{\vx}_2^2\big),\\
&\phi_1(\vx):=2\pi\norm{\vx}_2^2\cos\big(2\pi\norm{\vx}_2^2\big)+\sin\big(2\pi\norm{\vx}_2^2\big),\\
&\phi_2(\vx):=4\pi\norm{\vx}_2^2\sin\big(4\pi\norm{\vx}_2^2\big)-\cos\big(4\pi\norm{\vx}_2^2\big),
\end{align*}
for $\vx:=\left(x_1,x_2\right)\in\Domain$.

Now we look at the stochastic Poisson equation with the trivial Dirichlet's boundary conditions such as
\begin{equation}\label{eq:Poission-SPDE}
\begin{cases}
-\Delta u=f+\vartheta,&\text{in }\Domain,\\
~~~~~u=g,&\text{on }\partial\Domain,
\end{cases}
\end{equation}
where $f:=2\pi^2\psi_1+8\pi^2\psi_2$, $g:=(\psi_1+\psi_2)|_{\partial\Domain}$, and
$\vartheta$ is a Gaussian field with the mean $0$ and the covariance kernel
$\Phi(\vx,\vy):=64\pi^2\phi_1(\vx)\phi_1(\vy)+64\pi^2\phi_2(\vx)\phi_2(\vy)$.
Then the solution of the SPDE~\eqref{eq:Poission-SPDE} can be represented as
\[
u(\vx)=\psi_1(\vx)+\psi_2(\vx)
+\zeta_1\varphi_1(\vx)+\frac{\zeta_2}{2}\varphi_2(\vx),
\]
where $\zeta_1,\zeta_2\sim\text{i.i.d.}\Normal\left(0,1\right)$.

Let the data points $X\subseteq\Domain$ and $Z\subseteq\partial\Domain$ be the Halton points and the evenly spaced points, respectively.
By Equations~(\ref{eq:kernel-based-sol-elliptic-SPDE-1}-\ref{eq:kernel-based-sol-elliptic-SPDE-2}),
we construct the kernel-based approximate solutions $\hat{u}_{XZ}$ of the SPDE~\eqref{eq:Poission-SPDE}
by a compactly supported kernel with a shape parameter $\theta>0$
\[
K_{\theta}(\vx,\vy):=\big(3+18\theta\norm{\vx-\vy}_2+35\theta^2\norm{\vx-\vy}_2^2\big)\big(1-\theta\norm{\vx-\vy}_2\big)_{+}^{6},
\quad\text{for }\vx,\vy\in\Domain,
\]
where $\left(\cdot\right)_{+}$ is the cutoff function, that is, $\left(r\right)_{+}=r$ when $r\geq0$ otherwise $\left(r\right)_{+}=0$.

\begin{figure*}
\begin{minipage}{0.5\textwidth}
\center
\includegraphics[width=\textwidth, height=\textwidth]{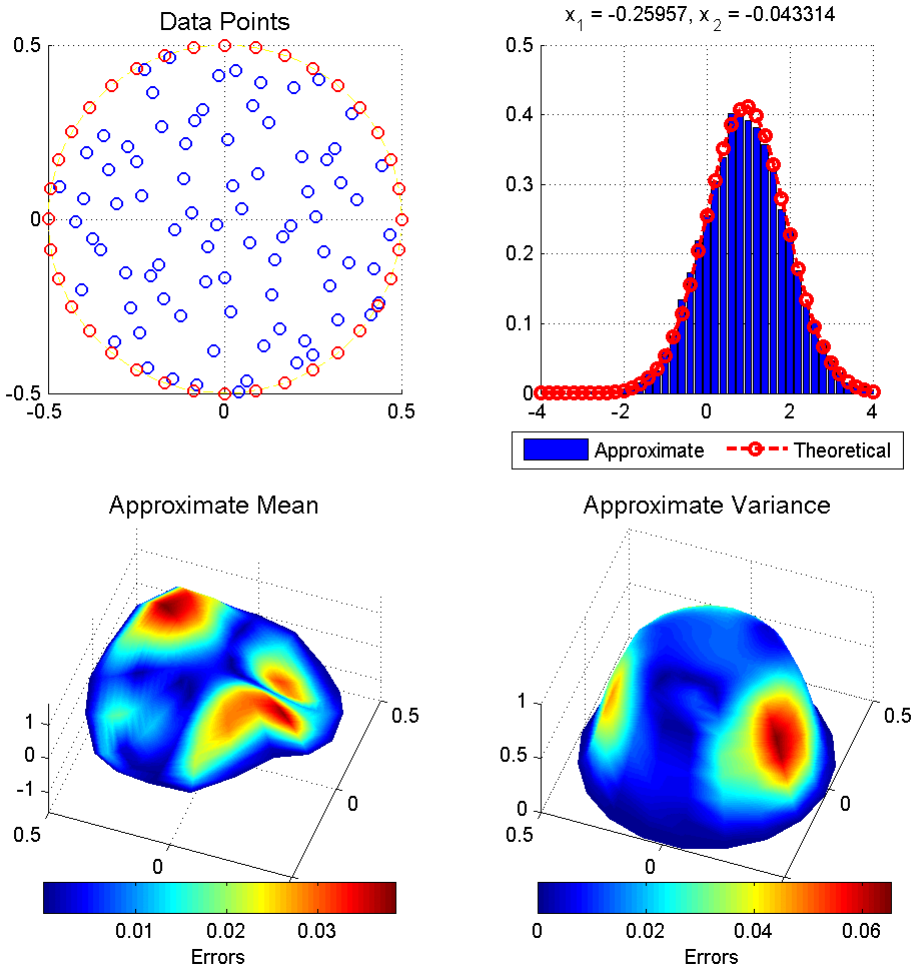}
{\small Approximate probability distributions.}
\end{minipage}
\begin{minipage}{0.5\textwidth}
\center
\includegraphics[width=\textwidth, height=\textwidth]{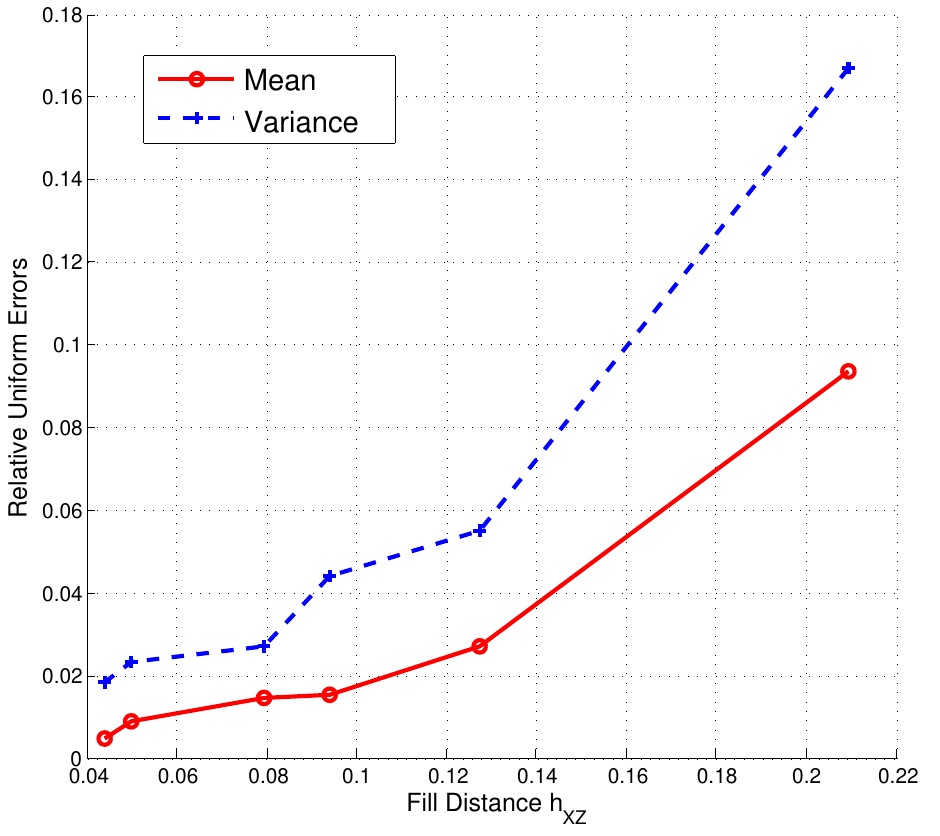}
{\small Convergence of approximate means and variances.}
\end{minipage}
\caption{The kernel-based approximate solutions $\hat{u}_{XZ}$ of the SPDE~\eqref{eq:Poission-SPDE}:
The kernel-based approximate solutions $\hat{u}_{XZ}$ are constructed by the compactly supported kernels $K_{\theta}$ with the different shape parameter $\theta=0.1$ (left) and $\theta=0.9$ (right) in Equations~(\ref{eq:kernel-based-sol-elliptic-SPDE-1}-\ref{eq:kernel-based-sol-elliptic-SPDE-2}).
The approximate probability distributions of $\hat{u}_{XZ}$ are simulated by the $p=10000$ simulated sample paths.
The left parts are the numerical experiments including the $N=80$ Halton Points $X$ shown in blue points and the $M=36$ evenly spaced points $Z$ shown in red points (the top left panel), the approximate and theoretical probability density functions of $\hat{u}_{XZ}$ and $u$ at $\vx=(-0.259,-0.043)$ (the top right panel), and the approximate means and variances of $\hat{u}_{XZ}$ (the bottom panels), where the bottom color bars represent the values of absolute errors. The right parts are the numerical experiments of the relative uniform errors of the approximate means and variances of $\hat{u}_{XZ}$ for different Halton points and evenly spaced points.}\label{fig:Kernel-Poisson-SPDE}
\end{figure*}

Comparing with the theoretical probability distributions of $u(\vx)$ in Figure~\ref{fig:Kernel-Poisson-SPDE},
the approximate probability distributions of $\hat{u}_{XZ}(\vx)$ are well-posed for $\vx\in\Domain$. Moreover, the approximate means and variances of $\hat{u}_{XZ}$ are convergent to the theoretical means and variances of $u$ uniformly on $\Domain$ when $h_{XZ}\to0$.
In Section~\ref{sec:ell-SPDE}, we require $m>2+d/2=3$ according to the conditions of Theorem~\ref{t:Gauss-PDK-L}.
Here, we find that $K_{\theta}$ belongs to $\Cont^{4,1}\left(\Domain\times\Domain\right)$ but not $\Cont^{6,1}\left(\Domain\times\Domain\right)$.
However, the kernel-based approximate solution $\hat{u}_{XZ}$ still works well for the approximations.
This indicates that the smooth conditions of the positive definite kernels could be weakened.

%---------------------------------------------------------------------------------------------------------------------
\subsection{Stochastic Heat Equations}\label{sec:num-heat-SPDE}
%---------------------------------------------------------------------------------------------------------------------

In this section, the time and space white noise
\[
W_t:=\sum_{n=1}^{\infty}W_{n,t}q_n\phi_n,
\]
is composed of a sequence of the i.i.d. standard scalar Bownian motions $\left\{W_n\right\}_{n=1}^{\infty}$.
Let
\[
q_n:=\frac{1}{n^2\pi^2},\quad
\phi_n(x):=\sqrt{2}\sin\left(n\pi x\right),\quad\text{for }n\in\NN.
\]
Since
the spatial covariance kernel $\Phi$ of the white noise $W$ has the form
\[
\Phi(x,y):=\sum_{n=1}^{\infty}q_n^2\phi_n(x)\phi_n(y),
\]
we have
\[
\Phi(x,y):=
\begin{cases}
-\frac{1}{6}x^3+\frac{1}{6}x^3y+\frac{1}{6}xy^3-\frac{1}{2}xy^2+\frac{1}{3}xy,&0\leq x\leq y\leq 1,\\
-\frac{1}{6}y^3+\frac{1}{6}xy^3+\frac{1}{6}x^3y-\frac{1}{2}x^2y+\frac{1}{3}xy,&0\leq y\leq x\leq 1.
\end{cases}
\]

Now we study with the stochastic heat equation
\begin{equation}\label{eq:heat-SPDE}
\ud U_t=\Delta U_t\ud t+\ud W_t,\quad\text{in }[0,1],~0\leq t\leq0.1,
\end{equation}
driven by the time and space white noise $W$.
If the SPDE~\eqref{eq:parabolic-SPDE} is endowed with the trivial Dirichlet's boundary conditions and the the initial condition $u^0(x):=\sqrt{2}\left(\sin(\pi x)+\sin(2\pi x)+\sin(3\pi x)\right)$,
then the solution of the SPDE~\eqref{eq:parabolic-SPDE} can be written as
\[
U_t(x)=\sum_{n=1}^{\infty}\eta_{n}(t)\phi_n(x),
\]
where
\[
\eta_{n}(t):=\mu_{n}e^{-n^2\pi^2t}+\frac{1}{q_n}\int_0^t
e^{n^2\pi^2(s-t)}\ud W_{n,s},
\text{ for }\mu_{n}:=\int_{0}^1u_0(x)\phi_n(x)\ud x.
\]

Let $X$ be the uniformly distributed points in $[0,1]$.
Algorithm (A3) provides the discrete kernel-based approximate solutions $\hat{u}^{i}_k$
of the SPDE~\eqref{eq:heat-SPDE} by
the Sobolev spline kernel with the shape parameter $\theta>0$
\[
K_{\theta}(x,y):=\big(15+15\theta\abs{x-y}+6\theta^2\abs{x-y}^2+\theta^3\abs{x-y}^3\big)e^{-\theta\abs{x-y}},\quad
\text{for }x,y\in[0,1].
\]

\begin{figure*}
\begin{minipage}{0.5\textwidth}
\center
\includegraphics[width=\textwidth, height=\textwidth]{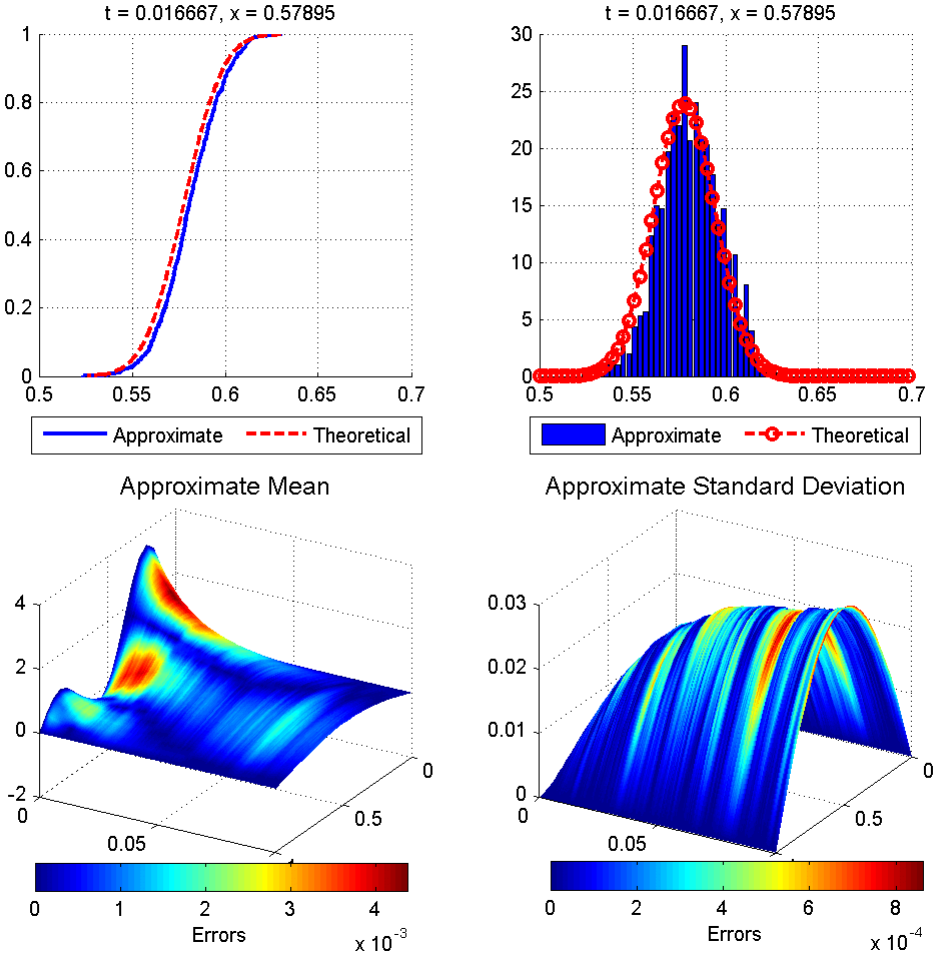}
{\small Approximate probability distributions.}
\end{minipage}
\begin{minipage}{0.5\textwidth}
\center
\includegraphics[width=\textwidth, height=\textwidth]{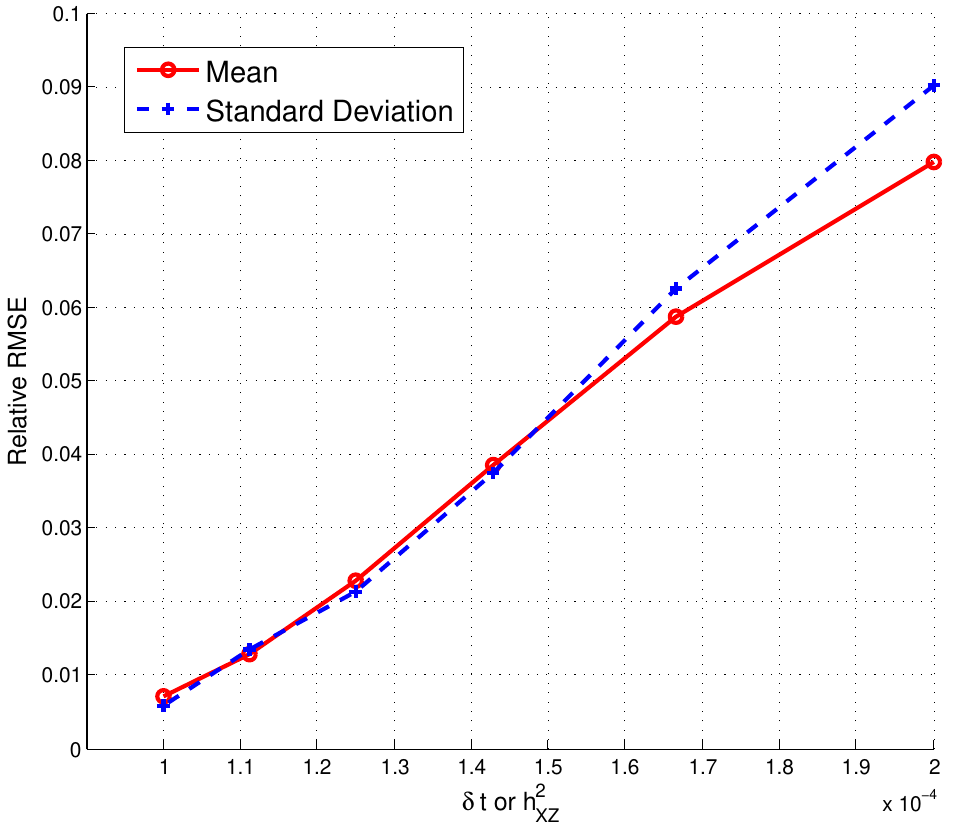}
{\small Convergence of approximate means and standard deviations.}
\end{minipage}
\caption{The discrete kernel-based approximate solutions $\hat{u}_j^{i}$ of the SPDE~\eqref{eq:parabolic-SPDE}:
The discrete kernel-based approximate solutions $\hat{u}_j^{i}$ are constructed by the Sobolev-spline kernels $K_{\theta}$ with the different shape parameter $\theta=2.6$ (left) and $\theta=30$ (right) in Algorithm (A3).
The approximate probability distributions of $\hat{u}_j^{i}$
are simulated by the $p=10000$ sample paths.
The left parts are the numerical experiments for $\delta t=0.000167$ and $h_{XZ}=0.0263$ including the approximate and theoretical cumulative distribution functions and probability density functions of $\hat{u}_j^{i}$ and $U_{t_i}(x_j)$ for $t_i=0.0167$ and $x_j=0.578$ (the top panels), and the approximate means and standard deviations of $\hat{u}_j^{i}$ (the bottom panels), where the bottom color bars represent the values of absolute errors. The right parts are the numerical experiments of the relative RMSE of the approximate means and standard deviations for different $\delta t$ and $h_{XZ}$ such that $\delta t\approx h_{XZ}^2$.}\label{fig:Kernel-Heat-SPDE}
\end{figure*}

By Figure~\ref{fig:Kernel-Heat-SPDE}, the probability distributions of $\hat{u}_j^{i}$ are the good approximations of the theoretical probability distributions of $U_{t_i}(x_j)$.
Moreover,
the approximate means and variances of $\hat{u}_j^{i}$ are convergent to the theoretical means and variances of $U_{t_i}(x_j)$ when the both $\delta t$ and $h_{XZ}$ tend to $0$. Here $\delta t$ and $h_{XZ}$ need to satisfy the Courant-Friedrichs-Lewy conditions. If not, the kernel-based approximate solutions will become unstable.

%---------------------------------------------------------------------------------------------------------------------
\section{Final Remarks}\label{sec:Final}
%---------------------------------------------------------------------------------------------------------------------

In this article, we try to combine approximation theory and statistical learning into one theoretical structure such that the best estimators have the both globally and locally geometrical meaning.
Here, we mainly focus on the connections of meshfree approximation and kriging interpolation by the Gaussian random variables defined on the Sobolev spaces. According to Theorem~\ref{t:Gauss-PDK-L} and Corollary~\ref{c:Gauss-PDK-multi-L}, the constructions of the multivariate normal random variables $L_1S,\ldots,L_NS$ give a connection of the interpolating data $L_1u=f_1,\ldots,L_Nu=f_N$ and the kernel basis $L_{1,\vy}K(\cdot,\vy),\ldots,L_{N,\vy}K(\cdot,\vy)$. Thus, we can use the statistics \& probability techniques to obtain the kernel-based estimators and the kernel-based approximate functions in Section~\ref{sec:ker-app}.
These kernel-based estimators are even consistent with the representations of the Hermite-Birkhoff interpolation in approximation theory.
Moreover, we obtain some new results in the both fields of meshfree approximation and kriging interpolation. But, these results have already been known well in one another field. Thus, we strongly believe that there could be some links between approximation theory and statistical learning such as the kernel-based methods discussed here.

%////////////////////////////////////////////////////////////////////////////////////////////////////////////////////////
\begin{remark}
In our original papers~\cite{CialencoFasshauerYe2012,FasshauerYe2012Bonn,FasshauerYe2012SPDEMCQMC,Ye2013NonSPDE}, we call the kernel-based methods the kernel-based collocation methods. But, some people may confuse the kernel-based collocation and the stochastic collocation in \cite{BabuvskaNobileTempone2010}.
In fact, the kernel-based collocation and the stochastic collocation are different, more precisely, the kernel-based collocation is the generalized interpolation in the deterministic domain $\Domain$ while the stochastic collocation focuses on the approximation of the finite-dimensional probability space $\Omega$.
Therefore, the kernel-based collocation methods are renamed the kernel-based methods in this article.
\end{remark}
%////////////////////////////////////////////////////////////////////////////////////////////////////////////////////////

{\bf Improvements:}
For reducing the complexity of this article, we mainly investigate the simple stochastic models here.
In fact, we can improve the above theorems, models, and algorithms in the follow ways.

i). In kriging interpolation, the estimators can be also computed by the Gaussian fields with the polynomial means.
Therefore, we improve Theorem~\ref{t:Gauss-PDK-L} to construct the probability measure $\PP_K^{\mu}$ centered at a function $\mu\in\Hilbert^m(\Domain)$ such that the Gaussian random variables defined on the Sobolev spaces also have the nonzero means.
Here $\mu$ can be viewed as the initial guess of the target function $u$.

%////////////////////////////////////////////////////////////////////////////////////////////////////////////////////////
\begin{theorem}[Improvement of Theorem~\ref{t:Gauss-PDK-L}]\label{t:Gauss-PDK-L-gen}
Suppose that the function $\mu\in\Hilbert^m(\Domain)$ and the positive definite kernel $K\in\Cont^{2m,1}\left(\Domain\times\Domain\right)$ for $m>d/2$.
Let $L$ be a bounded linear functional on the Sobolev space $\Hilbert^m(\Domain)$. Then there exists a probability measure $\PP_K^{\mu}$ on the measurable space
\[
\left(\Omega_m,\Filter_m\right):=\left(\Hilbert^m(\Domain),\Borel\left(\Hilbert^m(\Domain)\right)\right),
\]
such that the normal random variable
\[
LS(\omega):=L\omega,\quad\text{for }\omega\in\Omega_m,
\]
is well-defined on the probability space $\big(\Omega_m,\Filter_m,\PP_K^{\mu}\big)$
and this random variable $LS$ has the mean $L\mu$ and the variance $L_{\vx}L_{\vy}K(\vx,\vy)$.
Moreover, the probability measure $\PP_K^{\mu}$ is independent of the bounded linear functional $L$.
\end{theorem}
%////////////////////////////////////////////////////////////////////////////////////////////////////////////////////////
\begin{proof}
The key point of the proofs is to transfer the probability measure $\PP_K$ (given in Theorem~\ref{t:Gauss-PDK-L}) to another center at $\mu$.
Notes that $\mu\in\Hilbert^m(\Domain)$; hence we have $\mu+\Hilbert^m(\Domain)=\Hilbert^m(\Domain)$ and $\mu+\Borel\left(\Hilbert^m(\Domain)\right)=\Borel\left(\Hilbert^m(\Domain)\right)$. Then the probability measure
\[
\PP_{K}^{\mu}(A):=\PP_K\left(-\mu+A\right),
\quad\text{for }A\in\Filter_m,
\]
is well-defined on the measurable space $\left(\Omega_m,\Filter_m\right)$.

Moreover, Theorem~\ref{t:Gauss-PDK-L} guarantees that $L\mu+LS=L(\mu+S)$ is a normal random variable with the mean $L\mu$ and the variance $L_{\vx}L_{\vy}K(\vx,\vy)$ on the probability space $\big(\Omega_m,\Filter_m,\PP_K\big)$.
This assures that $\PP_K^{\mu}$ transfers the mean of the normal random variable $LS$ under $\PP_K$ from $0$ to $L\mu$.
Then the proofs are completed.
\end{proof}
%////////////////////////////////////////////////////////////////////////////////////////////////////////////////////////

%////////////////////////////////////////////////////////////////////////////////////////////////////////////////////////
\begin{remark}
Let the collection $\Gaussian$ be composed of all normal random variables $LS$ given in Theorem \ref{t:Gauss-PDK-L} or~\ref{t:Gauss-PDK-L-gen}, that is, $\Gaussian:=\left\{LS:L\in\Hilbert^m(\Domain)'\right\}$ where $\Hilbert^m(\Domain)'$ is the dual space of $\Hilbert^m(\Domain)$.
Clearly $\Hilbert^m(\Domain)'$ is a Hilbert space.
Then $\Gaussian$ is a Gaussian Hilbert space and the linear isometry $L\mapsto LS$ is a Gaussian field indexed by $\Hilbert^m(\Domain)'$
(see \cite[Definition~1.18 and 1.19]{Janson1997}). In this article, we do not consider the Gaussian Hilbert spaces and the Gaussian fields indexed by the Hilbert spaces because the theoretical formulas in \cite{Janson1997} are hard to connect to the classical kernel-based approximation.
\end{remark}
%////////////////////////////////////////////////////////////////////////////////////////////////////////////////////////

This improvement in Theorem~\ref{t:Gauss-PDK-L-gen} will give another colorful estimators to maximize the conditional probability similar as in Equation~\eqref{eq:deterministic-max-cond-prob},
that is,
\[
\max_{v\in\RR}\PP_K^{\mu}\left(\Aset_L(v)|\Aset_{\vL}(\vf)\right)
=\max_{v\in\RR}\PP_K^{\mu}\left(LS=v|\vL S=\vf\right);
\]
hence we can obtain the new kernel-based estimator
\[
L\hat{u}_{\vL}:=\underset{v\in\RR}{\text{argmax}}
~p_{L|\vL}\left(v-L\mu|\vf-\vL\mu\right)=L\mu+L\vk_{\vL}^T\vK_{\vL}^{\dag}\left(\vf-\vL\mu\right).
\]

ii). In approximation theory, the polynomials or the splines do not need to interpolate the given data exactly.
Thus, the interpolation $\Aset_{\vL}(\vf)$ (discussed in Section~\ref{sec:ker-app}) can be also improved to the oscillation $\Aset_{\vL}^{\delta}(\vf)$ for $\delta>0$, that is,
\[
\Aset_{\vL}^{\delta}(\vf):=\left\{\omega\in\Hilbert^m(\Domain): \norm{\vL\omega-\vf}_{\infty}\leq\delta\right\},
\]
and the estimate values can be measured by the sample paths oscillating around the error $\delta$ at the given data.
This indicates that Equation~\eqref{eq:deterministic-max-cond-prob} can be updated to
\[
\max_{v\in\RR}\PP_K\left(\Aset_L(v)\big|\Aset_{\vL}^{\delta}(\vf)\right)
=\max_{v\in\RR}\PP_K\left(LS=v\big|\norm{\vL S-\vf}_{\infty}\leq\delta\right);
\]
hence the best estimator $\hat{v}$ is solved by the maximum problem
\[
\hat{v}:=\underset{v\in\RR}{\text{argmax}}\frac{\int_{\norm{\vv-\vf}_{\infty}\leq\delta}p_{L,\vL}\left(v,\vv\right)\ud\vv}
{\int_{\norm{\vv-\vf}_{\infty}\leq\delta}p_{\vL}\left(\vv\right)\ud\vv}
=\underset{v\in\RR}{\text{argmax}}\int_{\norm{\vv-\vf}_{\infty}\leq\delta}p_{L,\vL}\left(v,\vv\right)\ud\vv.
\]

iii). In Sections~\ref{sec:Stochastic-Data-Approx}-\ref{sec:par-SPDE},
we only review the simple stochastic models for the comparisons of the deterministic models in \cite{Wendland2005,Fasshauer2007,HonSchabackZhong2013}.
Actually, the kernel-based methods can be applied to another complex stochastic models in the same ways. For example, we can generalize the derivative operator in Equation~\eqref{eq:L-derivative} to another differential operators
\[
L:=\sum_{\abs{\valpha}< m-d/2}\delta_{\vx}\circ a_{\valpha}D^{\valpha}\text{ or }
L:=\sum_{\abs{\valpha}\leq m}\int_{\Domain}a_{\valpha}D^{\valpha},
\quad\text{for }a_{\valpha}\in\Cont(\Domain),
\]
and the differential and boundary operators $P$ and $B$ of the SPDE~\eqref{eq:ellptic-SPDE} can be replaced by the high-order operators such as
\[
P:=\sum_{\abs{\valpha}<m-d/2}a_{\valpha}D^{\valpha}\text{ and }
B:=\sum_{\abs{\valpha}<m-d/2}b_{\valpha}D^{\valpha}|_{\partial\Domain},\quad
\text{for }a_{\valpha}\in\Cont^{\infty}(\Domain),~b_{\valpha}\in\Cont^{\infty}(\partial\Domain).
\]

\begin{wrapfigure}{r}{0.4\textwidth}
\vspace{-26pt}
\center
\includegraphics[width=0.4\textwidth, height=0.38\textwidth]{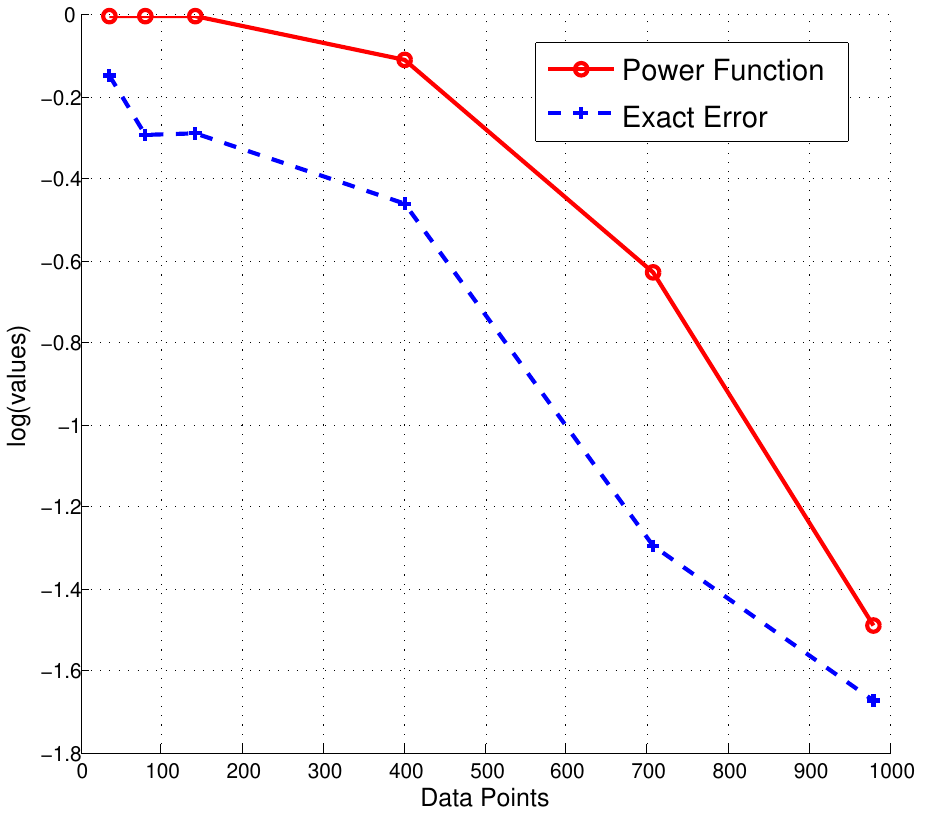}
\vspace{-30pt}
\end{wrapfigure}
iv). In Sections~\ref{sec:Stochastic-Data-Approx}-\ref{sec:par-SPDE}, we have already known that the convergence of the kernel-based estimators can be analyzed by the power functions. Moreover, according to the theorems in \cite{Wendland2005,Fasshauer2007,HonSchabackZhong2013}, we can compute the upper bounds of the power functions by the fill distances so that the convergent rates of the kernel-based estimators can be obtained by the fill distances. This technique is similar as the error estimates of the finite difference and finite element methods. Currently, the people have a great interest in the investigation of the convergence by the computational experiments only. For example in the right-hand-side figure, we can compare the power functions and the exact errors by the different data points.
In this numerical experiment, the kernel-based estimators are constructed by the Gaussian kernel with the shape parameter $\theta:=16$ and the interpolating data are evaluated by the 2D Franke's function at the Halton points.
We find that the exactly convergent rates of the kernel-based estimators follow the changes of the power functions.
This implies that the computer could learn the errors of the kernel-based estimators intelligently without the proofs by hands.
Thus, the kernel-based probability structures of the Sobolev spaces may provide another way of numerical analysis.

{\bf Advance Researches:}
The recent research of the SPDEs is still the active area including theoretical analysis and numerical algorithms.
We will continue to investigate the advance topics in our next works.

\begin{itemize}

\item For simplifying the proofs, we study with the strong conditions of the kernel-based methods such as the regularity of the domains and the smoothness of the positive definite kernels. Then we can directly apply the Sobolev imbedding theorem, the Mercer's theorem, the Kolmogorov-\v{C}entsov continuity theorem, and so on. But, the numerical examples given in Sections~\ref{sec:num-3D-model} and~\ref{sec:num-poisson-SPDE} shows that the kernel-based estimators or the kernel-based approximate solutions are still well-posed for the non-regular domains or the non-smooth kernels. Therefore, the weakened conditions could be still possible for kernel-based approximation. For example, the smooth conditions may be weakened to $K\in\Hilbert^{2m}(\Domain\times\Domain)$ because $L_{\vx}L_{\vy}K(\vx,\vy)$ is well-posed for any bounded linear functional $L$ on $\Hilbert^m(\Domain)$.

\item The kriging interpolation is a typical tool of statistical learning and the kriging predictions can be solved by the least-square loss $L(y,f(\vx)):=(y-f(\vx))^2$ for the linear models.
    Thus, we only discuss the stochastic linear models in this article.
    In fact, the kernel-based methods achieve a great success in statistical learning for the nonlinear models, for example, the minimum risks of the hinge loss $L(y,f(\vx)):=\max\left\{0,1-yf(\vx)\right\}$. In learning theory, the papers \cite{CuckerSmale2002,SmaleZhou2004} show the convergence of various loss functions for the spatial data.
    By the theorems in this article, the differential and integral data could be a new topic of statistical learning, for example, $L(y,D^{\valpha}f(\vx)):=(y-D^{\valpha}f(\vx))^2$.
    In our current researches, we also investigate the learning methods of the reproducing kernel Banach spaces induced by the positive definite kernels in \cite{FasshauerHickernellYe2013,Ye2014RKBS}. So, we will try to generalize the theorems and algorithms of the kernel-based methods to the Sobolev Banach spaces and the nonlinear stochastic models.

\item It is well known that there are still many time schemes for the SPDEs in \cite{KloedenPlaten1992,JentzenKloeden2011}.
Combing with various kinds of time schemes, we will design another kernel-based algorithms to solve the SPDEs.
Moreover, Algorithm (A3) for the white noise can be extended to the L\'evy noises in \cite{PeszatZabczy2007} such as the time and space Poisson noises.

\end{itemize}

{\bf Monographs:}
Finally, we recommend some nice books to learn the associated fields of the kernel-based methods and the SPDEs as follows:

$\bullet$ Meshfree approximation and radial basis functions: \cite{Buhmann2003,Wendland2005,Fasshauer2007}

$\bullet$ Statistical leaning and kriging interpolation: \cite{Stein1999,BerlinetThomas2004,HastieTibshiraniFriedman2009,SteinwartChristmann2008,Wahba1990}

$\bullet$ Stochastic analysis and probability theory: \cite{KaratzasShreve1991,Shiryaev1996,Janson1997}

$\bullet$ Stochastic differential equations and their numerical solutions: \cite{Oksendal2003,Chow2007,KloedenPlaten1992,JentzenKloeden2011}

{\bf Postscripts of the author:}
My researches mainly focus on approximation theory and meshfree approximation.
Now I join work with another research groups for statistical (machine) learning.
I find that the both fields are strongly connected for the kernel-based algorithms in the review papers~\cite{ScheuererSchabackSchlather2012,SchabackWendland2006}.
This inspires me to rethink the approximation theory for the stochastic data.
Moreover, the additional knowledge of stochastic analysis let me try to combine the both fields into one approach.
Just like the philosophical thoughts in Buddhism, I think that everything is correlated such as meshfree approximation and kriging interpolation discussed here.
This article may not be the perfect one to present the full connection of approximation theory and statistical learning.
But, I am sure that it is not the last one and this is just the beginning.

%---------------------------------------------------------------------------------------------------------------------
\section*{Acknowledgments}
%---------------------------------------------------------------------------------------------------------------------

The author would like to express his gratitude to Prof. Igor Cialenco and my advisor, Prof. Gregory E. Fasshauer, for their guide and assistance of this research topic at Illinois Institute of Technique, Chicago.




\bibliographystyle{elsarticle-num}
\bibliography{KernelSPDE}


%\begin{thebibliography}{00}
%%% Text of bibliographic item
%
%\bibitem{}
%
%\end{thebibliography}

\end{document}